\documentclass[invmat,numbook,envcountsame,francais]{svjour}
\usepackage{makeidx}
\usepackage{amsmath}
\usepackage{amssymb}
\usepackage{amscd}
\usepackage[all]{xypic}
\renewcommand{\hat}{\widehat}
\renewcommand{\tilde}{\widetilde}
\renewcommand{\geq}{\geqslant}
\renewcommand{\leq}{\leqslant} 
\renewcommand{\phi}{\varphi}
\renewcommand{\projlim}{\varprojlim}
\newcommand{\Qp}{\mathbf{Q}_p}
\newcommand{\Zp}{\mathbf{Z}_p}
\newcommand{\Cp}{\mathbf{C}_p}
\newcommand{\eps}{\varepsilon}
\newcommand{\ra}{\rightarrow}

\newcommand{\pibar}{\overline{\pi}}
\newcommand{\on}{\operatorname}
\newcommand{\OO}{\mathcal{O}}
\newcommand{\MM}{\mathfrak{m}}
\newcommand{\pscal}[1]{\langle #1 \rangle}
\newcommand{\ZZ}{\mathbf{Z}}
\newcommand{\NN}{\mathbf{N}}
\newcommand{\RR}{\mathbf{R}}
\newcommand{\logpi}{\log[\overline{\pi}]}
\newcommand{\btst}[2]{\widetilde{\mathbf{B}}^{\dagger #1}_{\mathrm{log} #2}}
\newcommand{\btrig}[2]{\widetilde{\mathbf{B}}^{\dagger #1}_{\mathrm{rig} #2}}
\newcommand{\bnst}[2]{\mathbf{B}^{\dagger #1}_{\mathrm{log} #2}}
\newcommand{\bnrig}[2]{\mathbf{B}^{\dagger #1}_{\mathrm{rig} #2}}
\newcommand{\btstplus}[1]{\widetilde{\mathbf{B}}^{+}_{\mathrm{log} #1}}
\newcommand{\btrigplus}[1]{\widetilde{\mathbf{B}}^{+}_{\mathrm{rig} #1}}

\newcommand{\atrig}[2]{\widetilde{\mathbf{A}}^{\dagger #1}_{\mathrm{rig} #2}}

\newcommand{\dnst}[1]{\mathbf{D}^{\dagger #1}_{\mathrm{log}}}
\newcommand{\dnrig}[1]{\mathbf{D}^{\dagger #1}_{\mathrm{rig}}}
\newcommand{\bmax}{\mathbf{B}_{\mathrm{max}}}
\newcommand{\bcont}{\mathbf{B}_{\mathrm{cont}}}
\newcommand{\bst}{\mathbf{B}_{\mathrm{st}}}
\newcommand{\bdR}{\mathbf{B}_{\mathrm{dR}}}
\newcommand{\btdag}[1]{\widetilde{\mathbf{B}}^{\dagger #1}}
\newcommand{\bdag}[1]{\mathbf{B}^{\dagger #1}}
\newcommand{\bplus}{\mathbf{B}^+}
\newcommand{\btplus}{\widetilde{\mathbf{B}}^+}
\newcommand{\bt}{\widetilde{\mathbf{B}}}
\newcommand{\amax}{\mathbf{A}_{\mathrm{max}}}

\newcommand{\atdag}[1]{\widetilde{\mathbf{A}}^{\dagger #1}}
\newcommand{\adag}[1]{\mathbf{A}^{\dagger #1}}
\newcommand{\aplus}{\mathbf{A}^+}
\newcommand{\atplus}{\widetilde{\mathbf{A}}^+}
\newcommand{\at}{\widetilde{\mathbf{A}}}
\newcommand{\bhol}[1]{\mathbf{B}^+_{\mathrm{rig} #1}}
\newcommand{\blog}[1]{\mathbf{B}^+_{\mathrm{log} #1}}
\newcommand{\e}{\mathbf{E}}
\newcommand{\eplus}{\mathbf{E}^+}
\newcommand{\etplus}{\widetilde{\mathbf{E}}^+}
\newcommand{\et}{\widetilde{\mathbf{E}}}
\newcommand{\dst}{\mathbf{D}_{\mathrm{st}}}
\newcommand{\dcris}{\mathbf{D}_{\mathrm{cris}}}
\newcommand{\ddR}{\mathbf{D}_{\mathrm{dR}}}
\newcommand{\dfont}{\mathbf{D}}

\newcommand{\ddif}{\mathbf{D}_{\mathrm{dif}}}
\renewcommand{\ddag}[1]{\mathbf{D}^{\dagger #1}}
\newcommand{\dsen}{\mathbf{D}_{\mathrm{Sen}}}
\newcommand{\mbf}{M}

\newcommand{\vale}{v_\mathbf{E}}
\newcommand{\ndr}{\mathbf{N}_{\mathrm{dR}}}
\makeindex
\begin{document}
\title{Repr{\'e}sentations $p$-adiques et {\'e}quations diff{\'e}rentielles}
\author{Laurent Berger\inst{1}}
\institute{ 
\email{laurent@brandeis.edu} \\ 
MS 050 Brandeis University \\  
PO Box 549110 \\ 
Waltham MA 02454-9110 \\ USA}
\date{\today}
\maketitle
\renewcommand{\abstractname}{Abstract}
\begin{abstract}
In this paper, we associate to every 
$p$-adic representation $V$ a $p$-adic 
differential equation 
$\mathbf{D}^{\dagger}_{\mathrm{rig}}(V)$,
that is to
say a module with a connection over the Robba ring. 
We do this via 
the theory of Fontaine's  
$(\varphi,\Gamma_K)$-modules.

This construction
enables us to relate the theory of $(\varphi,\Gamma_K)$-modules to 
$p$-adic Hodge theory. We explain how to 
construct $\mathbf{D}_{\mathrm{cris}}(V)$ 
and $\mathbf{D}_{\mathrm{st}}(V)$ 
from $\mathbf{D}^{\dagger}_{\mathrm{rig}}(V)$,
which allows us to
recognize
semi-stable or crystalline representations; 
the connection is then either unipotent or trivial
on $\mathbf{D}^{\dagger}_{\mathrm{rig}}(V)[1/t]$. 

In general, the connection has an infinite number of regular
singularities, but we show that $V$ is de Rham if and only if those
are apparent singularities. A structure theorem for modules over the
Robba ring allows us to get rid of all singularities at once, and to
obtain a ``classical'' differential equation, with a Frobenius
structure. A recent theorem of Y. Andr\'e gives a complete 
description of the structure of such an object. 
This allows us to prove Fontaine's $p$-adic monodromy
conjecture: every de Rham representation is potentially semi-stable.

As an application, we can extend to the 
case of arbitrary perfect residue fields 
some results of Hyodo ($H^1_g=H^1_{st}$), of 
Perrin-Riou (the semi-stability of ordinary representations), of Colmez 
(absolutely crystalline representations are of finite height), and
of Bloch and Kato (if the weights of $V$ are $\geq 2$, then
Bloch-Kato's exponential $\exp_V$ is an isomorphism).
\end{abstract}
\section*{Table des mati\`eres}
{\small
\begin{trivlist}
\item Introduction
\item {1. }Rappels et notations
\item {2. }Les anneaux $\btrig{}{}$ et $\btst{}{}$
\item {3. }Application aux repr{\'e}sentations $p$-adiques
\item {4. }Propri\'et\'es de $\bnrig{,r}{,K}$
\item {5. }Structures diff{\'e}rentielles sur les $(\varphi,\Gamma_K)$-modules 
et monodromie $p$-adique
\item {6. }Extensions de repr{\'e}sentations semi-stables
\item Diagramme des anneaux de p\'eriodes
\item Index des notations
\item Bibliographie
\end{trivlist}}

\section*{Introduction}

\subsection{G{\'e}n{\'e}ralit{\'e}s et notations}
Dans tout cet article $k$ est un corps parfait de caract{\'e}ristique
$p$, et $F$ est le corps des fractions de l'anneau des vecteurs de
Witt {\`a} coefficients dans $k$. Soit $K$ une extension finie  
totalement ramifi{\'e}e de $F$, 
ce qui fait que le corps r{\'e}siduel de
$\OO_K$ s'identifie {\`a} $k$ et que le corps $K$ est complet pour la
valuation $v_p$ qui {\'e}tend celle de $F$. Soit $\Cp$ le compl{\'e}t{\'e}
d'une cl{\^o}ture alg{\'e}brique $\overline{K}$ de $K$.
On s'int{\'e}ressera aux repr{\'e}sentations
$p$-adiques $V$ du groupe de Galois absolu
$G_K=\on{Gal}(\overline{K}/K)$ c'est-{\`a}-dire aux
$\Qp$-espaces vectoriels de dimension finie $d$ munis d'une action
lin{\'e}aire et continue de $G_K$.

J-M. Fontaine a construit dans \cite{Fo91}
une {\'e}quivalence de cat{\'e}gories 
$V \mapsto \dfont(V)$ entre la cat{\'e}gorie
de toutes les repr{\'e}sentations de $G_K$ 
et la cat{\'e}gorie des $(\phi,\Gamma_K)$-modules {\'e}tales (le foncteur
inverse est not{\'e} $D \mapsto \mathbf{V}(D)$).
Un $(\phi,\Gamma_K)$-module
est un espace vectoriel de dimension finie sur un corps local
$\mathbf{B}_K$ de dimension $2$ muni d'une action semi-lin{\'e}aire d'un
Frobenius $\phi$ et d'une action semi-lin{\'e}aire de $\Gamma_K$
commutant {\`a} celle de $\phi$. Un tel module est {\'e}tale si $\phi$ est
de pente $0$ (``unit root''). Le corps $\mathbf{B}_K$ est isomorphe
(non canoniquement) {\`a} l'anneau des s{\'e}ries
$\sum_{k \in \ZZ} a_k T^k$ 
en l'ind{\'e}termin{\'e}e $T$ 
o{\`u} la suite $a_k$ est une suite born{\'e}e
d'{\'e}l{\'e}ments de $F$ et $\lim_{k \ra + \infty} a_{-k}=0$. L'action de
$\phi$ et de $\Gamma_K$ est assez compliqu{\'e}e en g{\'e}n{\'e}ral mais si
$K$ est non-ramifi{\'e} ($K=F$), alors on peut choisir $T$ de telle sorte 
que $\phi(T)=(1+T)^p-1$ et $\gamma(T)=(1+T)^{\chi(\gamma)}-1$
($\chi:G_K \ra \Zp^*$ est le caract\`ere cyclotomique;
dans le corps du texte, $T$ est {\'e}gal {\`a} un {\'e}l{\'e}ment $\pi_K$
construit via la th{\'e}orie du corps des normes et
l'action de $\phi$ et de $\Gamma_K$ provient d'une action naturelle;
si $K=F$, on a $\pi_K=\pi=[\eps]-1$).

L'anneau $\mathbf{B}_K$ est assez d{\'e}sagr{\'e}able mais il contient
l'anneau $\bdag{}_K$ des s{\'e}ries surconvergentes, c'est-{\`a}-dire
l'anneau des s{\'e}ries
$A(T)=\sum_{k \in \ZZ} a_k T^k$ o{\`u} la suite $a_k$ est une suite born{\'e}e
d'{\'e}l{\'e}\-ments de $F$ et il existe $r<1$ tel que la s{\'e}rie  
$A(T)$ converge sur une couronne non-vide 
du type $\{ z \in \Cp, r < |z| < 1 \}$.
Le th{\'e}or{\`e}me principal de \cite{CC98} est que tout
$(\phi,\Gamma_K)$-module {\'e}tale provient par extension des scalaires
d'un ``module surconvergent''; plus pr{\'e}cis{\'e}ment on a le r{\'e}sultat 
suivant 
\begin{theorem}
Si $D$ est un $(\phi,\Gamma_K)$-module {\'e}tale, alors 
l'ensemble des sous-$\bdag{}_K$-modules libres de type fini 
stables par $\phi$ et  $\Gamma_K$ admet un plus grand {\'e}l{\'e}ment 
$D^{\dagger}$ et on a 
$D=\mathbf{B}_K 
\otimes_{\bdag{}_K} D^{\dagger}$.
\end{theorem}
Le fait que l'ensemble des sous-$\bdag{}_K$-modules libres de type fini 
stables par $\phi$ et  $\Gamma_K$ admet un plus grand {\'e}l{\'e}ment est un
r{\'e}sultat de Cherbonnier \cite{Ch96}, et moyennant ce r{\'e}sultat, 
le th{\'e}or{\`e}me est {\'e}nonc{\'e} dans \cite{CC98}
de mani{\`e}re duale (via le foncteur $D \mapsto \mathbf{V}(D)$)
sous la forme ``toute repr{\'e}sentation de $G_K$ est surconvergente''.
Gr{\^a}ce {\`a} ce r{\'e}sultat, on va pouvoir tensoriser au-dessus de
$\bdag{}_K$ avec l'anneau de Robba $\bnrig{}{,K}$
(apparaissant dans la th{\'e}orie
des {\'e}quations diff{\'e}rentielles $p$-adiques),
constitu{\'e} des s{\'e}ries 
$A(T)=\sum_{k \in \ZZ} a_k T^k$ telles que $a_k \in F$ et
qu'il existe $r<1$ tel que la s{\'e}rie  $A(T)$
converge sur une couronne non-vide du type
$\{ z \in \Cp, r < |z| < 1 \}$ (sans condition de croissance au
voisinage de $\{ z \in \Cp, |z| = 1 \}$). 
Ceci va nous permettre, si $V$ est une repr{\'e}sentation $p$-adique de
$G_K$, d'une part de retrouver les invariants $\dcris(V)$ et $\dst(V)$
associ{\'e}s {\`a} $V$ par la th{\'e}orie de Hodge $p$-adique et, d'autre
part, en utilisant l'action infinit{\'e}simale de $\Gamma_K$,
d'associer {\`a} $V$ un module {\`a} connexion $\dnrig{}(V)$ sur
l'anneau de Robba $\bnrig{}{,K}$.
Cet article 
rassemble un certain nombre de r{\'e}sultats que l'on peut obtenir
via l'{\'e}tude de ce module {\`a} connexion: caract{\'e}risation des
repr{\'e}sentations absolument cristallines, 
semi-stables, de de Rham et $\Cp$-admissibles et d\'emonstration de 
la conjecture de ``monodromie $p$-adique''.

\subsection{Th{\'e}orie de Hodge
$p$-adique et $(\phi,\Gamma_K)$-modules} 
Soient $\bnst{}{,K}=\bnrig{}{,K}[\log(T)]$ (anneau muni des actions
{\'e}videntes de $\phi$ et de $\Gamma_K$), 
\[ \dnrig{}(V)=\bnrig{}{,K} \otimes_{\bdag{}_K} \ddag{}(V)
\ \ \text{et} \ \ \dnst{}(V)=\bnst{}{,K} \otimes_{\bdag{}_K} \ddag{}(V). \]
Si $V$ est une repr{\'e}sentation $p$-adique de
$G_K$, soient $\dcris(V)$, $\dst(V)$ et $\ddR(V)$ les modules 
associ{\'e}s {\`a} $V$ par la th{\'e}orie de Hodge $p$-adique.
Le premier r{\'e}sultat
est que l'on peut retrouver ces modules 
\begin{theorem}\label{introcomp}
Si $V$ est une repr{\'e}sentation $p$-adique de $G_K$, alors 
\[ \dst(V)=(\dnst{}(V)[1/t])^{\Gamma_K} \ \text{ et }\  
\dcris(V)=(\dnrig{}(V)[1/t])^{\Gamma_K}. \]
De plus:
\begin{enumerate}
\item si $V$ est une repr{\'e}sentation semi-stable, alors \[ \ddag{}(V) 
\otimes_{\bdag{}_K}
\bnst{}{,K}[1/t] = \dst(V) \otimes_F \bnst{}{,K}[1/t] \]
\item si $V$ est une repr{\'e}sentation cristalline, alors \[ \ddag{}(V) 
\otimes_{\bdag{}_K}
\bnrig{}{,K}[1/t] = \dcris(V) \otimes_F \bnrig{}{,K}[1/t] \]
\end{enumerate}
\end{theorem}
Dans l'{\'e}nonc{\'e} ci-dessus le $t$ qui appara{\^\i}t est un
{\'e}l{\'e}ment de $\bnrig{}{,K}$ sur lequel $\Gamma_K$ agit via
$\gamma(t)=\chi(\gamma)t$ et tel que $\phi(t)=pt$; 
si $K=F$, on a $t=\log(1+T)$.

Les modules 
$\dcris(V)$ et $\dst(V)$ 
sont naturellement munis d'un Frobenius $\phi$, d'un op{\'e}rateur de
monodromie $N$ et d'une filtration. Le th{\'e}or{\`e}me pr{\'e}c{\'e}dent 
permet de retrouver les actions de $\phi$ et $N$, mais la recette
permettant de retrouver la filtration est suffisament peu rago{\^u}tante
pour ne pas {\^e}tre explicit{\'e}e dans cet article.

Un corollaire du th{\'e}or{\`e}me \ref{introcomp} est le r{\'e}sultat suivant: 
\begin{corollary}\label{introhf}
Si $V$ est une repr{\'e}sentation cristalline de $G_F$, alors $V$ est de
hauteur finie. 
\end{corollary}
Ce r{\'e}sultat avait {\'e}t{\'e} conjectur{\'e} par Fontaine (voir
\cite{Fo91,Wa96}) et d{\'e}montr{\'e} par Colmez de mani{\`e}re tr{\`e}s 
d{\'e}tourn{\'e}e (la d{\'e}monstration utilisait deux versions
\cite{Co98,CC99} de la 
d{\'e}monstration de la loi de r{\'e}ciprocit{\'e} conjectur{\'e}e par
Perrin-Riou \cite{BP94}) dans le cas d'un corps r{\'e}siduel fini.
La d{\'e}monstra\-tion donn{\'e}e dans cet article utilise {\`a} la place 
le th{\'e}or{\`e}me ci-dessus et un r{\'e}sultat d'analyse $p$-adique 
d{\'e}montr{\'e} par Kedlaya \cite{KK00}. Signalons que Benois \cite{Be00}
a par ailleurs d{\'e}montr{\'e} la loi de r{\'e}ciprocit{\'e} de Perrin-Riou
pour les repr{\'e}sentations cristallines de hauteur finie ce qui,
combin{\'e} avec le r{\'e}sultat ci-dessus, fournit une d{\'e}monstration de
cette loi de r{\'e}ciprocit{\'e} dans le cas g{\'e}n{\'e}ral, n'utilisant que
la th{\'e}orie des $(\phi,\Gamma_K)$-modules, r{\'e}alisant ainsi un
programme de Fontaine.

\subsection{Structures diff{\'e}rentielles sur les $(\phi,\Gamma_K)$-modules et monodromie
$p$-adique} 

Soit $\gamma$ un {\'e}l{\'e}ment de $\Gamma_K$ assez proche de $1$. 
La s{\'e}rie qui d{\'e}finit $\frac{\log(\gamma)}{\log(\chi(\gamma))}$ 
converge vers un op{\'e}rateur diff{\'e}rentiel $\nabla$; si $K=F$ 
on a $\nabla=\log(1+T)(1+T)\frac{d}{dT}$.
Si $\bnrig{,r_n}{,K}$ est l'anneau des s{\'e}ries convergeant sur
la couronne $\{ p^{-1/e_K r_n} \leq |z| < 1 \}$
(o{\`u} $e_K=[K_{\infty}:F_{\infty}]$ et $r_n=(p-1)p^{n-1}$),
alors $\nabla$ stabilise $\bnrig{,r_n}{,K}$ pour $n$ assez grand et on
dispose de plus de morphismes injectifs $\iota_n: \bnrig{,r_n}{,K} \ra K_n
[[t]]$ qui v{\'e}rifient $\iota_n \circ \nabla= t \frac{d}{dt}\circ \iota_n$.
Soit aussi $\partial=t^{-1}\nabla$; si $K=F$ 
on a donc $\partial=(1+T)\frac{d}{dT}$, c'est une autre d{\'e}rivation
de $\bnrig{}{,K}$, qui est d'ailleurs une base du $\bnrig{}{,K}$-module
de rang $1$ des d{\'e}rivations de $\bnrig{}{,K}$.

Ce qui pr{\'e}c{\`e}de correspond au cas de la repr{\'e}sentation triviale
et, plus g{\'e}n{\'e}ralement, si $V$ est 
une repr{\'e}sentation $p$-adique de $G_K$ et si 
$\gamma$ est un {\'e}l{\'e}ment de $\Gamma_K$ assez proche de $1$,
la s{\'e}rie qui d{\'e}finit $\log(\gamma)/\log(\chi(\gamma))$ 
converge vers un op{\'e}rateur diff{\'e}rentiel
$\nabla_V$ de $\dnrig{}(V)$ 
au-dessus de $\bnrig{}{,K}$ muni de 
$\nabla$. Plus pr{\'e}cis{\'e}ment, si $n$ est assez grand, $\nabla_V$
stabilise le sous-module  $\dnrig{,r_n}(V)$ ce qui permet d'associer
{\`a} toute  repr{\'e}sentation $p$-adique de $G_K$ un module diff{\'e}rentiel
sur une couronne.

Une cons{\'e}quence presque imm{\'e}diate de \ref{introcomp} est le th\'eor\`eme suivant:
\begin{theorem}\label{introconn}
Si $V$ est une repr{\'e}sentation $p$-adique, alors il existe $n$ tel que
la restriction de $V$ {\`a} $G_{K(\mu _{p^n})}$ est
semi-stable (respectivement cristalline)
si et seulement si $\nabla_V$ est une connexion unipotente
(respectivement triviale) sur $\dnrig{}(V)[1/t]$.
\end{theorem}

Le module diff{\'e}rentiel $\dnrig{}(V)$ 
est {\`a} points singuliers r{\'e}guliers, mais a
une infinit{\'e} de singularit{\'e}s dans la couronne; en effet la s{\'e}rie
$\log(1+T)$ s'annule en tous les $\zeta-1$, o{\`u} $\zeta$ est une
racine d'ordre $p^n$ de $1$.
On n'est donc pas dans le cadre classique des {\'e}quations
diff{\'e}rentielles $p$-adiques, mais si la repr{\'e}sentation $V$ est de
de Rham, on peut supprimer les singularit{\'e}s et retomber sur 
un objet classique (voir le th{\'e}or{\`e}me \ref{ndrintro} ci-dessous).

En localisant en $\zeta-1$ (o{\`u} $\zeta$ est un racine primitive
d'ordre $p^n$ de l'unit{\'e}), ce qui revient {\`a} consid{\'e}rer
l'application
$\iota_n:  \dnrig{,r_n}(V) \ra (\bdR^+ \otimes_{\Qp} V)^{H_K}$,
on retombe sur le module diff{\'e}rentiel consid{\'e}r{\'e} par Fontaine
dans \cite[chap. 3]{Fo00}, ce qui permet en particulier de retrouver
le module $\ddR(V)$ {\`a} partir du noyau de cette connexion, et de
montrer que la repr{\'e}sentation $V$ est de de Rham si et seulement si
$\nabla_V$ n'a que des singularit{\'e}s apparentes en les $\zeta-1$.
Il est donc naturel d'esp{\'e}rer que l'on peut supprimer toutes les
singularit{\'e}s simultan{\'e}ment ou, ce qui revient au m{\^e}me,
construire un sous-objet de rang maximal, 
stable par $\partial_V=t^{-1}\nabla_V$.
De fait, on a le r{\'e}sultat suivant:
 
\begin{theorem}\label{ndrintro}
Soit $V$ une repr{\'e}sentation de de Rham de $G_K$, dont
les poids de Hodge-Tate sont n{\'e}gatifs.
Alors il existe un unique sous
$\bnrig{}{,K}$-module $\ndr(V)$ de $\dnrig{}(V)$, libre de rang $d$ et
stable par $\partial_V$. 
De plus, $\ndr(V)$ est stable par $\phi$ et $\Gamma_K$.
\end{theorem}

Le $\bnrig{}{,K}$-module $\ndr(V)$ est donc muni d'une connexion et
d'un op{\'e}rateur de Frobenius, et Andr\'e a r\'ecemment montr\'e
qu'un tel module est quasi-unipotent.
Une cons{\'e}quence de ceci, alli{\'e} au th{\'e}or{\`e}me \ref{introconn},
est une d\'emonstration de la conjecture
de monodromie pour les repr{\'e}sentations $p$-adiques \cite[6.2]{Bu88sst}: 

\begin{theorem}\label{intromono}
Si $V$ est une repr{\'e}sentation de de Rham, 
alors $V$ est potentiellement semi-stable.
\end{theorem}

Enfin, en utilisant l'application $\theta: \bdR^+ \ra \Cp$, on retombe 
sur le module de Sen (le $K_{\infty}$-espace vectoriel associ{\'e} {\`a}
$V$ via la th{\'e}orie de Sen, voir \cite{Sn73,Co94,Fo00}),
la connexion devenant un op{\'e}rateur
$K_n$-lin{\'e}aire dont les valeurs propres sont les ``poids de
Hodge-Tate g{\'e}n{\'e}ralis{\'e}s''.
En particulier, $V$ est
$\Cp$-admissible 
(ce qui {\'e}quivaut {\`a} $V$ de Hodge-Tate et tous ses poids de
Hodge-Tate sont nuls)
si et seulement si 
cet op{\'e}rateur est nul. 
On en d{\'e}duit le fait que $\nabla_V$ est divisible par $t$, et donc que
$\ndr(V)=\dnrig{}(V)$ est un $(\phi,\partial)$-isocristal
surconvergent, avec un Frobenius {\'e}tale.
Utilisant ce fait et des techniques
d'{\'e}quations diff{\'e}rentielles $p$-adiques (plus
pr{\'e}cis{\'e}ment un th{\'e}or{\`e}me de Tsuzuki \cite{GC00,TS99})
on obtient une d{\'e}monstration du r{\'e}sultat suivant, d{\^u} {\`a} Sen 
\cite{Sn73}, qui
caract{\'e}rise les repr{\'e}sentations $\Cp$-admissibles:
\begin{theorem}\label{introsen}
Si $V$ est une repr{\'e}sentation $p$-adique de $G_K$, alors les deux conditions
suivantes sont {\'e}quivalentes:
\begin{enumerate}
\item $V$ est $\Cp$-admissible;
\item le sous-groupe d'inertie de $G_K$ agit sur $V$ {\`a} travers un
quotient fini.
\end{enumerate}
\end{theorem}
Ce th\'eor\`eme est d'ailleurs \'equivalent au th\'eor\`eme de Tsuzuki, 
voir la remarque \ref{remasentsu}.

\subsection{Extensions de repr{\'e}sentations semi-stables}
En utilisant le th\'eor\`eme de monodromie $p$-adique, on peut montrer que:
\begin{theorem}\label{introordi}
\begin{enumerate}
\item Toute repr{\'e}sentation ordinaire de $G_K$ est semi-stable;
\item une repr\'esentation de de Rham, qui est une
extension de deux repr{\'e}sentations semi-stables,  
est elle-m{\^e}me semi-stable;
\item si $V$ est une repr\'esentation semi-stable dont les 
poids de Hodge-Tate sont tous $\geq 2$, alors 
$\exp_V : \ddR(V) \ra H^1(K,V)$ est un isomorphisme.
\end{enumerate}
\end{theorem}
Ces trois r{\'e}sultats {\'e}taient connus dans le cas d'un corps
r{\'e}siduel fini, o{\`u} on peut les d{\'e}duire de calculs de dimensions de
groupes de cohomologie galoisienne (ce qui n'est plus possible si le
corps r{\'e}siduel n'est pas fini).
Le (1) avait {\'e}t{\'e} d{\'e}montr{\'e} dans ce cas l{\`a} par
Perrin-Riou \cite{Bu88ord,BP94,BP99}  
comme corollaire 
des calculs de Bloch et Kato 
\cite{BK91}, le (2) par Hyodo
\cite{Hy95,Ne95}, et le (3) par
Bloch et Kato.

\subsection{Plan de l'article}
Cet article comporte six chapitres  subdivis{\'e}s en paragraphes. 
Le premier chapitre est consacr{\'e} {\`a} des
rappels sur les repr{\'e}sentations $p$-adiques et les anneaux de
Fontaine. Dans le deuxi{\`e}me on donne la construction des anneaux
$\btrig{}{}$ et $\btst{}{}$, qui sont fondamentaux pour ce qui suit. On
donne dans le troisi{\`e}me chapitre
une application de ces constructions: comment retrouver $\dcris(V)$ et
$\dst(V)$ {\`a} partir du $(\phi,\Gamma_K)$-module $\dfont(V)$. Le
quatri{\`e}me chapitre est consacr{\'e} {\`a} l'{\'e}tude de la connexion
$\nabla$ sur l'anneau $\bnrig{}{,K}$, ce qui permet de d{\'e}finir dans
le cinqui{\`e}me chapitre l'{\'e}quation diff{\'e}rentielle associ{\'e}e {\`a}
une repr{\'e}sentation $p$-adique; on donne des applications {\`a} la
th{\'e}orie de Fontaine, {\`a} la th{\'e}orie de Sen, et aux 
repr{\'e}sentations de de Rham
(preuve de la conjecture de monodromie $p$-adique). Enfin dans le
sixi{\`e}me chapitre on donne quelques applications du th\'eor\`eme 
de monodromie $p$-adique.

Le th{\'e}or{\`e}me \ref{introcomp} est la r{\'e}union du th{\'e}or{\`e}me 
\ref{isomcomp} et de la proposition \ref{isomcomp2}. Le corollaire
\ref{introhf} est l'objet du th{\'e}or{\`e}me \ref{hf}, le th{\'e}or{\`e}me 
\ref{introsen} correspond {\`a} la proposition \ref{sen} et le th{\'e}or{\`e}me
\ref{introconn} {\`a} la proposition \ref{conilp}. 
Le th{\'e}or{\`e}me \ref{ndrintro} et sa r{\'e}ciproque sont d{\'e}montr{\'e}s dans le paragraphe
\ref{ndrsec}. Le th{\'e}or{\`e}me \ref{intromono} correspond au th\'eor\`eme
\ref{monofont}.
Le premier point de \ref{introordi} est le corollaire
\ref{ordisst}, le (2) est le th\'eor\`eme \ref{extsst}, 
et le (3) est le th\'eor\`eme \ref{expbk}.
Cet article comporte deux appendices pour le rendre plus lisible: 
un diagramme des anneaux de p{\'e}riodes, et un index des notations.

\begin{acknowledgement}
Cet article est une version am{\'e}lior{\'e}e de ma th{\`e}se avec Pierre
Colmez. Je tiens {\`a} le remercier d'avoir partag{\'e} avec moi ses id{\'e}es
et ses connaissances, ainsi que pour le temps et l'{\'e}nergie qu'il m'a
consacr{\'e}s.

Lors de la r{\'e}daction de cet article, j'ai eu des conversations
enrichissantes avec Pierre Berthelot, Bruno Chiarellotto, Gilles Christol, 
Richard Crew, Jean-Marc Fontaine et Kiran Kedlaya. Je remercie le referee 
pour ses nombreuses suggestions, qui m'ont permis de grandement am\'eliorer 
cet article.

Je tiens aussi \`a remercier les organisateurs 
du ``Dwork Trimester'' \`a Padova pour leur
hospitalit\'e. Mon s\'ejour \`a Padova a \'et\'e 
une occasion unique de discuter du contenu de cet article. 
\end{acknowledgement}

\section{Rappels et notations}
Ce chapitre est enti{\`e}rement constitu{\'e} de rappels sur la th{\'e}orie
des repr{\'e}\-sentations $p$-adiques. On se reportera {\`a} 
\cite{CC98,CC99,Co98,Bu88per,Bu88sst,Fo91} 
ou aussi \cite{Co00} pour la d{\'e}monstration des
faits qui sont rappel{\'e}s ici. Pour ce qui est de la th{\'e}orie des
{\'e}quations diff{\'e}rentielles $p$-adiques et
des isocristaux surconvergents, le lecteur pourra se
reporter {\`a} \cite{CM00,CR95,RC98,KK00,TS97,TS99}.

La principale strat{\'e}gie 
(due {\`a} Fontaine, voir par exemple \cite{Bu88sst})
pour {\'e}tudier les
repr{\'e}sentations $p$-adiques d'un groupe $G$ est de cons\-truire des 
$\Qp$-alg{\`e}bres topologiques $B$ munies d'une action du groupe
$G$ et de structures suppl{\'e}\-mentaires de telle mani{\`e}re que si $V$
est une repr{\'e}sentation $p$-adique, 
alors $D_B(V)=(B \otimes_{\Qp} V)^G$  est un
$B^G$-module qui h{\'e}rite de ces structures, et que le foncteur
qui {\`a} $V$ associe $D_B(V)$ fournisse des invariants int{\'e}ressants
de $V$. On dit qu'une repr{\'e}sen\-tation $p$-adique $V$ de $G$ est
\emph{$B$-admissible} si on a $B\otimes_{\Qp} V \simeq B^d$ en tant que
$G$-modules. 

\subsection{Le corps $\et$ et ses sous-anneaux}
Soit $k$\index{k@$k$} un corps parfait de caract{\'e}ristique $p$,
$F$\index{F@$F$} le corps des fractions de l'anneau des vecteurs de Witt sur $k$
et $K$\index{K@$K$} une extension finie totalement ramifi{\'e}e de $F$. Soit
$\overline{F}$\index{Fbar@$\overline{F}$} une
cl{\^o}ture alg{\'e}brique de $F$ et $\Cp=\hat{\overline{F}}$\index{Cp@$\Cp$} sa compl{\'e}tion
$p$-adique. 
On pose $G_K=\on{Gal}(\overline{K}/K)$\index{GK@$G_K$}, c'est aussi le groupe
des automorphismes continus $K$-lin{\'e}aires de $\Cp$.
Le corps $\Cp$ est un corps complet alg{\'e}briquement clos
de corps r{\'e}siduel $\OO_{\Cp}/\MM_{\Cp}=\overline{k}$. 
On pose aussi $K_n=K(\mu _{p^n})$\index{Kn@$K_n$} et 
$K_{\infty}$\index{Kinfini@$K_{\infty}$} est d{\'e}fini comme 
{\'e}tant la r{\'e}union des $K_n$.  
Soit $H_K$\index{HK@$H_K$} le noyau du caract{\`e}re cyclotomique 
$\chi: G_K \ra \Zp^*$\index{chi@$\chi$} et $\Gamma_K=G_K/H_K$\index{GammaK@$\Gamma_K$}
le groupe de Galois de $K_{\infty}/K$
qui s'identifie via le caract{\`e}re cyclotomique {\`a} un sous groupe ouvert
de $\Zp^*$.
Soient\footnote{Le lecteur est invit\'e \`a consulter
l'appendice ``Diagramme des anneaux de p\'eriodes'' tout au long de la lecture de ces
constructions.}   
\[ \et=\projlim_{x\mapsto x^p} \Cp 
=\{ (x^{(0)},x^{(1)},\cdots) \mid (x^{(i+1)})^p = x^{(i)} \} \]\index{Etilde@$\et$}
et $\etplus$\index{Etildeplus@$\etplus$} l'ensemble des $x \in \et$ tels que $x^{(0)} \in \OO_{\Cp}$.
Si $x=(x^{(i)})$ et $y=(y^{(i)})$ sont deux {\'e}l{\'e}ments de $\et$,
alors on d{\'e}finit leur somme $x+y$ et leur produit $xy$ par:
\[ (x+y)^{(i)}= \lim_{j \ra \infty} (x^{(i+j)}+y^{(i+j)})^{p^j} \text{ et }
(xy)^{(i)}=x^{(i)}y^{(i)} \] ce qui fait de
$\et$ un corps de caract{\'e}ristique $p$ dont on peut montrer qu'il est
alg{\'e}briquement clos.
Si $x=(x^{(n)})_{n \geq 0} \in \et$ soit $\vale(x)=v_p(x^{(0)})$\index{ve@$\vale$}. C'est
une valuation sur $\et$ pour laquelle celui-ci est complet;
l'anneau des entiers de $\et$ est $\etplus$ et l'id{\'e}al maximal 
est $\MM_{\et}=\{ x \in \et, \vale(x)>0 \}$.
Soit $\atplus$\index{Atildeplus@$\atplus$} l'anneau $W(\etplus)$ des vecteurs de Witt {\`a}
coefficients dans $\etplus$ et \[ \btplus =\atplus[1/p] 
=\{ \sum_{k\gg -\infty} p^k [x_k],\ x_k \in \etplus
\}\]\index{Btildeplus@$\btplus$} 
o{\`u} $[x] \in \atplus$ est le
rel{\`e}vement de Teichm{\"u}ller de $x \in \etplus$\index{[]@$[\cdot]$}.
Cet anneau est muni d'un morphisme d'anneaux
$\theta: \btplus \ra \Cp$\index{theta@$\theta$} d{\'e}fini de la mani{\`e}re suivante:
\[ \theta\left(\sum_{k\geq 0} p^k[x_k]\right)
=\sum_{k\geq 0} p^k x_k^{(0)} \]
Soient $\eps=(\eps^{(i)})\in\etplus$\index{epsilon@$\eps$} avec $\eps^{(0)}=1$ 
et $\eps^{(1)}\neq 1$, $\pi=[\eps]-1$, $\pi_n=
[\eps^{p^{-n}}]-1$\index{pi@$\pi$,$\pi_n$}, $\omega=\pi/\pi_1$\index{omega@$\omega$}
et $q=\phi(\omega)=\phi(\pi)/\pi$\index{q@$q$}. Alors $\ker(\theta)$
est l'id{\'e}al principal engendr{\'e} par $\omega$.
De m{\^e}me soit $\tilde{p}=(p^{(n)})\in\etplus$\index{ptilde@$\tilde{p}$} avec $p^{(0)}=p$, alors
$\ker(\theta)$ est aussi engendr{\'e} par $[\tilde{p}]-p$.

Remarquons que $\eps$ est un {\'e}l{\'e}ment de $\etplus$ tel que
$\vale(\eps-1)=p/(p-1)$. On pose $\e_F=k((\eps-1))$\index{EF@$\e_F$} et on d{\'e}finit
$\e$\index{E@$\e$} comme {\'e}tant 
la cl{\^o}ture s{\'e}parable de $\e_F$ dans $\et$ ainsi que $\eplus=\e
\cap \etplus$\index{Eplus@$\eplus$} et $\MM_{\e}=\e \cap \MM_{\et}$ l'anneau des entiers
et l'id{\'e}al maximal de $\e$. 
Remarquons que, par d{\'e}finition, $\e$ est s{\'e}parablement clos, et que 
l'on retrouve $\et$ {\`a} partir de $\e$ en prenant le compl{\'e}t{\'e} de
sa cl{\^o}ture radicielle.
On renvoie {\`a} \cite[p. 243]{CC99} pour
une application de la th{\'e}orie du corps de normes \cite{Wi83} {\`a} 
la construction d'une application $\iota_K : \projlim \OO_{K_n} \ra
\etplus$\index{iotaK@$\iota_K$} de la limite projective des $\OO_{K_n}$ relativement aux
applications normes dans $\etplus$ dont la principale propri{\'e}t{\'e}
est la suivante:
\begin{proposition}\label{youpla}
L'application $\iota_K$ induit une bijection de $\projlim  \OO_{K_n}$
sur l'anneau des entiers $\eplus_K$ de $\e_K=\e^{H_K}$\index{EK@$\e_K$,$\eplus_K$}. 
\end{proposition}

On remarquera que la proposition ci-dessus est {\'e}nonc{\'e}e 
dans \cite[I.1.1]{CC99} avec la
restriction suppl{\'e}\-mentaire que $K/\Qp$ est finie mais ceci 
n'est pas n{\'e}cessaire (en revanche il est important que $K/F$ soit finie).

De plus on peut montrer de la m{\^e}me
mani{\`e}re que $\e_K$ est une extension
finie s{\'e}parable de $\e_F$ de degr{\'e} $e_K=[K_{\infty}:F_{\infty}]$\index{eK@$e_K$}
et de groupe de Galois $H_F/H_K$ 
(si $K/F$ est galoisienne), 
et que le groupe de Galois
$\on{Gal}(\e/\e_K)$ s'identifie {\`a} $H_K$. Enfin $\e_K^+$ est un
anneau de valuation discr{\`e}te de la forme $k[[\pibar_K]]$ o{\`u}
$\pibar_K=\iota_K(\varpi_K)$\index{pibar@$\pibar_K$} 
est l'image d'une suite $\varpi_K$ 
d'uniformisantes
compatibles pour les normes 
des $\OO_{K_n}$.

\subsection{L'anneau $\bdR$ et ses sous-anneaux}
L'anneau $\bdR^+$\index{BdRplus@$\bdR^+$} est d{\'e}fini comme {\'e}tant le compl{\'e}t{\'e} pour la 
topologie $\ker(\theta)$-adique de $\btplus$ (on remarquera que
$\atplus$ est complet pour cette topologie):
\[ \bdR^+=\projlim_{n\geq 0} \btplus/(\ker(\theta)^n) \]
c'est un anneau de valuation discr{\`e}te, d'id{\'e}al maximal $\ker(\theta)$;
la s{\'e}rie qui d{\'e}finit $\log([\eps])$ converge
dans $\bdR^+$ vers un {\'e}l{\'e}ment $t$\index{t@$t$}, qui est un g{\'e}n{\'e}rateur de l'id{\'e}al
maximal, ce qui fait que $\bdR=\bdR^+[1/t]$\index{BdR@$\bdR$} est un corps, muni d'une
action de $G_F$ et d'une filtration d{\'e}finie par 
$\on{Fil}^i(\bdR)=t^i \bdR^+$ pour $i \in \ZZ$.

On dit qu'une repr{\'e}sentation $V$ de $G_K$
est de de Rham si 
elle est $\bdR$-admissible ce qui {\'e}quivaut {\`a} ce que
le $K$-espace
vectoriel $\ddR(V)=(\bdR\otimes_{\Qp} V)^{G_K}$\index{DdR@$\ddR(V)$}
est de dimension $d=\dim_{\Qp}(V)$\index{d@$d$}.

L'anneau $\bmax^+$\index{Bmaxplus@$\bmax^+$} est d{\'e}fini comme {\'e}tant
\[ \bmax^+= \{ \sum_{n \geq 0} a_n \frac{\omega^n}{p^n}
\text{ o{\`u} $a_n\in \btplus$ est une suite qui tend vers $0$} \} \]
et $\bmax=\bmax^+[1/t]$\index{Bmax@$\bmax$}. 
On peut d'ailleurs remplacer $\omega$ par
n'importe quel g{\'e}n{\'e}rateur de $\ker(\theta)$, par exemple
$[\tilde{p}]-p$.
Cet anneau se plonge canoniquement dans $\bdR$ (les s{\'e}ries
d{\'e}finissant ses {\'e}l{\'e}ments convergent dans $\bdR$) 
et en particulier il est muni de l'action
de Galois et de la filtration induites par celles de $\bdR$, ainsi que
d'un Frobenius $\phi$\index{phi@$\phi$}, qui {\'e}tend l'application $\phi:\atplus \ra
\atplus$ d{\'e}duite de $x \mapsto x^p$ dans $\etplus$.
On remarquera que $\phi$ ne se prolonge pas par continuit{\'e} {\`a} $\bdR$.
On pose $\btrigplus{}=\cap_{n=0}^{+\infty}
\phi^n(\bmax^+)$\index{Btildeplusrig@$\btrigplus{}$}
(on remarquera que l'anneau $\btrigplus{}$ est
l'anneau $\bcont^+$ de \cite{Co98}).
La notation s'explique par le fait que l'on a \cite{Bt01} un isomorphisme:
$\btrigplus{}= H^0_{\mathrm{rig}}(\OO_{\overline{F}}/F)$.
 
On dit qu'une repr{\'e}sentation $V$ de $G_K$
est cristalline si 
elle est $\bmax$-admissible ou, ce qui revient au m{\^e}me,
$\btrigplus{}[1/t]$-admissible
(les p{\'e}riodes des repr{\'e}sentations cristallines vivent dans des
sous $F$-espaces vectoriels de dimension finie,
stables par $\phi$, 
de $\bmax$, et donc en fait
dans l'anneau $\cap_{n=0}^{+\infty}
\phi^n(\bmax^+)[1/t]$); 
ceci {\'e}quivaut {\`a} ce que
le $F$-espace
vectoriel \[ \dcris(V)=(\bmax \otimes_{\Qp} V)^{G_K}=(\btrigplus{}[1/t]
\otimes_{\Qp} V)^{G_K} \]\index{Dcris@$\dcris(V)$} 
est de dimension 
$d=\dim_{\Qp}(V)$. Alors $\dcris(V)$ est muni d'un Frobenius et d'une 
filtration induits par ceux de $\bmax$, 
et $(\bdR\otimes_{\Qp} V)^{G_K}=\ddR(V)=
K \otimes_F \dcris(V)$ ce qui fait qu'une repr{\'e}sentation cristalline
est aussi de de Rham.

La s{\'e}rie qui d{\'e}finit $\log(\overline{\pi}^{(0)})+
\log([\overline{\pi}]/\overline{\pi}^{(0)})$ 
(apr{\`e}s avoir choisi $\log(p)$ et o{\`u} $\overline{\pi}=\eps-1$)
converge dans $\bdR^+$ vers
un {\'e}l{\'e}ment $\logpi$\index{logpi@$\logpi$} qui est 
transcendant sur $\on{Frac} \bmax^+$
et on pose
$\bst=\bmax[\logpi]$\index{Bst@$\bst$}
et $\btstplus{}=\btrigplus{}[\logpi]$\index{Btildepluslog@$\btstplus{}$}. 
On dit qu'une repr{\'e}sentation $V$
est semi-stable si 
elle est $\bst$-admissible ou, ce qui revient au m{\^e}me,
$\btstplus{}[1/t]$-admissible
(m{\^e}me raison que ci-dessus); 
ceci {\'e}quivaut {\`a} ce que
le $F$-espace
vectoriel \[ \dst(V)=(\bst \otimes_{\Qp} V)^{G_K} 
=  (\btstplus{}[1/t] \otimes_{\Qp} V)^{G_K}  \]\index{Dst@$\dst(V)$} 
est de dimension 
$d=\dim_{\Qp}(V)$. 
La d\'efinition de $\bst$ donn\'ee ci-dessus 
est l\'eg\`erement diff\'erente de celle de
Fontaine, mais le foncteur $\dst$ est le-m\^eme.
Alors $\dst(V)$ est muni d'un Frobenius, d'une 
filtration et d'un op{\'e}rateur de monodromie $N=-d/d\logpi$\index{N@$N$}
qui v{\'e}rifie $N\phi=p\phi N$
(voir \cite{Ts99} pour une justification du signe ``$-$''), 
et $(\bdR\otimes_{\Qp} V)^{G_K}=\ddR(V)=K \otimes_F \dst(V)$.
De plus $V$ est cristalline si et seulement si
elle est semi-stable et $N=0$ sur $\dst(V)$. On utilisera aussi le $F$-espace
vectoriel $\dst^+(V)=(\btstplus{} \otimes_{\Qp} V)^{G_K}$. 

\subsection{L'anneau $\bt$ et ses sous-anneaux}
Soit $\at$\index{Atilde@$\at$} l'anneau des vecteurs de Witt {\`a} coefficients dans
$\et$ et $\bt=\at[1/p]$\index{Btilde@$\bt$}.
Soit $\mathbf{A}_F$\index{AF@$\mathbf{A}_F$}
le compl{\'e}t{\'e} de $\OO_F[\pi,\pi^{-1}]$ dans
$\atplus$ pour la topologie de celui-ci, c'est aussi le compl{\'e}t{\'e}
$p$-adique de $\OO_F[[\pi]][\pi^{-1}]$. 
C'est un anneau de valuation discr{\`e}te complet
dont le corps r{\'e}siduel est $\mathbf{E}_F$.
Soit $\mathbf{B}$\index{B@$\mathbf{B}$} le compl{\'e}t{\'e}
pour la topologie $p$-adique
de l'extension maximale non ramifi{\'e}e de 
$\mathbf{B}_F=\mathbf{A}_F[1/p]$\index{BF@$\mathbf{B}_F$}
dans $\bt$. On d{\'e}finit alors 
$\mathbf{A}=\mathbf{B}
\cap \at$\index{A@$\mathbf{A}$} 
et $\aplus=\mathbf{A} \cap \atplus$. Ces anneaux sont 
munis d'une action de Galois et d'un Frobenius d{\'e}duits de ceux de $\et$. 
On pose $\mathbf{A}_K=\mathbf{A}^{H_K}$ et
$\mathbf{B}_K=\mathbf{A}_K[1/p]$\index{AK@$\mathbf{A}_K$,
$\mathbf{B}_K$}.
Quand $K=F$ les deux d{\'e}finitions co{\"\i}ncident.

Si $V$ est une repr{\'e}sentation $p$-adique de $G_K$ soit $\dfont(V)
=(\mathbf{B} \otimes_{\Qp} V)^{H_K}$\index{D@$\dfont(V)$}. 
On sait \cite{Fo91} que $\dfont(V)$ est un
$\mathbf{B}_K$-espace vectoriel de dimension $d=\dim(V)$ muni d'un Frobenius
et d'une action r{\'e}siduelle de $\Gamma_K$ 
qui commutent (c'est un $(\phi,\Gamma_K)$
module) et que l'on peut r{\'e}cup{\'e}rer $V$ gr{\^a}ce {\`a} la formule
$V=(\dfont(V) \otimes_{\mathbf{B}_K} \mathbf{B})^{\phi=1}$. 

Le corps $\mathbf{B}$ est une extension totalement 
ramifi{\'e}e (l'extension r{\'e}siduelle est purement
ins{\'e}pa\-rable
(elle est ``radicielle'')) de degr{\'e} $p$ de $\phi(\mathbf{B})$.
Le Frobenius $\phi: \mathbf{B} \ra \mathbf{B}$ 
est injectif, mais n'est donc pas surjectif, et il
est utile d'en d{\'e}finir un inverse {\`a} gauche par la formule:
\[ \psi(x)=\phi^{-1}(p^{-1} \on{Tr}_{\mathbf{B}/\phi(\mathbf{B})}(x)).\]\index{psi@$\psi$}

Tout {\'e}l{\'e}ment $x\in\bt$ peut s'{\'e}crire de mani{\`e}re unique
sous la forme $x=\sum_{k \gg -\infty} p^k [x_k]$ o{\`u} les $x_k$ sont des
{\'e}l{\'e}ments de $\et$. Si $r>0$, on pose:
\[\btdag{,r}=\{ x \in \bt,\ 
\lim_{k \ra +\infty} \vale(x_k)+\frac{pr}{p-1}k 
= +\infty \}\]\index{Btildedagr@$\btdag{,r}$}
Si $r\in\RR$ on d{\'e}finit  $n(r)$\index{nr@$n(r)$}
comme {\'e}tant
le plus petit entier $n$ tel que
$r \leq r_n = p^{n-1}(p-1)$\index{rn@$r_n$}.
La s{\'e}rie $\sum_{k \gg -\infty} p^k [x_k]$ converge dans 
$\bdR$ si et seulement si
la s{\'e}rie $\sum_{k \gg -\infty} p^k x_k^{(0)}$ converge dans $\Cp$
(voir \cite[II.25]{Co99}, et la construction qui pr{\'e}c{\`e}de \ref{kertheta00}). 
On en d{\'e}duit notamment pour $n$ entier tel que $p^{n-1}(p-1) \geq r$ une
application injective 
(on montrera cela et plus en \ref{kertheta00}, \ref{iotainj})
$\iota_n=\phi^{-n}: \btdag{,r} 
\ra \bdR^+$\index{iotan@$\iota_n$} 
qui envoie
$\sum_{k \gg -\infty} p^k [x_k]$ sur la somme de la s{\'e}rie
$\sum_{k \gg -\infty} p^k [x_k^{p^{-n}}]$ dans $\bdR^+$. 
Soient 
$\bdag{,r} = \mathbf{B} \cap \btdag{,r}$\index{Bdagr@$\bdag{,r}$}, 
$\btdag{} =\cup_{r \geq 0} \btdag{,r}$\index{Btildedagr@$\btdag{}$} et
$\bdag{} = \cup_{r\geq 0} \bdag{,r}$\index{Bdag@$\bdag{}$}. 
Enfin $\atdag{,r}$\index{Atildedagr@$\atdag{,r}$} 
est l'ensemble des $x \in \btdag{,r} \cap \at$ tels qu'en plus,
$\vale(x_k)+\frac{pr}{p-1}k \geq 0$ pour tout $k \geq 0$, 
$\adag{,r}=\atdag{,r} \cap
\mathbf{A}$\index{Adagr@$\adag{,r}$}, 
et $\adag{}=\atdag{} \cap
\mathbf{A}$\index{Adag@$\adag{}$} o{\`u} $\atdag{}=\btdag{} \cap
\at$\index{Atildedag@$\atdag{}$}. On pose
$\bdag{,r}_K=(\bdag{,r})^{H_K}$\index{BdagrK@$\bdag{,r}_K$},
$\adag{,r}_K=(\adag{,r})^{H_K}$\index{AdagrK@$\adag{,r}_K$},
$\btdag{,r}_K=(\btdag{,r})^{H_K}$\index{BtdagrK@$\btdag{,r}_K$} et
$\atdag{,r}_K=(\atdag{,r})^{H_K}$\index{AtdagrK@$\atdag{,r}_K$}
\index{BtdagK@$\btdag{}_K$}\index{AtdagK@$\atdag{}_K$}
(en g\'en\'eral, si $R$ est un anneau muni d'une action de $H_K$, 
on pose $R_K=R^{H_K}$).

Pour situer ces anneaux, le lecteur est invit\'e \`a se reporter au diagramme des anneaux de
p\'eriodes qui se trouve en appendice.

\begin{remark}
La notation adopt{\'e}e ici diff{\`e}re un peu de celle de
\cite{CC98,CC99}. Ce que nous appelons $\btdag{,r}$ (repectivement
$\btdag{,r_n}$) y est not{\'e} 
$\widetilde{\mathbf{B}}^{\dagger}_{r^-}$ (respectivement 
$\widetilde{\mathbf{B}}^{\dagger,n}$).
\end{remark}

On dit qu'une repr{\'e}sentation $p$-adique $V$ est surconvergente si 
$\dfont(V)$ poss{\`e}de une base sur $\mathbf{B}_K$ constitu{\'e}e
d'{\'e}l{\'e}ments de $\ddag{}(V) = (\bdag{} \otimes_{\Qp} V)^{H_K}$\index{Ddag@$\ddag{}(V)$}.
Rappelons le r{\'e}sultat principal \cite{CC98,CC99} {\`a} ce sujet:
\begin{theorem}\label{toutsur}
Toute repr{\'e}sentation $V$ de $G_K$ est
surconvergente, c'est-{\`a}-dire qu'il existe $r(V)$\index{rV@$r(V)$} tel que
$\dfont(V)=\mathbf{B}_K \otimes_{\bdag{,r(V)}_K} \ddag{,r(V)}(V)$.
\end{theorem}

On a d'autre part une description assez pr{\'e}cise des anneaux
$\bdag{,r}_K$, comme le montre la proposition suivante:

\begin{proposition}\label{rappeldag}
Il existe $n(K) \in \NN$\index{nK@$n(K)$} et
$\pi_K \in \adag{,r_{n(K)}}_K$\index{piK@$\pi_K$} dont l'image modulo
$p$ est une uniformisante $\pibar_K$ de $\e_K$. 
Si $r \geq r_{n(K)}$\index{rnK@$r_{n(K)}$},
alors tout {\'e}l{\'e}ment $x\in \bdag{,r}_K$ peut s'{\'e}crire
$x=\sum_{k \in \ZZ} a_k \pi_K^k$ o{\`u} $a_k \in F$ et o{\`u}
la s{\'e}rie
$\sum_{k \in \ZZ} a_k T^k$ est
holomorphe et born{\'e}e
sur la couronne 
$\{ p^{-1/e_K r} \leq |T| < 1 \}$.
\end{proposition}

\section{Les anneaux $\btrig{}{}$ et $\btst{}{}$}
Soit $\alpha<1$  et $\mathcal{B}_F^{\alpha}$ l'ensemble des s{\'e}ries
$\sum_{k \in \ZZ} a_k X^k$ avec $a_k \in F$ une suite born{\'e}e 
telle que tout $\rho \in [\alpha;1[$, on ait 
$\lim_{k \ra \pm \infty} |a_k|\rho^k =
0$, et $\mathcal{B}_F = \cup_{\alpha < 1} \mathcal{B}_F^{\alpha}$.
Les rappels du chapitre pr{\'e}c{\'e}dent montrent que l'application
$f \mapsto f(\pi_K)$ de $\mathcal{B}_F$ dans $\bdag{}_K$ est
un isomorphisme. 

Soit $\mathcal{H}_F^{\alpha}$ l'ensemble des s{\'e}ries
$\sum_{k \in \ZZ} a_k X^k$ avec $a_k \in F$ 
et telles que tout $\rho \in [\alpha;1[$, 
on ait $\lim_{k \ra \pm \infty} |a_k|\rho^k =
0$. Alors $\mathcal{H}_F=\cup_{\alpha<1}\mathcal{H}_F^{\alpha}$ 
est l'anneau de Robba 
{\`a} coefficients dans $F$ dont il a {\'e}t{\'e} question dans l'introduction:
\[ \mathcal{H}_F = \cup_{r \geq 0} \cap_{s \geq r} \OO_F \left\{
\frac{p}{T^r},
\frac{T^s}{p} \right\} \left[ \frac{1}{p} \right]. \]
Afin de faire le lien entre repr{\'e}sentations $p$-adiques et {\'e}quations
diff{\'e}rentielles,
on ``d{\'e}finit'' donc un anneau $\btrig{}{} \supset \btdag{}$ tel
que $(\btrig{}{})^{H_K}$ contienne $\mathcal{H}_F(\pi_K)$:
\[  \btrig{}{} = \cup_{r \geq 0} \cap_{s \geq r} \atplus \left\{
\frac{p}{[\pibar^r]},
\frac{[\pibar^s]}{p} \right\} \left[ \frac{1}{p} \right]. \]
Cet anneau contient $\btdag{}$, et aussi $\btrigplus{}$: il va donc
servir comme interm{\'e}diaire entre les $(\phi,\Gamma_K)$-modules 
et la th{\'e}orie de Hodge $p$-adique.
Il y a un certain nombre de choses techniques {\`a} d{\'e}montrer pour
v{\'e}rifier que sa d{\'e}finition a bien un sens (les trois premiers
paragraphe de ce chapitre), et le lecteur est invit{\'e} {\`a} ne les consulter qu'en
cas de besoin.

\subsection{Les anneaux $\at_I$}
Dans tout ce chapitre, $r$ et $s$ sont deux {\'e}l{\'e}ments de 
$\NN[1/p] \cup \{+\infty\}$ tels que 
$r \leq s$. Rappelons que pour $n \geq 0$ on a pos{\'e}
$r_n=(p-1)p^{n-1}$. Dans toute la suite, on notera $[x]^r$ pour
$[x^r]$, m{\^e}me si $r$ n'est pas entier.
On convient que $p/[\pibar]^{+\infty}=1/[\pibar]$ et que
$[\pibar]^{+\infty}/p=0$. Soient 
\begin{align*} 
\at_{[r;s]} &= \atplus 
\{ \frac{p}{[\pibar]^r} , \frac{[\pibar]^s}{p} \} \\ 
:&=
\atplus\{X,Y\}/([\pibar]^r X-p, pY-[\pibar]^s, XY-[\pibar]^{s-r}) \\
\text{et}\ \ \bt_{[r;s]} &=
\at_{[r;s]}[1/p]
\end{align*}
\index{Atilders@$\at_{[r;s]}$} 
\index{Btilders@$\bt_{[r;s]}$}
o{\`u}, si $A$ est un anneau complet pour la topologie $p$-adique,
$A\{X,Y\}$ d{\'e}note la compl{\'e}tion $p$-adique de $A[X,Y]$
c'est-{\`a} dire que $A\{X,Y\}=\{ \sum_{i,j \geq 0}a_{ij} X^iY^j \}$ o{\`u}
$a_{ij}$ est une suite qui tend vers $0$ selon le filtre des 
compl{\'e}mentaires des parties finies. 

\begin{lemma}
Soit $A=\atplus\{X,Y\}$ et $I=([\pibar]^r X-p, pY-[\pibar]^s,
XY-[\pibar]^{s-r})$. 
Alors
\begin{enumerate}
\item $I \cap p^n A = p^n I$;
\item $I$ est ferm{\'e} dans $A$ pour la topologie $p$-adique.
\end{enumerate}
\end{lemma}

\begin{proof}
Le (2) suit du (1) puisque $I$ est complet pour la topologie
$p$-adique et que le (1)
montre qu'une suite d'{\'e}l{\'e}ments de
$I$ qui tend vers $0$ dans $A$ tend vers $0$ dans $I$. Montrons donc
le (1). Soit $f(X,Y) =  \sum_{i,j \geq 0}a_{ij} X^iY^j$, $f(X,Y) \in p^n A$ ce
qui revient {\`a} dire que $p^n | a_{ij}$. Quitte {\`a} modifier $f$ par
des {\'e}l{\'e}ments de $p^n(XY-[\pibar]^{s-r}) \subset p^n I$ on peut
supposer que $f(X,Y)=\sum_i a_i X^i + \sum_j b_j Y^j$ o{\`u} $p^n$ divise
$a_i$ et $b_j$. Si $f(X,Y) \in I$ c'est que l'on peut {\'e}crire
\begin{multline*}
f(X,Y)=a(X,Y)([\pibar]^r
X-p) \\
+b(X,Y)(pY-[\pibar]^s)+c(X,Y)(XY-[\pibar]^{s-r})
\end{multline*}
et les relations 
\begin{align*}
Y([\pibar]^r X -p) &= [\pibar]^r(XY-[\pibar]^{s-r})-(pY-[\pibar]^s) \ \text{et} \\
X(pY-[\pibar]^s) &= p(XY-[\pibar]^{s-r})-([\pibar]^r X-p)[\pibar]^{s-r}
\end{align*}
montrent que quitte {\`a} modifier $c(X,Y)$ on peut supposer que $a(X,Y)=a(X,0)$
et $b(X,Y)=b(0,Y)$. On voit alors que $c(X,Y)=0$ et donc finalement
que $f(X,Y)=a(X)([\pibar]^r X-p)+b(Y)(pY-[\pibar]^s)$. Montrons que
$a(X) \in p^n A$ (la preuve pour $b(Y)$ est la m{\^e}me). Posons
$a(X)=\sum_i c_i X^i$. On a alors $a_0=-p c_0$ et $a_i=[\pibar]^r
c_{i-1}-pc_i$ ce qui fait si $i \geq 0$ et $0 \leq j \leq n-1$,
alors $p^n$ divise $c_{i+j}[\pibar]^r-pc_{i+j+1}$ et donc aussi
$\sum_{j=0}^{n-1}[\pibar]^{r(n-1-j)}
p^j(c_{i+j}[\pibar]^r-pc_{i+j+1})=[\pibar]^{rn} c_i - p^n
c_{i+n}$ et donc $p^n$ divise $c_i$ dans $\atplus$. 
\qed\end{proof}

\begin{corollary}
L'anneau $\at_{[r;s]}$ est s{\'e}par{\'e} pour la topologie $p$-adique (il
est clairement complet). De plus on a une application naturelle 
surjective $\at_{[r;s]}$ dans le compl{\'e}t{\'e} $p$-adique de 
$\atplus[p/[\pibar]^r,[\pibar]^s/p]$ dont le noyau est l'image de l'adh{\'e}rence
de $I$ dans $\at_{[r;s]}$ et donc nul, ce qui fait que $\at_{[r;s]}$ s'identifie aussi au
compl{\'e}t{\'e} $p$-adique de $\atplus[p/[\pibar]^r,[\pibar]^s/p]$.
\end{corollary}

\begin{lemma}
Tout {\'e}l{\'e}ment de $\at_{[r;s]}$ peut s'{\'e}crire de la mani{\`e}re suivante:
\[ \sum_{k \geq 0} \left( \frac{p}{[\pibar]^r} \right)^k a_k + \sum_{k > 0}
 \left( \frac{[\pibar]^s}{p} \right)^k b_k \]
avec $(a_k)$, $(b_k)$ deux suites de $\atplus$ qui convergent vers $0$
(cette {\'e}criture est bien s{\^u}r non-unique).
\end{lemma}

\begin{proof}
C'est une cons{\'e}quence imm{\'e}diate de la d{\'e}finition. On se
contentera de remarquer que \[ \left(\frac{p}{[\pibar]^r}\right)^k \cdot 
\left(\frac{[\pibar]^s}{p}\right)^{\ell} =  \begin{cases} 
{[\pibar]}^{k(s-r)} \left(\frac{[\pibar]^s}{p}\right)^{\ell-k} \text{ si $k
    \leq \ell$} \\
{[\pibar]}^{\ell(s-r)} \left(\frac{p}{[\pibar]^r}\right)^{k - \ell}
 \text{ si $k
    \geq \ell$}
\end{cases} \]
\qed\end{proof}

\begin{lemma}
Si $\rho$ et $\sigma$ sont deux {\'e}l{\'e}ments de $\etplus$ qui
v{\'e}rifient
$\vale(\rho)=pr/(p-1)$ et $\vale(\sigma)=ps/(p-1)$, alors
$\at_{[r;s]}=\atplus\{p/[\rho],[\sigma]/p\}$.
\end{lemma}

\begin{proof}
C'est {\'e}vident. 
\qed\end{proof}

Remarquons que si on a $r_1 \leq r_2 \leq s_2 \leq s_1$, alors il y a une
inclusion (les deux anneaux en pr{\'e}sence sont int{\`e}gres et ont
m{\^e}me corps des fractions):
\[\atplus[\frac{p}{[\pibar]^{r_1}},\frac{[\pibar]^{s_1}}{p}] 
\hookrightarrow
\atplus[\frac{p}{[\pibar]^{r_2}},\frac{[\pibar]^{s_2}}{p}] \]

\begin{lemma}
Cette inclusion se prolonge en un morphisme 
$\at_{[r_1;s_1]} \ra \at_{[r_2;s_2]}$ qui est toujours injectif.
\end{lemma}
 
\begin{proof}
Le morphisme du haut se factorise en $\at_{[r_1;s_1]} \ra
\at_{[r_1;s_2]} \ra \at_{[r_2;s_2]}$ et on peut donc supposer
$r_1=r_2$ ou $s_1=s_2$. Supposons par exemple $r_1=r_2=r$
(l'autre cas se traite de la m{\^e}me mani{\`e}re). Alors 
il suffit de montrer que le morphisme compos{\'e} 
$\at_{[r;s_1]} \ra
\at_{[r;s_2]} \ra \at_{[r;r]}$ est injectif (tout ceci pour simplifier
la notation). On va donc montrer que si $s \geq r$, alors 
$\at_{[r;s]} \ra \at_{[r;r]}$ est injectif. Soient
$\alpha=[\overline{\alpha}]$ et 
$\beta=[\overline{\beta}]$ avec $\overline{\alpha},\overline{\beta}
\in \et$ tels que $\vale(\overline{\alpha})=r$ et 
$\vale(\overline{\beta})=s-r$ de telle
sorte que 
\begin{align*}
\at_{[r;s]} &= \at\{X,Y\}/(\alpha X-p,pY-\alpha\beta,XY-\beta) \ \text{et} \\
\at_{[r;r]} &= \at\{X,Y\}/(\alpha X-p,pY-\alpha,XY-1).
\end{align*}
Le fait que l'application naturelle $f(X,Y) \mapsto f(X,\beta Y)$ du
premier anneau dans le second est injective est {\'e}quivalent au fait
que si $f(X,\beta Y) \in (\alpha X-p,pY-\alpha,XY-1)$, alors
$f(X,\beta Y)$ appartient {\`a} l'id{\'e}al de $\at\{X,\beta Y\}$
engendr{\'e} par  
$(\alpha X-p, pY\beta -\alpha\beta, XY\beta -\beta)$,
ce que nous allons maintenant d{\'e}montrer.

Commen{\c{}}{c}ons par remarquer que quitte {\`a} modifier $f(X,Y)$
par des {\'e}l{\'e}\-ments de l'id{\'e}al $(\alpha X-p,pY -\alpha,XY-1)$ 
on peut supposer (en rempla{{\c c}}ant $XY$ par $1$)
que \[ f(X,\beta Y)=\sum_{i=0}^{+\infty} \gamma_i
X^i+\sum_{j=1}^{+\infty} \delta_j \beta^j Y^j. \] 
Supposons  que l'on ait
\[ f(X,\beta Y)=a(X,Y)(\alpha X -p)+b(X,Y)(pY-\alpha)+c(X,Y)(XY-1). \] 
Les relations 
\begin{align*}
Y(\alpha X -p) &= \alpha(XY-1)-(pY-\alpha)\ \text{et} \\
X(pY-\alpha) &= p(XY-1)-(\alpha X-p) 
\end{align*} 
montrent que l'on peut supposer,
quitte {\`a} modifier $c(X,Y)$, que $a(X,Y)=a(X,0):=a(X)$ et $b(X,Y)=b(0,Y):=b(Y)$.
Comme $f(X,\beta Y)$ ne contient pas de terme en $XY$ c'est alors que
$c(X,Y)=0$. On a donc montr{\'e} que $f(X,\beta Y)=a(X)(\alpha X
-p)+b(Y)(pY-\alpha)$. Posons $b(Y)=\sum_{j=0}^{+\infty} b_j Y^j$. Pour
terminer la d{\'e}monstration il faut montrer que $b_j$ est un multiple
de $\beta^{j+1}$.
Un calcul facile montre que si $j\geq 1$, alors $\delta_j
\beta^j=(pb_{j-1}-\alpha b_j)$  ce qui fait que $\beta^j$ divise
$pb_{j-1}-\alpha b_j$ dans $\atplus$ et donc que $\beta^{j+1}$ divise
$\sum_{k=0}^{n} p^{n-k}\alpha^k (pb_{j+k}-\alpha
b_{j+k+1})=p^{n+1}b_j-\alpha^{n+1}b_{j+n+1}$. Reste {\`a} choisir $n$
assez grand pour que $\beta^{j+1}$ divise $\alpha^{n+1}$ ce qui montre
alors que $\beta^{j+1}$ divise $b_j$.
\qed\end{proof}

On utilise ces injections
pour d{\'e}finir, si $I$ est un intervalle de $\RR \cup \{ +\infty \}$:
$\bt_I = \cap_{[r;s] \subset I \cap \RR} \bt_{[r;s]}$\index{BtildeI@$\bt_I$}. 

Soient $I \subset J$ deux intervalles ferm{\'e}s, ce qui fait que $\bt_J
\subset \bt_I$, on d{\'e}finit une valuation 
$p$-adique $V_I$\index{VI@$V_I$}
sur $\bt_J$ en d{\'e}cidant
que $V_I(x)= 0$ si et seulement si $x \in \at_I-p \at_I$ et que l'image
de $V_I$ est $\ZZ$.
Par d{\'e}finition $\bt_I$ muni de $V_I$ est un espace de Banach $p$-adique.
De plus le compl{\'e}t{\'e} de $\bt_J$ pour $V_I$ s'identifie a $\bt_I$.
\begin{remark}\label{valanno}
Comme $\at_I$ est un anneau on a 
notamment $V_I(xy) \geq V_I(x)+V_I(y)$.
\end{remark}

Le but de cette partie est de d{\'e}gager quelques propri{\'e}t{\'e}s de ces
anneaux.
Commen{{\c c}}ons par remarquer que le groupe de Galois $G_F$ agit sur
$\atplus$ et que cette action s'{\'e}tend {\`a} l'anneau
$\atplus[p/[\pibar]^r,[\pibar]^s/p]$ et le stabilise, ce qui fait que
l'action de $G_F$ s'{\'e}tend par continuit{\'e} {\`a} une action par
isom{\'e}tries sur son compl{\'e}t{\'e}
$p$-adique et par suite {\`a} tous les $\at_I$ et $\bt_I$. 

De m{\^e}me le Frobenius $\phi$ s'{\'e}tend en un morphisme 
\[ \phi: 
\atplus[\frac{p}{[\pibar]^r},\frac{[\pibar]^s}{p}] \ra
\atplus[\frac{p}{[\pibar]^{pr}},\frac{[\pibar]^{ps}}{p}] \]
et se prolonge donc en une application de $\at_I$ dans $\at_{pI}$ pour
tout $I$.

\begin{lemma}\label{vw}
Si $I \subset [r;+\infty]$, alors $\btdag{,r} \subset \bt_I$ et si $x
\in \btdag{,r}$ s'{\'e}crit $x=\sum_{k \gg 0} p^k [x_k]$, alors
la valuation
\[ W_I(x) = 
\inf_{\alpha \in I} \inf_{k \in \ZZ} 
k+ \frac{p-1}{p\alpha} \vale(x_k) \]\index{WI@$W_I$}
v{\'e}rifie $V_I(x)=\lfloor W_I(x) \rfloor$ o{\`u} $\lfloor a \rfloor$ est
le plus grand entier $\leq a$.
\end{lemma}

\begin{proof}
Le premier point suit de la d{\'e}finition et de plus si $x
\in \btdag{,r}$ v{\'e}rifie $W_I(x) \geq 0$, alors la somme qui le
d{\'e}finit converge dans $\at_I$ ce qui fait que $V_I(x) \geq 0$. 
Reste {\`a} montrer que si $x \in \at_I$, alors $W_I(x) \geq 0$. Comme
$\at_I$ est le compl{\'e}t{\'e} $p$-adique de
$\atplus[p/[\pibar]^r,[\pibar]^s/p]$ il suffit de montrer que $W_I(x)
\geq 0$ si $x \in \atplus$, si $x=p/[\pibar]^r$ et si $x=[\pibar]^s/p$
ce qui est clair. Comme $W_I(p)=1$ on en d{\'e}duit que $\lfloor
W_I(\cdot ) \rfloor$ est une valuation $p$-adique dont l'image est
$\ZZ$ et telle que $W_I(x) \geq 0$ si et seulement si $x \in \at_I$ ce
qui fait que $V_I(x)=\lfloor W_I(x) \rfloor$.
\qed\end{proof}

\begin{example}
Beaucoup de ces anneaux sont d{\'e}j{\`a} connus:
\begin{enumerate}
\item $\at_{[0;r_0]}=\amax^+$ et $\bt_{[0;r_0]}=\bmax^+$;
\item $\btrigplus{}=\bt_{[0;+\infty[}$;
\item $\atplus=\at_{[0;+\infty]}$ et $\btplus=\bt_{[0;+\infty]}$;
\item $\at=\at_{[+\infty;+\infty]}$ et $\bt=\bt_{[+\infty;+\infty]}$;
\item $\atdag{,r}=\at_{[r;+\infty]}$ et  $\btdag{,r}=\bt_{[r;+\infty]}$.
\end{enumerate}
\end{example}

\begin{proof}
Le $(2)$ est une cons{\'e}quence du $(1)$ et du fait que par d\'efinition, 
$\btrigplus{} = \cap_{n
= 0}^{+\infty} \phi^n(\bmax^+)$. 
Les $(3)$ et $(4)$ sont {\'e}vidents, et le $(5)$ est
contenu dans \cite[remarque II.1.3]{CC98}. Reste {\`a} montrer le $(1)$
qui suit du fait que par d{\'e}finition $\amax^+=\atplus\{
[\tilde{p}]/p-1\}$ et $\at_{[0;r_0]}=\atplus\{
[\tilde{p}]/p\}$ (et $A\{X\}=A\{X-1\}$ puisque $A[X]=A[X-1]$).
\qed\end{proof}

\begin{lemma}\label{modp}
On a $\at_{[r;s]}/(p) = \etplus[X,\pibar^{s-r}X^{-1}]/(\pibar^s, \pibar^r
X)$. Notamment si $r=s$, alors
$\at_{[r;r]}/(p) = \etplus/(\pibar^r)[X,X^{-1}]$.
\end{lemma}

\begin{proof}
Soit $A=\atplus \{ X,Y \}$ et $I=(XY-[\pibar]^{s-r}, p-X[\pibar]^r, [\pibar]^s-p Y)$ de
telle sorte que $\at_{[r;s]}$ s'identifie {\`a}
$A/I$ et donc que $\at_{[r;s]}/(p)=(A/I)/(p)$.
On a une suite exacte $0 \ra I \ra A \ra A/I \ra 0$ et la
multiplication par $p$ induit un diagramme:
\[ \begin{CD}
0 @>>> I @>>> A @>>> A/I @>>> 0 \\
 & &     @V{p}VV @V{p}VV @V{p}VV \\
0 @>>> I @>>> A @>>> A/I @>>> 0 \\
 & &     @VVV @VVV @VVV \\
 & &  I/p @>>> A/p @>>> (A/I)/p @>>> 0 \\
\end{CD} \] 
et comme $A/I$ est sans $p$-torsion, le lemme du serpent montre que
$(A/I)/p$ s'identifie au quotient du $A/p$ par l'image de $I$ dans ce
dernier.
Dans notre situation on a $A/p=\etplus[X,Y]$ et l'image de $I$
s'identifie {\`a} $(XY-\pibar^{s-r}, -X\pibar^r, \pibar^s)$ d'o{\`u} le lemme.
\qed\end{proof}

\begin{remark}
Attention au fait que dans le lemme pr\'ec\'edent, $(\pibar^s, \pibar^r X)$ est l'id\'eal 
de l'anneau $\etplus[X,\pibar^{s-r}X^{-1}]$ 
engendr\'e par 
$\pibar^s$ et $\pibar^r X$.
\end{remark}

\subsection{Plongement des $\at_I$ dans $\bdR^+$}
On va maintenant d{\'e}finir des morphismes de $\at_{[r;s]}$
dans $\bdR^+$. On va montrer que si $r_n \in I$, alors 
l'application $\phi^{-n}$ r{\'e}alise une injection de $\bt_I$ dans
$\bdR^+$. Pour cela,
soit $J_n=[r_n;r_n]$\index{Jn@$J_n$} 
avec $n \geq 0$. Si $x\in
\at_{J_0}$, alors on peut {\'e}crire 
\begin{multline*} 
x = \sum_{k \geq 0} a_k \left( \frac{p}{[\tilde{p}]} \right)^k
+ \sum_{j \geq 0} b_j \left( \frac{[\tilde{p}]}{p} \right)^j  \\
= \sum_{k \geq 0} a_k \left[\left( \frac{p}{[\tilde{p}]}-1 \right) +1
\right]^k
+ \sum_{j \geq 0} b_j \left[ \left( 
\frac{[\tilde{p}]}{p}-1 \right) +1  \right] ^j  \\
= \sum_{\ell \geq 0} \left( \frac{p}{[\tilde{p}]} - 1 \right)^{\ell} 
\sum_{k \geq \ell} \binom{k}{\ell} a_k
+ \sum_{m \geq 0} \left( \frac{[\tilde{p}]}{p} - 1 \right)^m
\sum_{j \geq m} \binom{j}{m} b_j 
\end{multline*}
et comme les $a_j$ et $b_k$ tendent vers $0$, les s{\'e}ries
\[ \sum_{k \geq \ell} \binom{k}{\ell} a_k \ \text{et}\ 
\sum_{j \geq m} \binom{j}{m} b_j\] convergent dans $\atplus$
et la s{\'e}rie du haut est donc convergente pour 
la topologie $\ker(\theta)$-adique 
et converge dans $\bdR^+$ vers un {\'e}l{\'e}ment 
$\iota_0(x) \in \bdR^+$.\index{iota0@$\iota_0$}

\begin{proposition}\label{kertheta00}
L'application $x \mapsto
\iota_0(x)$ est un morphisme injectif de $\at_{J_0}$ dans $\bdR^+$
et si $r_0 \in I$, alors
le noyau du morphisme compos{\'e} $\theta \circ
\iota_0 : \at_I \ra \Cp$ est
$\ker(\theta\circ\iota_0: \at_I \ra \Cp)=([\tilde{p}]/p-1)\at_I$.
\end{proposition} 

\begin{proof}
Montrons tout d'abord que $\ker(\theta \circ \iota_0: \at_I \ra \Cp) \subset 
([\tilde{p}]/p-1)\at_I$.
Soit donc
$I=[rr_0;sr_0]$ avec $r \leq 1 \leq s$ et $x \in \at_I$
tel que $\theta \circ \iota_0(x)=0$. 
On peut {\'e}crire
\begin{multline*} 
x = \sum_{k \geq 0} a_k \left( \frac{p}{[\tilde{p}]^r} \right)^k
+ \sum_{j \geq 0} b_j \left( \frac{[\tilde{p}]^s}{p} \right)^j  \\
= \sum_{\ell \geq 0} \left( \frac{p}{[\tilde{p}]^r} - [\tilde{p}^{1-r}] \right)^{\ell} 
\sum_{k \geq \ell} \binom{k}{\ell}  [\tilde{p}^{1-r}]^{k-\ell} a_k \\
+ \sum_{m \geq 0} \left( \frac{[\tilde{p}]^s}{p} -  [\tilde{p}^{s-1}] \right)^m
\sum_{j \geq m} \binom{j}{m}  [\tilde{p}^{s-1}]^{j-m} b_j \\
=  \sum_{\ell \geq 1} \left( \frac{p}{[\tilde{p}]^r} -  [\tilde{p}^{1-r}] \right)^{\ell} 
\sum_{k \geq \ell} \binom{k}{\ell} [\tilde{p}^{1-r}]^{k-\ell} a_k    \\
+ \sum_{m \geq 1} \left( \frac{[\tilde{p}]^s}{p} -  [\tilde{p}^{s-1}] \right)^m     
\sum_{j \geq m} \binom{j}{m}  [\tilde{p}^{s-1}]^{j-m} b_j \\
+ \sum_{k \geq 0} (a_k [\tilde{p}^{1-r}]^k + b_k [\tilde{p}^{s-1}]^k).   
\end{multline*}
les deux premiers termes {\'e}tant des s{\'e}ries convergeant 
dans $\bdR^+$ (ne pas oublier que $a_k \ra 0$ et $b_j \ra 0$)
et dont la somme est dans le noyau de $\theta\circ\iota_0$, et
le troisi{\`e}me terme {\'e}tant un {\'e}l{\'e}ment de $\atplus$ 
(m{\^e}me argument pour la convergence) 
qui est annul{\'e}
par $\theta$ et qui s'{\'e}crit donc $([\tilde{p}]-p) y$ avec $y \in \atplus$.
Montrons que  \[ x-\sum_{k \geq 0} (a_k [\tilde{p}^{1-r}]^k + 
b_k [\tilde{p}^{s-1}]^k)=([\tilde{p}]/p-1)z \]
avec $z \in \at_I$.
On a 
\begin{multline*} 
x-\sum_{\ell \geq 0} (a_{\ell} [\tilde{p}^{1-r}]^{\ell} + 
b_{\ell} [\tilde{p}^{s-1}]^{\ell})=
\\ 
\sum_{\ell \geq 1} a_{\ell} \left(\left( \frac{p}{[\tilde{p}]^r} \right)^{\ell}-  
[\tilde{p}^{1-r}]^{\ell}\right)
+ \sum_{\ell \geq 1} b_{\ell} \left(\left( \frac{[\tilde{p}]^s}{p}
\right)^{\ell} -[\tilde{p}^{s-1}]^{\ell} \right)
= 
\\ \left(\frac{[\tilde{p}]-p}{p}\right) 
\Biggl( -\sum_{k \geq 1}\left( \frac{p}{[\tilde{p}]^r} \right)^{k} \sum_{\ell \geq k}
a_{\ell} [\tilde{p}^{1-r}]^{\ell-k} \\
+[\tilde{p}^{s-1}]
\sum_{j \geq 0}\left( \frac{[\tilde{p}^s]}{p} \right)^j \sum_{\ell \geq j+1}   
b_{\ell}[\tilde{p}^{s-1}]^{\ell-j-1}\Biggr)
\end{multline*}
comme le montre un petit calcul, ce qui montre l'assertion quant au
noyau de $\theta \circ \iota_0$ sur $\at_I$.

Enfin montrons que $\iota_0$ est injectif. 
Si $\iota_0(x)=0$ avec 
$x \in \at_{J_0}$,
alors $\theta \circ \iota_0(x) =0$ et donc $x$ est divisible
par $[\tilde{p}]/p-1$. Comme $\iota_0$ est un morphisme d'anneau
et que $\bdR^+$ est int{\`e}gre, c'est que $x=([\tilde{p}]/p-1) x_1$
avec $\iota_0(x_1)=0$ et $x_1 \in \at_{J_0}$. 
En it{\'e}rant ce proc{\'e}d{\'e} on en d{\'e}duit que
$x \in \cap_{n=0}^{+\infty} ([\tilde{p}]/p-1)^n 
\at_{J_0}$. Reste donc {\`a} montrer que
$\cap_{n=0}^{+\infty}([\tilde{p}]/p-1)^n 
\at_{J_0}=0$.
L'image de cette intersection dans $\at_{J_0}/(p)$
s'identifie {\`a} $\cap_{n=0}^{+\infty} (X-1)^n \etplus/(\tilde{p}) [X,X^{-1}]$ qui est
nulle. On en d{\'e}duit qu'un {\'e}l{\'e}ment de l'intersection
$\cap_{n=0}^{+\infty} ([\tilde{p}]/p-1)^n 
\at_{J_0}$ est infiniment divisible par $p$ et donc nul.
\qed\end{proof}

\begin{proposition}\label{kertheta0}
Si $x\in \at_{J_n}$ on pose
$\iota_n(x)=\iota_0(\phi^{-n}(x))$. L'application $x \mapsto
\iota_n(x)$ est un morphisme injectif de $\at_{J_n}$ dans $\bdR^+$
et si $r_n \in I$, alors
le noyau du morphisme compos{\'e} $\theta \circ
\iota_n : \at_I \ra \Cp$ est
$\ker(\theta\circ\iota_n: \at_I \ra \Cp)=([\tilde{p}^{p^n}]/p-1)\at_I$.
\end{proposition} 

\begin{proof}
C'est une cons{\'e}quence imm{\'e}diate de la proposition pr{\'e}c{\'e}dente,
{\'e}tant donn{\'e} que l'application $\phi^{-n}: \at_{J_n} \ra 
\at_{J_0}$ est une bijection, et que 
$\phi^{-n}([\tilde{p}^{p^n}]/p-1)=[\tilde{p}]/p-1$.
\qed\end{proof}

\begin{corollary}
Si $r_n \in I$, alors $\iota_n$ r{\'e}alise une
injection de $\at_I$ et $\bt_I$ dans $\bdR^+$.
\end{corollary}

\begin{remark}\label{phiqn}
Le noyau de $\theta \circ \iota_n: \bt_I \ra \Cp$ est
$\ker(\theta \circ \iota_n) = \phi^{n-1}(q)\bt_I$.
En effet, $([\tilde{p}^{p^n}]-p)=\phi^n([\tilde{p}]-p)$, 
$\phi^{n-1}(q)=\phi^n(\omega)$ et on sait que $[\tilde{p}]-p$ et
$\omega$ engendrent le m{\^e}me id{\'e}al de $\atplus$.
\end{remark} 

\begin{lemma}\label{unicite}
Les inclusions naturelles de $\amax^+$ et $\atdag{,r_0}$
dans $\at_{J_0}$ induisent
une suite exacte $0 \ra \atplus \ra \amax^+ \oplus \atdag{,r_0}
\ra \at_{J_0} \ra 0$.
\end{lemma}

\begin{proof}
La fl{\`e}che est $\amax^+ \oplus \atdag{,r_0}
\ra \at_{J_0}$ est surjective car il suffit de 
d{\'e}composer une {\'e}criture d'un {\'e}l{\'e}ment de $\at_{J_0}$ en deux. 
Ensuite $\atplus$ est contenu dans l'intersection de $\amax^+$ et de
$\atdag{,r_0}$ et il reste donc {\`a}  montrer que l'inclusion 
\[ \atplus \ra \atdag{,r_0} \cap \amax^+ \]
est aussi une surjection. 
On va d'abord montrer que c'est vrai
modulo $p \at_{J_0}$ (on remarquera que modulo $p$ la fl{\`e}che n'est plus injective).
Rappelons que les anneaux $\amax^+$ et $\atdag{,r_0}$
s'identifient {\`a}
$\atplus\{X\}/(pX-[\tilde{p}])$ et
{\`a} $\atplus\{Y\}/([\tilde{p}]Y-p)$ et que
$\at_{J_0}/(p)=\etplus/(\tilde{p})[X,X^{-1}]$.
L'image de $\atdag{,r_0}$ dans cet anneau s'identifie {\`a}
$\etplus/(\tilde{p})[1/X]$ et celle de $\amax^+$ {\`a} 
$\etplus/(\tilde{p})[X]$ ce qui fait que l'image de leur
intersection (qui est un sous-ensemble de l'intersection de leurs images)
est un sous-anneau de $\etplus/(\tilde{p})$ et donc que la fl{\`e}che 
$\atplus \ra \atdag{,r_0} \cap \amax^+$
est surjective modulo $p\at_{J_0}$. 
Si l'on prend $x$ dans $\atdag{,r_0} \cap \amax^+$ il existe
donc $y \in \atplus$ tel que $x-y \in p\at_{J_0}$. Cela veut dire que 
$x-y \in p \amax^+$ d'une part et 
$\in p\atdag{,r_0}+[\tilde{p}]\atplus$ 
d'autre part (il suffit
d'appliquer le lemme \ref{modp} {\`a} tous ces anneaux). Comme $p$ divise
$[\tilde{p}]$ dans $\amax^+$ il existe $z \in [\tilde{p}]\atplus$ tel
que $x-y-z \in p(\atdag{,r_0} \cap \amax^+)$. On conclut en
it{\'e}rant ce proc{\'e}d{\'e} (comme $\atplus$ est complet pour la
topologie $p$-adique).
\qed\end{proof}

\subsection{L'anneau $\btrig{}{}$}
Dans ce paragraphe on introduit l'anneau $\btrig{}{}$.

\begin{definition}
Les anneaux $\btrig{,r}{}$ et $\btrig{}{}$
sont d{\'e}finis par 
$\btrig{,r}{} = \bt_{[r;+\infty[}$ 
et  $\btrig{}{} = \cup_{r \geq 0} \btrig{,r}{}$.
\index{Btildedagrigr@$\btrig{,r}{}$}\index{Btildedagrig@$\btrig{}{}$} 
On munit $\btrig{,r}{}$ 
de la topologie de Fr{\'e}chet d{\'e}finie par
l'ensemble des $V_I$ o{\`u} $I$ parcourt l'ensemble des 
intervalles ferm{\'e}s de $[r;+\infty[$.
On d{\'e}finit aussi $\atrig{,r}{}$ comme {\'e}tant l'anneau des entiers de
$\btrig{,r}{}$ pour la valuation $V_{[r;r]}$.
\end{definition}

\begin{proposition}\label{kertheta1}
On a 
$\ker(\theta\circ\iota_n: \btrig{,r_n}{} \ra
\Cp)=\phi^{n-1}(q)\btrig{,r_n}{}$.
\end{proposition}

\begin{proof}
{\'E}tant
donn{\'e}e la remarque \ref{phiqn} il suffit de montrer que
$\ker(\theta\circ\iota_n: \bt_I \ra \Cp)=([\tilde{p}^{p^n}]/p-1)\bt_I$
pour tout $I \subset [r_n;+\infty[$ ce qui suit de \ref{kertheta0}.
\qed\end{proof}

\begin{lemma}
On a une suite exacte \[ 0 \ra \btplus \ra 
\btdag{,r} \oplus
\btrigplus{} \ra \btrig{,r}{} \ra 0 \]
\end{lemma}

\begin{proof}
On va d'abord montrer que si $r_n \geq r$, alors on a 
\[ 0 \ra \btplus \ra \bt_{[r;+\infty]} \oplus
\bt_{[0;r_n]} \ra \bt_{[r;r_n]} \ra 0 \]
il est clair que tout {\'e}l{\'e}ment de
$\bt_{[r;r_n]}$
s'{\'e}crit comme somme
d'{\'e}l{\'e}ments des deux autres et il faut v{\'e}rifier que deux {\'e}critures
diff{\'e}rentes diff{\`e}rent par un {\'e}l{\'e}ment de $\btplus$ ce qui revient {\`a}
montrer que $\bt_{[r;+\infty]} \cap \bt_{[0;r_n]} = \btplus$, ou encore
en appliquant $\phi^{-n}$ que 
$\bt_{[rp^{-n};+\infty]} \cap \bt_{[0;r_0]} = \btplus{}$ ce qui suit du
lemme \ref{unicite}.

Montrons maintenant le lemme. Si $x \in \btrig{,r}{}$, alors pour tout
$n$ on peut {\'e}crire (puisque $\btrig{,r}{} \subset \bt_{[r;r_n]}$):
$x=a_n+b_n$ avec $a_n \in \bt_{[r;+\infty]}$ et
$b_n \in \bt_{[0;r_n]}$. Remarquons que $x=a_{n+1}+b_{n+1}$ est une
autre {\'e}criture de ce type (puisque $a_{n+1} \in \bt_{[r;+\infty]}$ et
$b_{n+1} \in \bt_{[0;r_n]}$) et donc que $b_{n+1}-b_n \in \btplus$ ce
qui fait que quitte {\`a} modifier $a_{n+1}$ et $b_{n+1}$ par des {\'e}l{\'e}ments
de $\btplus$ on peut supposer que $a_n=a_{n+1}$ et $b_n=b_{n+1}$ ce qui fait
que $x=a+b$ avec $a \in \btdag{,r}$ et $b \in 
\cap_{n=0}^{+\infty} \bt_{[0;r_n]}
=\btrigplus{}$. 
\qed\end{proof}

\begin{proposition}\label{complet}
L'anneau $\btrig{,r}{}$ est complet pour sa topologie de Fr{\'e}chet et
contient $\btdag{,r}$ comme sous-anneau dense.
\end{proposition}

\begin{proof}
Le fait que $\btrig{,r}{}$ est complet suit du fait que chacun des
$\bt_{[r;s]}$ est complet pour $V_{[r;s]}$. Ensuite montrons que 
$\btdag{,r}$ est dense. Soient $x \in \btrig{,r}{}$ et 
$r < s < t$ trois r{\'e}els. 
Alors comme $x \in \bt_{[r;t]}$, on peut l'{\'e}crire comme
\[ x_n+\sum_{k>n} b_k \left(\frac{[\pibar]^t}{p}\right)^k \]
avec $x_n \in \btdag{,r}$ et $b_k \in\atplus_K$, si $n \gg 0$. 
On a alors 
\[ x-x_n \in \left(\frac{[\pibar]^t}{p}\right)^n \at_{[r;t]}
\subset \left(\frac{[\pibar]^t}{p}\right)^n \at_{[r;s]} \]
et un petit calcul montre qu'alors $V_{[r;s]}(x-x_n) \geq n(t/s-1)$.
Un argument d'extraction diagonale permet de trouver une suite qui converge
vers $x$ pour la topologie de Fr{\'e}chet. 
\qed\end{proof}

\begin{corollary}[Principe du maximum]\label{princmax}
Pour $x \in \btrig{,r}{}$ et $I=[s;t] \ni r$, on a $V_I(x)=\inf \{
V_{[s;s]}(x); V_{[t;t]}(x) \}$.
\end{corollary}

\begin{proof}
Un petit argument montre que c'est 
vrai avec $W_I$ {\`a} la place de $V_I$ pour $x \in
\btdag{,r}{}$; on conclut par densit{\'e} et continuit{\'e}.
\qed\end{proof}

\begin{lemma}
Si $a$ est un {\'e}lement de $\etplus$ qui v{\'e}rifie
$\lambda=\vale(a)>0$,
alors la topologie d{\'e}finie par $V_I$ sur $\atrig{,r}{}$ est
plus fine que la topologie $[a]$-adique c'est-{\`a}-dire
que si $V_I(y_i) \ra +\infty$, alors $y_i \ra 0$ pour la topologie
$[a]$-adique. De plus les topologies $[a]$-adiques et $V_{[r;r]}$-adiques
sont {\'e}quivalentes.
\end{lemma}

\begin{proof}
Soit $(y_i)$ une suite telle que $V_I(y_i) \geq n$ pour $i \geq i_0$ et
soit $m \leq \frac{nrp}{\lambda(p-1)}$. Alors un petit calcul montre que
$V_I(y_i [a]^{-m}) \geq 0$ pour $i \geq i_0$ et donc
$V_{[r;r]}(y_i [a]^{-m}) \geq 0$ pour $i \geq i_0$ ce qui revient {\`a}
dire que $y_i \in [a]^m \atrig{,r}{}$ pour $i \geq i_0$ et donc que $y_i$
tend vers $0$ pour la topologie $[a]$-adique.

De m{\^e}me si $y \in [a]^m \atrig{,r}{}$, 
alors $V_{[r;r]}(y) \geq \frac{(p-1)m \lambda}{pr}$
et donc les deux topologies sont {\'e}quiva\-lentes. On prendra garde au
fait que cela n'est plus vrai si $I$ n'est plus r{\'e}duit {\`a} $[r;r]$. 
\qed\end{proof}

\begin{corollary}
L'anneau $\atdag{,r}$ est complet pour la valuation $V_{[r;r]}$.
\end{corollary}

\begin{proof}
Dans \cite[II.1.2]{CC98} on montre que 
$\atdag{,r}$ est s{\'e}par{\'e} complet pour la
topologie $[a]$-adique si $\vale(a)>0$.
\qed\end{proof}

\subsection{Les anneaux $\btst{,r}{}$ et leurs plongements dans $\bdR^+$}
Ce paragraphe  est consacr{\'e} {\`a} la construction d'un anneau
$\btst{}{}$ qui est {\`a} $\btrig{}{}$ ce que $\bst$ est {\`a} $\bmax$.
On commence par construire une application logarithme.

\begin{proposition}
Il existe une et une seule
application $x \mapsto \log[x]$\index{log@$\log$} 
de $\et$ dans $\btrigplus{}[X]$
qui v{\'e}rifie $\log[xy]=\log[x]+\log[y]$,
$\log[x]=0$ si $x\in \overline{k}$,
$\logpi=X$ et
\[ \log[x]=\sum_{n>0} (-1)^{n-1}\frac{([x]-1)^n}{n}
\text{ si } v_p(x^{(0)}-1) \geq 1  \]
\end{proposition}

\begin{proof}
Soit $U_1$ l'ensemble des $x\in \et$
tels que $v_p(x^{(0)}-1) \geq 1$. Pour $x\in U_1$ la s{\'e}rie
$\log[x]=\sum_{n>0} (-1)^{n-1}([x]-1)^n/n$ converge dans $\bmax^+$
et $\log[xy]=\log[x]+\log[y]$ par un argument de s{\'e}ries formelles.
On en d{\'e}duit notamment que $\log[x]=\phi(\log[x]/p)$ ce qui
fait que l'image de $U_1$ par $\log$ est en fait incluse dans $\btrigplus{}$.
Si $x\in\et$ est tel que $v_p(x^{(0)}-1) >0$, alors il existe $n$ tel
que $x^{p^n} \in U_1$ et le $\log$ s'{\'e}tend donc {\`a} l'ensemble des
$x$ tels que $v_p(x^{(0)}-1) >0$.
Ensuite soit $x\in(\etplus)^*$. On peut {\'e}crire $x=x_0(1+y)$ avec
$x_0 \in \overline{k}$ et $y \in \MM_{\etplus}$ de mani{\`e}re unique
ce qui montre que $\log$ s'{\'e}tend {\`a} $(\etplus)^*$.
Enfin $\etplus$ est un anneau de valuation et le choix de $\logpi$ ach{\`e}ve
de d{\'e}terminer l'application $\log: \et \ra \btrigplus{}[X]$.
\qed\end{proof}

\begin{proposition}\label{log}
Il existe une et une seule application $\log : \atplus \ra \btrig{}{}[X]$
telle que $\log([x])=\log [x]$,
$\log(p)=0$ et $\log(xy)=\log(x)+\log(y)$.
\end{proposition}

\begin{proof}
Si $x\in\atplus$ est exactement divisible par
$p^a$, alors il existe $r$ tel que  
l'on peut {\'e}crire $x=p^a[\overline{x/p^a}](1-pz)$ avec $z\in\atdag{,r}$, o\`u 
$\overline{x/p^a}$ est l'image de $x/p^a$ dans $\etplus$,
et la s{\'e}rie
\[ \log(1-pz)=-\sum_{n \geq 1} \frac{(pz)^n}{n} \]
converge dans $\atdag{,r}$ qui
est complet pour la topologie
$p$-adique ce qui permet d'{\'e}tendre $\log$ par multiplicativit{\'e} {\`a}
$\atplus$ (et aussi {\`a} $\btplus$).
\qed\end{proof}

On pose $\btst{,r}{}=\btrig{,r}{}[X]$ muni de l'action de $G_F$ et du
Frobenius donn{\'e}s par $\phi(X)=pX$ et $g(X)=X+\log([g(\pibar)/\pibar])$
ce qui fait que l'on peut prolonger l'application
$\iota_n : \btrig{,r}{} \ra \bdR^+$
pour $n$ assez grand par
$\iota_n(X)=p^{-n}\logpi$. La proposition suivante montre que $\iota_n$ est
injectif et commute aux actions de $G_F$ et de Frobenius (l{\`a} o{\`u}
ce-dernier est d{\'e}fini),
et d{\`e}s lors on {\'e}crira
$\btst{,r}{}=\btrig{,r}{}[\logpi]$\index{Btildedaglogr@$\btst{,r}{}$}.
Soit aussi $\btst{}{} = \cup_{r \geq 0} 
\btst{,r}{}$\index{Btildedaglog@$\btst{}{}$}.

\begin{proposition}\label{iotainj}
L'application $\iota_n: \btrig{,r_n}{}[X] \ra \bdR^+$ qui {\'e}tend
$\iota_n$\index{iotan@$\iota_n$} 
par $\iota_n(X)=p^{-n} \logpi$ est injective, commute {\`a}
l'action de Galois et sa restriction {\`a} $\btstplus{}$ est $\phi^{-n}$.
\end{proposition}

\begin{proof}
Les deux derniers points sont triviaux. Pour montrer le premier,
il suffit de montrer que
l'{\'e}l{\'e}ment $\logpi$ est transcendant sur le corps des fractions de
$\iota_n(\btrig{,r_n}{})$ et cela
revient au m{\^e}me de montrer que $u=\log([\tilde{p}])$ est transcendant
sur $\iota_n(\on{Frac} \btrig{,r_n}{})$.
Montrons tout d'abord que $u \notin \iota_n(\on{Frac} \btrig{,r_n}{})$. Soit
$\beta=1-[\tilde{p}]/p$ et $S$ l'anneau des {\'e}l{\'e}ments de $\bdR^+$
qui appartiennent {\`a}
$F \otimes_{\OO_F} \atplus[[\beta]]$. Rappelons
que Fontaine a montr{\'e} que $u \notin \on{Frac} S$ (cf
\cite[4.3.2]{Bu88per}).
La d{\'e}monstration de la proposition \ref{kertheta00} montre que
si $x\in\btrig{,r_0}{}$, alors $\iota_0(x) \in S$ 
(en regardant la s\'erie qui donne l'image de $x$ dans $\bdR^+$)
et donc  
si $u \in \iota_n(\on{Frac} \btrig{,r_n}{})$, alors on
a $x,y \in \btrig{,r_0}{}$ tels que $\iota_0(y)=\iota_0(x)u$ 
(puisque $\phi^{-n}(\btrig{,r_n}{})=\btrig{,r_0}{}$) et le
r{\'e}sultat de 
Fontaine montre que cela n'est pas possible et donc que
$u \notin \iota_n(\on{Frac} \btrig{,r_n}{})$.

Montrons maintenant que $u$ est transcendant sur
$\iota_n(\on{Frac} \btrig{,r_n}{})$ pour tout $n$.
Un petit calcul montre qu'il existe $\eta: G_F
\ra \Qp$ tel que $g(u)=u+\eta(g)t$. Soit $u^d+x_{d-1}u^{d-1}+\cdots+x_0=0$
le polyn{\^o}me minimal de $u$. Alors en appliquant $g$ et en comparant les
coefficients il vient $g(x_{d-1})=x_{d-1}+d\eta(g)t$
ce qui fait que $x_{d-1}-du$ s'identifie
{\`a} un {\'e}l{\'e}ment $c$ de $\bdR$ stable par $G_F$ et donc que
$u=d^{-1}(x_{d-1}-c) \in \iota_n(\on{Frac} \btrig{,r_n}{})$, avec 
$c \in F$, et on vient de voir que
cela est impossible.
\qed\end{proof}

Par la proposition \ref{log} il existe un {\'e}l{\'e}ment $\log(\pi) \in
\btst{}{}$\index{logdepi@$\log(\pi)$}
(remarquons que l'on a
$\btst{}{} = \btrig{}{}[\log(\pi)]$
puisque la s{\'e}rie qui d{\'e}finit $\log([\overline{\pi}]/\pi)$ converge
dans $\atdag{,r_0}$).
On munit $\btst{}{}$
de l'op{\'e}rateur de monodromie $N$\index{N@$N$}
d{\'e}fini par
\[ N \left( \sum_{k=0}^d a_k \log(\pi)^k \right) =
-\sum_{k=0}^d k a_k \log(\pi)^{k-1} \]
c'est-{\`a}-dire que $N=-d/d \log(\pi)$. 
Un calcul facile montre que $N$ commute {\`a} l'action de $G_F$.

\begin{remark}\label{logpietro}
L'{\'e}l{\'e}ment $\log(\pi)$ est plus agr{\'e}able que $\logpi$, par
exemple $\iota_n(\log(\pi)) \in F_n[[t]]$ si $n \geq 1$. 
\end{remark}

\subsection{Action de $H_K$ sur $\btrig{}{}$}
Dans ce paragraphe on d{\'e}crit les invariants de $\btrig{}{}$ sous
l'action de $H_K$. Soit $\btrig{}{,K}=(\btrig{}{})^{H_K}$
\index{BtildedagrigK@$\btrig{}{,K}$}.

\begin{lemma}\label{decomp}
Soit $I$ un intervalle qui contient $[0;r]$ et $J=I \cap [r;+\infty]$.
On a une suite exacte $0 \ra \btplus_K \ra \btdag{,r}_K \oplus \bt_I^{H_K}
\ra \bt_J^{H_K} \ra 0$.
\end{lemma}
 
\begin{proof}
On a vu que l'on a une suite exacte $0 \ra \btplus \ra \btdag{,r} \oplus \bt_I
\ra \bt_J \ra 0$ et en prenant les invariants par $H_K$ on trouve
\[ 0 \ra \btplus_K \ra \btdag{,r}_K \oplus \bt_I^{H_K}
\ra \bt_J^{H_K} \overset{\delta}{\ra} H^1(H_K,\btplus) \] et pour
montrer le lemme il suffit de montrer que $\delta=0$. Si $x \in
\bt_J^{H_K}$, alors comme $[\pibar]$ est inversible dans 
$\bt_J^{H_K}$ (sauf si $r=0$, auquel cas le lemme est trivial), on a $x[\pibar]^{-1} \in
\bt_J^{H_K}$ et 
il existe $n$ tel que $\delta(x[\pibar]^{-1}) \in H^1(H_K,p^{-n}
\atplus)$ et donc $\delta(p^n (x [\pibar]^{-1}) [\pibar]) \in
H^1(H_K,W(\MM_{\etplus}))$. 
Or ce dernier espace de cohomologie est nul (\cite[IV.2.4]{Co98}
appliqu{\'e} {\`a} la repr{\'e}sentation triviale) et on en d{\'e}duit que 
$\delta(x)=p^{-n} \delta(p^n x)=0$.
\qed\end{proof}

\begin{corollary}\label{youpi}
Dans le cas o{\`u} $I=[0;+\infty[$, on obtient la suite exacte: 
\[ 0 \ra \btplus_K \ra \btdag{,r}_K \oplus (\btrigplus{})^{H_K}
\ra \btrig{,r}{,K} \ra 0 \]
\end{corollary}

\begin{lemma}\label{bmaxhk}
Si $x \in (\bmax^+)^{H_K}$, alors il existe une suite $a_i$
d'{\'e}l{\'e}ments de $\btplus_K$, qui tend vers $0$,
telle que $x=\sum_{i \geq 0} a_i (\omega/p)^i$.
\end{lemma}

\begin{proof}
On se ram{\`e}ne imm{\'e}diatement {\`a} montrer que si 
$x \in (\amax^+)^{H_K}$, alors il existe une suite  $a_i$
d'{\'e}l{\'e}ments de $\atplus_K$, 
telle que $x=\sum_{i \geq 0} a_i (\omega/p)^i$.
On va d'abord montrer que $\theta:\atplus_K \ra \OO_{\hat{K}_{\infty}}$
est surjective. Si $K=F$, c'est bien connu (on se ram{\`e}ne {\`a} montrer
que c'est vrai modulo $p$, et cela r{\'e}sulte alors du fait que les 
$\{\eps^i, i \in \ZZ[1/p] \cap [0;1[ \}$ forment une base de 
$\etplus_F/(\eps-1)\etplus_F$ (voir \cite[III.2.1]{CC98})); ensuite soient
$\iota_K(\varpi_K)$ l'{\'e}l{\'e}ment construit dans la proposition
\ref{youpla}, 
$\varpi_n=\iota_K(\varpi_K)^{(n)}$, et
$\mathfrak{a}= \{ x \in \OO_{\hat{K}_{\infty}}, v_p(x) \geq 1/p \}$.
Rappelons \cite[p. 243]{CC99}
que si $n \gg 0$ et $x \in \OO_{K_{n+1}}$, 
alors $N_{K_{n+1}/K_n}(x)-x^p \in \mathfrak{a}$.
Comme 
$\vale(\iota_K(\varpi_K))=\frac{p}{(p-1)e_K}$ on a, pour $n \gg 0$,
$v_p(\varpi_n)
=\frac{1}{p^{n-1}(p-1)e_K}$ et donc, si $n \gg 0$, alors
$\varpi_n$ est 
{\'e}gal modulo $\mathfrak{a}$ {\`a}
une uniformisante de $\OO_{K_n}$. Ceci
implique que l'application 
$\OO_{F_n}[\varpi_n] \ra \OO_{K_n}/\mathfrak{a}$
est surjective, et donc que
$\OO_{\hat{F}_{\infty}}[\varpi_n]_{n \geq 0} \ra
\OO_{\hat{K}_{\infty}}$ est surjective (puisqu'elle l'est modulo
$\mathfrak{a}$ et que les deux
$\OO_{\hat{F}_{\infty}}$-modules en pr{\'e}sence sont
complets pour la topologie $\mathfrak{a}$-adique). 
Comme $\varpi_n=\theta\circ\phi^{-n}([\iota_K(\varpi_K)])$,
ceci montre que $\theta:\atplus_K \ra \OO_{\hat{K}_{\infty}}$
est surjective. Soit maintenant $x \in (\amax^+)^{H_K}$. On d{\'e}finit
deux suites $x_i$ et $a_i$ de $(\amax^+)^{H_K}$ et $\atplus_K$ de
la mani{\`e}re suivante: $x_0=x$, et si
$i\geq 0$, alors  comme
$\theta(x_i) \in \OO_{\hat{K}_{\infty}}$, il
existe $a_i \in \atplus_K$ tel que $\theta(x_i)=\theta(a_i)$. 
Ceci montre que $x_i-a_i$ est dans le noyau de $\theta$ et est donc
divisible par $\omega/p$ dans $\amax^+$. On pose
$x_{i+1}=(p/\omega)(x_i-a_i)$.
Il est clair qu'alors, $x=\sum_{i \geq 0} a_i (\omega/p)^i$.

Montrons que $a_i \ra 0$: si $x \in
\amax^+$ alors il existe une suite $\alpha_i=o(i)$,
telle que 
$x \in p^{-\alpha_i}\atplus+(\omega/p)^i \amax^+$, et 
un petit calcul montre que
cela force $a_i \ra 0$.
\qed\end{proof}

\begin{proposition}\label{dagrig}
L'anneau $\btdag{,r}_K$ est dense dans $\btrig{,r}{,K}$ pour la
topologie de Fr{\'e}chet de ce dernier.
\end{proposition}

\begin{proof}
Cela revient {\`a} montrer que si l'on se fixe un intervalle $I=[r;r_n]$
avec $n \gg 0$, alors pour tout $x \in \btrig{,r}{,K}$ il existe une
suite $x_j \in \btdag{,r}_K$ telle que $V_I(x-x_j) \ra +\infty$. 
{\'E}tant donn{\'e}e la d{\'e}composition du corollaire \ref{youpi} il suffit de
le faire pour $x \in \bt_{[0;r_{n+m}]}^{H_K}$
avec $m \gg 0$. 
On se fixe $m$ tel que $V_I(\phi^{n+m-1}(q/p)-1)>0$ (c'est possible
car pour tout $I$ compact, $\phi^{n+m-1}(q/p) \ra 1$ pour $V_I$).
Le lemme \ref{bmaxhk} auquel on
applique $\phi^{n+m}$ montre que tout {\'e}l{\'e}ment $x \in \bt_{[0;r_{n+m}]}^{H_K}$
s'{\'e}crit $x=\sum_{i \geq 0} z_i (\phi^{n+m-1}(q/p))^i$ o{\`u} $z_i$ est
une suite de $\btplus_K$ qui tend vers $0$, et donc que 
$x=\sum_{i \geq 0} z_i (\phi^{n+m-1}(q/p)-1)^i$ o{\`u} $y_i$ est
une suite born{\'e}e de $\btplus_K$. Comme $V_I(\phi^{n+m-1}(q/p)-1)>0$,
la suite $(\phi^{n+m-1}(q/p)-1)^i$ tend vers $0$ pour $V_I$ et donc,  
pour montrer la proposition, il suffit de 
prendre $x_j=\sum_{i=0}^j z_i (\phi^{n+m-1}(q/p)-1)^i$.
\qed\end{proof}

\subsection{Les anneaux $\bnrig{}{,K}$ et $\bnst{}{,K}$}\label{rmdef}
Soit $\bnrig{,r}{,K}$\index{BdagrigrK@$\bnrig{,r}{,K}$} 
le compl{\'e}t{\'e} de
$\bdag{,r}_K$ pour la topologie de Fr{\'e}chet. 
On va donner une description nettement plus agr{\'e}able des $\bnrig{,r}{,K}$.
On a d{\'e}j{\`a} vu qu'il existe $n(K) \in \NN$ et
$\pi_K \in \adag{,r_{n(K)}}_K$ dont l'image modulo
$p$ est une uniformisante de $\mathbf{E}_K$ et que 
si $r \geq r_{n(K)}$,
alors tout {\'e}l{\'e}ment $x\in \bdag{,r}_K$ peut s'{\'e}crire
$x=\sum_{k \in \ZZ} a_k \pi_K^k$ o{\`u} $a_k \in F$ et o{\`u} la s{\'e}rie
$\sum_{k \in \ZZ} a_k T^k$ est holomorphe 
et born{\'e}e sur la couronne 
$\{ p^{-1/e_K r} \leq |T| < 1\}$.

L'anneau $\bdag{,r}_K$ est alors muni de la topologie induite
par celle de $\btrig{,r}{}$ qui devient la topologie de la convergence
sur les couronnes d{\'e}finies par un
intervalle compact
ce qui fait que $\bnrig{,r}{,K}$ est le compl{\'e}t{\'e}
de $\bdag{,r}_K$ pour cette topologie et donc que
\begin{proposition}
Soit $\mathcal{H}_F^{\alpha}(X)$ l'ensemble des s{\'e}ries
$\sum_{k \in \ZZ} a_k X^k$ avec $a_k \in F$ 
et telles que tout $\rho \in [\alpha;1[$, 
$\lim_{k \ra \pm \infty} |a_k|\rho^k = 0$ et
soit $\alpha(K,r)=p^{-1/e_Kr}$\index{alphaKr@$\alpha(K,r)$}. 
Alors l'application 
$\mathcal{H}_F^{\alpha(K,r)} \ra \bnrig{,r}{,K}$ 
qui {\`a} $f$ associe $f(\pi_K)$ est un isomorphisme.
\end{proposition}

La fin de ce paragraphe est consacr{\'e} {\`a} la d{\'e}monstration de la
proposition suivante:
\begin{proposition}\label{rm}
Si $r$ est un entier assez grand, alors
il existe des applications 
$R_k : \btrig{,r}{,K} \ra 
\phi^{-k}(\bnrig{,p^k r}{,K})$, telles que:
\begin{enumerate}
\item $R_k$ est une section continue de l'inclusion
$\phi^{-k}(\bnrig{,p^k r}{,K}) \subset  \btrig{,r}{,K}$;
\item $R_k$ est $\phi^{-k}(\bnrig{,p^k r}{,K})$-lin{\'e}aire;
\item si $x \in \btrig{,r}{,K}$, alors 
$\lim_{k \ra + \infty}R_k(x)=x$;
\item si $x\in\btrig{,r}{,K}$, et $\gamma \in \Gamma_K$, alors $\gamma
  \circ R_k (x) = R_k \circ \gamma (x)$.
\end{enumerate}
\end{proposition}

Soit $I=p^{-\infty}\ZZ \cap [0;1[$, si $i \in I$, soit
$\eps^i=(\eps^{(n)})^{p^n i}$, pour $n$ assez grand tel que $p^n i
  \in \ZZ$.
Rappelons (voir \cite[III.2]{CC98}, par exemple) que tout {\'e}l{\'e}ment
$x \in \etplus_F$ s'{\'e}crit de mani{\`e}re unique sous la forme:
$x=\sum_{i \in I} \eps^i \overline{a}_i(x)$, o{\`u} $(\overline{a}_i(x))_i$
est une suite de $\eplus_F$ qui tend vers $0$. 
On en d{\'e}duit que 
tout {\'e}l{\'e}ment
$x \in \atplus_F$ s'{\'e}crit de mani{\`e}re unique sous la forme:
$x=\sum_{i \in I} [\eps^i] a_i(x)$, o{\`u} $(a_i(x))_i$
est une suite de $\aplus_F$ qui tend vers $0$. 
On pose $R_k(x)= \sum_{i \in 
p^{-k} \ZZ \cap I} [\eps^i] a_i(x)$. 
En particulier, $R_k(\atplus_F) \subset \phi^{-k}(\aplus_F)$.
Rappelons que si $r$ est un entier $\geq 2$, alors  
$\btdag{,r}_F=\btplus_F\{p/\pi^r\}$, et 
que pour d{\'e}finir les $R_k$ sur $\btdag{,r}_F$ dans \cite[III.2]{CC98},
on {\'e}tend $R_k$ {\`a} $\btdag{,r}_F$ par la formule
$R_r(\sum a_i (p/\pi^r)^i) = \sum R_k(a_i) (p/\pi^r)^i$.
Pour construire les $R_k$ sur $\btrig{,r}{,K}$, et d{\'e}montrer 
la proposition \ref{rm}, on va
suivre un chemin semblable.

\begin{lemma}\label{riginter}
Soit $r$ un entier $\geq 2$.
On a $\btrig{,r}{,F}=\cap_{s \geq r} \btplus_F\{p/\pi^r,\pi^s/p\}$.
\end{lemma}

\begin{proof}
On va montrer que $\btrig{,r}{,F}=\cap_{s \geq r} 
\btplus_F\{p/[\pibar^r],[\pibar^s]/p\}$, le lemme en suit car si
$r>1$, alors $\pi/[\pibar]$ est une unit{\'e} de $\atdag{,r}$.

Soit $x \in \btrig{,r}{,F}$, et $n \gg 0$. On peut {\'e}crire $x=a+b$,
avec $a \in \btdag{,r}_K$, et $b \in (\bt_{[0;r_n]})^{H_K}$ par le
lemme \ref{decomp}. Le lemme \ref{bmaxhk} montre que $b$ peut
s'{\'e}crire $b=\sum b_i (\phi^n(\omega)/p)^i$, o{\`u} $b_i$ est une suite
de $\btplus_F$ qui tend vers $0$. On peut d'ailleurs remplacer $\omega$
par n'importe quel {\'e}l{\'e}ment $u$ de $\atplus_F$, tel que les id{\'e}aux 
$(\omega,p)$ et $(u,p)$ de $\atplus_F$ soient {\'e}gaux. On peut
notamment {\'e}crire $x=\sum c_i ([\pibar]^{r_n}/p)^i$, o{\`u} $c_i$ est une suite
de $\btplus_F$ qui tend vers $0$. Ceci {\'e}tant vrai pour tout $n \gg
0$, le lemme en r{\'e}sulte.
\qed\end{proof}

\begin{lemma}
Si on {\'e}tend $R_k$ {\`a} $\btplus_F \{ p/\pi^r, \pi^s/p \}$,
par \[ R_k \left(\sum_{i,j \geq 0} a_{i,j}   \left(\frac{p}{\pi^r}\right)^i 
\left(\frac{\pi^s}{p}\right)^j\right)=
 \left(\sum_{i,j \geq 0} R_k(a_{i,j})   \left(\frac{p}{\pi^r}\right)^i 
\left(\frac{\pi^s}{p}\right)^j\right) \]
alors, si $x \in  \btplus_F \{ p/\pi^r, \pi^s/p \}$, on a 
$\lim_{k \ra + \infty}R_k(x)=x$.
\end{lemma}

\begin{proof}
Cela suit du fait que $\lim_{k \ra + \infty}R_k(a_{i,j})=a_{i,j}$,
d'une part, et que $a_{i,j} \ra 0$, d'autre part.
\qed\end{proof}

\begin{proof}[de la proposition \ref{rm}]
Comme on a d{\'e}fini les applications $R_k$ de $\btplus_F \{ p/\pi^r, \pi^s/p
\}$ dans lui-m{\^e}me, et que $\btrig{,r}{,F} = \cap_{s \geq r} 
\btplus_F \{ p/\pi^r, \pi^s/p \}$ (par le lemme \ref{riginter}), 
on en d{\'e}duit des
applications $R_k : \btrig{,r}{,F} \ra \btrig{,r}{,F}$.
Les deux premiers points de la proposition \ref{rm} 
sont clairs sur la d{\'e}finition
pour $K=F$. Le troisi{\`e}me suit (toujours si $K=F$) du lemme pr{\'e}c{\'e}dent.

Cela d{\'e}finit $R_k$ dans le cas o{\`u} $K=F$. Dans le cas g{\'e}n{\'e}ral,
$\bdag{,r}_K / \bdag{,r}_F$ est une extension de degr{\'e} $e_K$ si
$r$ est assez grand. 
Soient $T_{K/F}=\sum_{\sigma \in H_F / H_K} \sigma$,
$\{e_i\}$ une base de $\bdag{,r}_K$  sur
$\bdag{,r}_F$, et $e_i^*$ la base duale pour la forme lin{\'e}aire
$(x,y) \mapsto T_{K/F}(xy)$. 
Alors, si $x \in \btrig{,r}{,K}$, on a $x=\sum
T_{K/F}(xe_i^*)e_i=\sum a_i e_i$
avec $a_i \in \btrig{,r}{,F}$, et on pose $R_k(x)=\sum R_k(a_i)e_i$.
Comme $R_k$ est $\bdag{,r}_F$-lin{\'e}aire, cette d{\'e}finition ne
d{\'e}pend pas du choix de la base $e_i$.

Les deux premiers points de la proposition sont clairs sur la
d{\'e}finition. 
On a d{\'e}j{\`a} montr{\'e} le troisi{\`e}me, dans le cas o{\`u}
$K=F$, et dans le cas g{\'e}n{\'e}ral, si on {\'e}crit $x= \sum_i a_i e_i$,
alors $R_k(x) \ra x$ car $R_k(a_i) \ra a_i$ pour chaque $i$.

Reste {\`a} voir que $R_k$ commute {\`a} $\gamma \in \Gamma_K$. C'est vrai
sur $\atplus_F$, car l'{\'e}criture que l'on utilise pour d{\'e}finir $R_k$
est unique. La mani{\`e}re dont on a construit $R_k$ montre qu'il en est
de m{\^e}me sur $\btrig{,r}{,F}$. Comme $\gamma(e_i)=\sum g_{i,j}
e_j$, avec $g_{i,j} \in \bdag{,r}_F$, on a
$R_k(\gamma(e_i))=\gamma(e_i)=\gamma(R_k(e_i))$: 
$R_k$ commute donc {\`a} $\gamma \in \Gamma_K$.
\qed\end{proof}

On d{\'e}finit aussi $\bnst{}{,K}=\bnrig{}{,K}[\log(\pi)]$\index{BdaglogK@$\bnst{}{,K}$}, 
cet anneau
est stable par les actions de $\phi$ et de $\Gamma_K$ {\'e}tant donn{\'e} que 
$\phi(\log(\pi))=\log(\phi(\pi))=p\log(\pi)+\log(\phi(\pi)/\pi^p)$
et $\gamma(\log(\pi))=\log(\pi)+\log(\gamma(\pi)/\pi)$
et que les s{\'e}ries qui d{\'e}finissent $\log(\phi(\pi)/\pi^p)$ et
$\log(\gamma(\pi)/\pi)$ convergent dans $\bnrig{}{,K}$.

Remarquons que si $n \geq 1$, alors $\iota_n(\bnst{,r_n}{,K}) \subset K_n[[t]]$.

\begin{definition}\label{rmlog}\index{Rk@$R_k$}
On prolonge les $R_k$ en une section
$\phi^{-k}(\bnrig{,p^k r}{,K})$-lin{\'e}aire
de l'inclusion de $\phi^{-k}(\bnst{,p^k r}{,K}[1/t])$ dans
$\btst{,r}{,K}[1/t]$: ils commutent toujours {\`a} $\Gamma_K$, et 
$\lim_{k \ra + \infty}R_k(x)=x$ si $x\in \btst{,r}{,K}[1/t]$.
\end{definition}

\section{Application aux repr{\'e}sentations $p$-adiques}
Ce chapitre est consacr{\'e} {\`a} l'application des constructions de la
section pr{\'e}c{\'e}dente {\`a} la caract{\'e}risation des repr{\'e}sentations
$p$-adiques semi-stables (et cristallines) en termes du
$(\phi,\Gamma_K)$-module qui leur est associ{\'e}. On montre 
en particulier que
si $V$ est une repr{\'e}sentation $p$-adique, alors 
\begin{align*}
 \dst(V) &= (\bnst{}{,K} \otimes_{\bdag{}_K}
\ddag{}(V)[1/t])^{\Gamma_K} 
\ \text{ et } \\
\dcris(V) &= (\bnrig{}{,K} \otimes_{\bdag{}_K} \ddag{}(V)[1/t])^{\Gamma_K}
\end{align*}

\subsection{R{\'e}gularisation par le Frobenius}
On commence par
{\'e}tablir un r{\'e}sultat de r{\'e}gularisation par le Frobenius, qui est
{\`a} la base des applications suivantes.
Le Frobenius $\phi$ est une bijection de $\bt_I$ sur $\bt_{pI}$ et
induit donc une bijection de $\btrig{,r}{}$ sur $\btrig{,pr}{}$ 
ainsi que de $\btst{,r}{}$ sur $\btst{,pr}{}$
puisque $\phi(\logpi)=p\cdot\logpi$.

\begin{lemma}\label{regul}
Soit $h$ un entier positif. Alors
\[ \cap_{s=0}^{+\infty} p^{-hs} \atdag{,p^{-s}r} = \atplus 
\text{ et }
\cap_{s=0}^{+\infty} p^{-hs} \atrig{,p^{-s}r}{} \subset \btrigplus{} \]
\end{lemma}

\begin{proof}
Montrons le premier point: comme $x \in \atdag{,r}$ il s'{\'e}crit de
mani{\`e}re unique sous la forme $\sum_{k \geq 0} p^k [x_k]$ et de
m{\^e}me 
$p^{hs}x=\sum p^{k+hs} [x_k]$. 
Comme $p^{hs}x \in
\atdag{,p^{-s}r}$ c'est que
\[ \vale(x_k)+\frac{rp^{1-s}}{p-1}(k+hs) \geq 0 \]
ce qui implique que
\[ \vale(x_k) \geq - \frac{(k+hs)r}{p^{s-1}(p-1)} \]
et donc (en laissant tendre $s$ vers $+\infty$) que $\vale(x_k) \geq 0$
ce qui fait que $x \in \atplus$.

Passons au deuxi{\`e}me point. Pour tout $s$ on peut {\'e}crire $x=a_s+b_s$
avec $a_s \in p^{-hs}\atdag{,p^{-s}r}$ et $b_s \in \btrigplus{}$.
Par le lemme \ref{unicite} on a $a_s-a_{s+1} \in \btplus$ et d'autre part
$a_s-a_{s+1} \in p^{-h(s+1)}\atdag{,p^{-s}r}$ ce qui fait que
$a_s-a_{s+1} \in p^{-h(s+1)}\atplus$ et que quitte {\`a} modifier $a_{s+1}$
on peut supposer que $a_s=a_{s+1}=a$. On a alors
$a \in \cap_{s=0}^{+\infty} p^{-hs} \atdag{,p^{-s}r} = \atplus$ et donc 
$x \in \btrigplus{}$. 
\qed\end{proof}

\begin{proposition}[R{\'e}gularisation par le Frobenius]\label{phireg}
Soit $r$ et $u$ deux entiers po\-sitifs et
$A \in \on{M}_{u \times r}(\btst{}{})$. On suppose qu'il existe $P \in
\on{GL}_u(F)$ telle que $A=P \phi^{-1}(A)$. Alors
$A \in \on{M}_{u \times r}(\btstplus{})$.
\end{proposition}

\begin{proof}
Soit $A=(a_{ij})$ et $a_{ij}=\sum_{n=0}^d a_{ij,n} \logpi^n$. Soit $h_0 \in
\ZZ$ tel que $p^{h_0} P \in \on{M}_u(\OO_F)$ et $h=h_0+d$. L'hypoth{\`e}se reliant
$A$ et $P$ peut s'{\'e}crire:
\[ p_{i1} \phi^{-1}(a_{1j})+\cdots+p_{iu} \phi^{-1}(a_{uj})=a_{ij}\qquad\forall
i \leq u,\ j \leq r \]
et comme $\phi^{-1}(\logpi^n)=p^{-n} \logpi^n$, on en d{\'e}duit que
si $a_{ij,n} \in p^{-c} \atrig{,r}{}$, alors comme $p^{h_0} p_{ik} \in \OO_F$
et $\phi^{-1}(a_{ik,n}) \in p^{-c} \atrig{,r/p}{}$, on a $a_{ij,n} \in
p^{-h-c} \atrig{,r/p}{}$. On it{\`e}re ce proc{\'e}d{\'e} et il en sort que
$a_{ij,n} \in \cap_{s=0}^{+\infty} p^{-hs-c} \atrig{,rp^{-s}}{}$. 
On est en mesure
d'appliquer le lemme \ref{regul} {\`a} $p^c a_{ij,n}$ et la proposition suit.
\qed\end{proof}

\subsection{Repr{\'e}sentations semi-stables}
Soient $\bnst{,r}{,K}=\bnrig{,r}{,K}[\log(\pi)]$ et
\[ \dnrig{}(V)= \bnrig{}{,K} \otimes_{\bdag{}_K} \ddag{}(V) \ \text{
et }\ \dnst{}(V)= \bnst{}{,K} \otimes_{\bdag{}_K} \ddag{}(V) \]
\index{Ddagrig@$\dnrig{}(V)$}\index{Ddaglog@$\dnst{}(V)$}
Le th{\'e}or{\`e}me \ref{toutsur} montre que $\dnrig{}(V)$ et $\dnst{}(V)$
sont des $\bnrig{}{,K}$- et $\bnst{}{,K}$- modules libres de rang 
$d=\dim_{\Qp}(V)$. Si $M$ est un $G_F$-module soit
$M(i)$\index{Mi@$M(i)$} le tordu de
$M$ par $\chi^i$ (twist de Tate).

\begin{proposition}
On a \[ \{ x \in \btst{}{},\ g(x)=\chi^i(g)x,\ \forall g \in G_K \} = 
\begin{cases} \text{$Ft^i$ si $i \geq 0$;} \\
\text{$0$ si $i < 0$.} \end{cases} \]
\end{proposition}

\begin{proof}
Soit $V^n_i= (\btst{,r_n}{}(i))^{G_K}$. C'est un $F$-espace vectoriel de
dimension finie (puisque $\iota_n$ r{\'e}alise une injection 
de $V^n_i$ dans $(\bdR^+)^{G_K}=K$)
stable par Frobenius et 
la proposition \ref{phireg} implique
que $V^n_i = (\btstplus{}(i))^{G_K}$ ce qui
fait \cite{Bu88per} que  $V^n_i=F$. Comme $V_i=\cup_{n=0}^{+\infty} V^n_i$ cela
d{\'e}montre le r{\'e}sultat.
\qed\end{proof}

Soit $\dst^+(V)=(\bst^+ \otimes_{\Qp} V)^{G_K}$\index{Dstplus@$\dst^+(V)$}, 
rappelons
que $\dst^+(V)=(\btstplus{} \otimes_{\Qp} V)^{G_K}$. 
Si $V$ a ses poids
de Hodge-Tate n{\'e}gatifs, alors $\dst^+(V)=\dst(V)$ et en g{\'e}neral
$\dst(V)=t^{-d} \dst^+(V(-d))$ pour $d$ assez grand.

\begin{proposition}\label{comp1}
Si $V$ est une repr{\'e}sentation $p$-adique,
alors $(\btst{}{} \otimes_{\Qp} V)^{G_K}$ est un $F$-espace vectoriel
de dimension finie, et
le morphisme induit par 
l'inclusion de $\btstplus{}$ dans $\btst{}{}$ 
\[ \dst^+(V) \ra (\btst{}{} \otimes_{\Qp} V)^{G_K} \] est un
isomorphisme de $(\phi,N)$-modules.
\end{proposition}

\begin{proof}
Si $n \in \NN$, alors $D_n=(\btst{,r_n}{} \otimes_{\Qp} V)^{G_K}$ 
est un $F$-espace vectoriel
de dimension finie $\leq [K:F]d$, 
car $\iota_n$ r{\'e}alise une injection de $D_n$ dans $\ddR(V)$, qui est un
$K$-espace vectoriel de dimension finie $\leq d$. 
Si l'on prend 
$[K:F]d+1$ {\'e}l{\'e}ments de $(\btst{}{} \otimes_{\Qp} V)^{G_K}$, 
ils vivent dans $D_n$ pour $n \gg 0$, et 
v{\'e}rifient donc une relation de d{\'e}pendance $F$-lin{\'e}aire. C'est donc que
$(\btst{}{} \otimes_{\Qp} V)^{G_K}$ est un $F$-espace vectoriel de
dimension $\leq [K:F]d$.

Passons maintenant au deuxi{\`e}me point.
Soient $v_1,\cdots,v_r$ et $d_1,\cdots,d_u$ des $\Qp$- et $F$- bases
de $V$ et $(\btst{}{} \otimes_{\Qp} V)^{G_K}$.
Il existe une matrice $A \in \on{M}_{r \times u}
(\btst{}{})$ telle que $(d_i)=A(v_i)$ (les $(d_i)$ et $(v_i)$ sont
des vecteurs colonnes). Soit $P \in \on{GL}_u(F)$ 
la matrice de $\phi$ dans la base $(d_i)$
(qui est inversible car $\phi: \btst{}{} \ra \btst{}{}$ 
est une bijection).
On a alors $\phi(A)=PA$ et donc $A=\phi^{-1}(P) \phi^{-1}(A)$; la proposition
\ref{phireg} montre que $A \in \on{M}_{r \times u}(\btstplus{})$
et donc que $(\btst{}{} \otimes_{\Qp}V)^{G_K} \subset
(\btstplus{} \otimes_{\Qp}V)^{G_K}= \dst^+(V)$, ce qui permet de conclure.
\qed\end{proof}

Une repr{\'e}sentation $V$ {\`a} poids n{\'e}gatifs
est donc semi-stable
si et seulement si
elle est $\btst{}{}$-admissible et elle est cristalline
si et seulement si
elle est $\btrig{}{}$-admissible, c'est-{\`a}-dire
si elle est $\btst{}{}$-admissible et que ses
p{\'e}riodes sont tu{\'e}es par $N$.

De plus, si les poids de Hodge-Tate de $V$ ne sont pas n{\'e}gatifs,
alors en tordant $V$ on en d{\'e}duit que 
$V$ est semi-stable si et seulement si
elle est $\btst{}{}[1/t]$-admissible et elle est cristalline
si et seulement si
elle est $\btrig{}{}[1/t]$-admissible.

\begin{proposition}\label{comparaison}
Si $V$ est semi-stable
on a un isomorphisme de comparaison:
\[ \btst{}{}[1/t] \otimes_F \dst(V) = \btst{}{}[1/t] \otimes_{\Qp} V \]
\end{proposition}

\begin{proof}
Ceci r{\'e}sulte du fait que dans ce cas on a d{\'e}j{\`a}:
\[ \btstplus{}[1/t] \otimes_F \dst(V) = \btstplus{}[1/t]
\otimes_{\Qp} V \] il suffit alors de tensoriser les deux membres par 
$\btst{}{}[1/t]$ au-dessus de $\btstplus{}[1/t]$.
\qed\end{proof}

\begin{theorem}\label{isomcomp}
Si $V$ est une repr{\'e}sentation $p$-adique, alors 
\[ \dst(V)=(\dnst{}(V)[1/t])^{\Gamma_K}  \ \text{ et }\  
\dcris(V)=(\dnrig{}(V)[1/t])^{\Gamma_K} \]
Notamment $V$ est semi-stable
(respectivement cristalline) si et seulement si
$(\dnst{}(V)[1/t])^{\Gamma_K}$ (respectivement
$(\dnrig{}(V)[1/t])^{\Gamma_K}$) est un $F$-espace vectoriel de
dimension $d=\dim_{\Qp}(V)$.
\end{theorem}

\begin{proof}
Le deuxi{\`e}me point est une cons{\'e}quence imm{\'e}diate du
premier. Ensuite comme 
$\dnst{}(V)[1/t] \subset (\btst{}{}[1/t]\otimes_{\Qp}
V)$ (et que $\dnrig{}(V)[1/t] \subset (\btrig{}{}[1/t]\otimes_{\Qp}
V)$) les r{\'e}sultats pr{\'e}c{\'e}dents montrent que 
$(\dnst{}(V)[1/t])^{\Gamma_K}$ (respectivement
$(\dnrig{}(V)[1/t])^{\Gamma_K}$) est inclus dans $\dst(V)$
(respectivement dans $\dcris(V)$). 

Montrons donc que $\dst(V) \subset  (\dnst{}(V)[1/t])^{\Gamma_K}$, et que 
$\dcris(V) \subset (\dnrig{}(V)[1/t])^{\Gamma_K}$.
Il suffit de s'occuper du cas semi-stable
car le cas cristallin en suit en faisant $N=0$. Soit $r=\dim_F(\dst(V))$.
On peut supposer 
(quitte {\`a} tordre)
les poids de Hodge-Tate de $V$ n{\'e}gatifs
puisque l'on a invers{\'e} $t$ partout. 
On sait qu'alors
$\dst(V)=(\btst{}{} \otimes_{\Qp} V)^{G_K}$
et de plus 
$(\btst{}{} \otimes _{\Qp}V)^{H_K}  
=\btst{}{,K} \otimes_{\bdag{}_K} \ddag{}(V)$
puisque $\ddag{}(V)$ a la bonne dimension.
On en d{\'e}duit 
que si l'on choisit une base $\{e_i\}$ de $\ddag{}(V)$ 
et $\{d_i\}$ une base de $\dst(V)$, alors la matrice $M \in 
\on{M}_{r \times d}(\btst{}{,K})$ d{\'e}finie par $(d_i)=M(e_i)$
est de rang $r$ et v{\'e}rifie
$\gamma_K(M)G-M=0$ o{\`u} $G \in \on{GL}_d(\bdag{}_K)$ est
la matrice de $\gamma_K$ dans la base $\{e_i\}$.

Les op{\'e}rateurs $R_m$ introduits au paragraphe \ref{rmdef} sont 
$\bnst{}{,K}$-lin{\'e}aires et commutent {\`a} 
$\Gamma_K$ (voir la proposition \ref{rm}
et la d{\'e}finition \ref{rmlog}) ce qui fait
que $\gamma_K( R_m(M) )G-R_m(M)=0$.
De plus, $R_m(M) \ra M$ et si $M \in \on{M}_{r \times d}(
\btst{,r_n}{,K})$,
alors $R_m(M) \in \on{M}_{r \times d}(\btst{,r_n}{,K})$.
Soit $N=\phi^m(R_m(M))$. On a alors $\gamma_K(N)\phi^m(G)=N$ et comme les
actions de $\phi$ et $\Gamma_K$ commutent sur $\dnrig{}(V)$ on a
$\phi(G)=\gamma_K(P)GP^{-1}$ ($P$ est la matrice de $\phi$ et est 
inversible car $\phi$ est surconvergent et $\bdag{}_K$ est un corps) 
ce qui fait que si 
$Q=\phi^{m-1}(P)\cdots\phi(P)P$, alors $\phi^m(G)=\gamma_K(Q)GQ^{-1}$
et donc $\gamma_K(NQ)G=(NQ)$. 
La matrice $NQ$ d{\'e}termine $r$ {\'e}l{\'e}ments de $\dnst{}(V)$ qui sont 
fix{\'e}s par $\gamma_K$. Il reste {\`a} montrer que ces {\'e}l{\'e}ments sont 
libres sur $F$ quand $m$ est assez grand.
Mais comme $R_m(M) \ra M$, la matrice $NQ$ va  
{\^e}tre de rang $r$
pour $m$ assez grand (puisque $M$ est de rang $r$)
et donc d{\'e}terminer un sous-module libre de rang $r$ de
$\dnst{}(V)$. A fortiori le $F$-espace vectoriel engendr{\'e} par les
{\'e}l{\'e}ments d{\'e}termin{\'e}s par $NQ$ va {\^e}tre de dimension $r$
et donc {\'e}gal {\`a} $\dst(V)$.
\qed\end{proof}

\begin{proposition}\label{isomcomp2}
On a les isomorphismes de comparaison suivants:
\begin{enumerate}
\item si $V$ est une repr{\'e}sentation semi-stable, alors \[ \ddag{}(V)
  \otimes_{\bdag{}_K}
\bnst{}{,K}[1/t] = \dst(V) \otimes_F \bnst{}{,K}[1/t] \]
\item si $V$ est une repr{\'e}sentation cristalline, alors \[ \ddag{}(V) 
\otimes_{\bdag{}_K}
\bnrig{}{,K}[1/t] = \dcris(V) \otimes_F \bnrig{}{,K}[1/t] \]
\end{enumerate}
Si de plus $V$ a ses poids de Hodge-Tate n{\'e}gatifs, alors
$\dst(V) \subset \bnst{}{,K} \otimes_{\bdag{}_K}  \ddag{}(V)$.
\end{proposition}

\begin{proof}
Encore une fois on ne s'occupe que du cas semi-stable, le cas
cristallin s'obtenant en faisant $N=0$.
On peut supposer que $V$ a ses poids de 
Hodge-Tate n{\'e}gatifs car cela revient 
{\`a} multiplier par une puissance de $t$ le terme de gauche. On sait
qu'alors $\dst(V) \subset \btst{}{,K}\otimes_{\bdag{}_K} \ddag{}(V)$ et que
\[ \btst{}{,K}[1/t] \otimes_{\bdag{}_K} \ddag{}(V)=
\btst{}{,K}[1/t] \otimes_F \dst(V) \]
ce qui montre que si l'on choisit des bases $\{ d_i \}$ de
$\dst(V)$ et $\{e_i\}$ de $\ddag{}(V)$, alors 
$(e_i)=B(d_i)$ avec $B \in \on{M}_d( \btst{}{,K}[1/t])$; la
proposition \ref{isomcomp} implique d'autre part que
$(d_i)=A(e_i)$ avec 
$A \in \on{M}_d(\bnst{}{,K}[1/t])$; de plus $AB=\on{Id}$. On
peut alors appliquer l'op{\'e}rateur $R_0$ qui est
$\bnst{}{,K}[1/t]$-lin{\'e}aire pour trouver  $AR_0(B)=\on{Id}$ ce qui
fait que $B=R_0(B)$ et que $B$ a donc ses coefficients dans 
$\bnst{}{,K}[1/t]$ et $A\in\on{GL}_d(\bnst{}{,K}[1/t])$. Ceci permet
de conclure.
\qed\end{proof}

\begin{remark}
Les coefficients d'une matrice de l'isomorphisme de comparaison 
sont tr\`es li\'es \`a l'exponentielle de Perrin-Riou \cite{BP94,BP99,BP00}.
\end{remark}

\begin{proposition}\label{detstdag}
Soit $V$ une repr{\'e}sentation semi-stable et
soit $M$ la matrice de passage d'une base de $\dst(V)$ {\`a} une base de
$\ddag{}(V)$. Alors il existe $r \in \ZZ$ et $\lambda \in \bdag{}_K$ 
tels que $\det(M)=\lambda t^r$.
\end{proposition}

\begin{proof}
Le d{\'e}terminant de la matrice de passage est {\'e}gal au coefficient de
la matrice de passage pour le d{\'e}terminant de $V$ et il suffit donc
de montrer l'assertion en dimension $1$. Une repr{\'e}sentation semi-stable
de dimension $1$ est cristalline et est donc
de la forme $\omega \chi^r$ o{\`u} $\omega$ est un
caract{\`e}re non-ramifi{\'e} et $\chi$ est le caract{\`e}re
cyclotomique. La p{\'e}riode de $\omega$ est alors un {\'e}l{\'e}ment
$\beta\in W(\overline{k})$, 
ce qui fait que $\dst(V)=F\cdot \beta t^{-r}$ et
$\ddag{}(V)=\bdag{}_K \cdot \beta$
d'o{\`u} le r{\'e}sultat.
\qed\end{proof}

\subsection{Repr{\'e}sentations cristallines et repr{\'e}sentations de hauteur 
finie}
Dans ce paragraphe, on se place dans le cas $K=F$ et $V$ est une
repr{\'e}sentation cristalline de $G_F$.
On dit qu'une repr{\'e}sentation $p$-adique $V$ de $G_F$ 
est de hauteur finie si
$\dfont(V)$ poss{\`e}de une base sur $\mathbf{B}_F$ form{\'e}e
d'{\'e}l{\'e}ments de $\dfont^+(V) = (\bplus \otimes_{\Qp} V)^{H_F}$. Un
r{\'e}sultat de Fontaine 
\cite{Fo91} (voir aussi \cite[III.2]{Co99}) montre que
$V$ est de hauteur finie si et
seulement si $\dfont(V)$ poss{\`e}de un sous-$\bplus_F$-module libre de
type fini stable par $\phi$
de rang {\'e}gal {\`a} $d=\dim_{\Qp}(V)$.
L'objet de ce paragraphe est de d{\'e}montrer le r{\'e}sultat suivant:
\begin{theorem}\label{hf}
Si $V$ est une repr{\'e}sentation cristalline de $G_F$, 
alors $V$ est de hauteur finie.
\end{theorem}

\begin{proof}
Fixons une base $\{e_i\}$ de $\ddag{}(V)$ ainsi qu'une base 
$\{d_i\}$ de $\dcris(V)$, et soit $U$ la
matrice de passage de l'une {\`a} l'autre 
c'est-{\`a}-dire que $(e_i)=U(d_i)$. La matrice $U$
est {\`a} coefficients dans $\bnrig{}{,F}[1/t]$.
De plus si l'on remplace $V$ par $V(1)$, alors on peut remplacer 
$d_i$ par $t^{-1}d_i(1)$ et donc $U$ par $tU$ ce qui fait que quitte
{\`a} tordre suffisamment $V$ on peut supposer que $U$
est {\`a} coefficients dans $\bnrig{}{,F}$, ce que l'on fait maintenant.
Soit $\bhol{,F}$ l'anneau des s{\'e}ries formelles 
$\sum_{k
\geq 0} a_k \pi^k$ o{\`u} $a_k \in F$ et 
$\sum_{k
\geq 0} a_k X^k$ est
de rayon de convergence $1$
(c'est-\`a-dire qu'elles convergent sur le disque ouvert de rayon $1$).
Nous aurons besoin d'un r{\'e}sultat de Kedlaya \cite[5.3]{KK00}
(la formulation originale de Kedlaya
est $U=VW$ mais on s'y ram{\`e}ne en transposant):
\begin{proposition}
Si $U$ est une matrice {\`a} coefficients dans $\bnrig{}{,F}$, alors il existe
$V$ dans $\on{Id}+\pi\on{M}_d(\bhol{,F})$ et $W$ dans
$\on{M}_d(\bdag{}_F)$ telles que $U=WV$.
\end{proposition}

Remarquons qu'on a n{\'e}cessairement $\det(W) \neq 0$ et comme $\bdag{}_F$
est un corps cela implique que $W$ est inversible. 
Soit $P$ la matrice de $\phi$ dans la base
$\{e_i\}$ (qui a ses coefficients surconvergents),
$D$ la matrice de $\phi$ dans la base $\{d_i\}$ 
(qui est donc {\`a} coefficients dans 
$F$) et $G$ la matrice de $\gamma_F$ (un g{\'e}n{\'e}rateur de $\Gamma_F$) 
dans la base $\{e_i\}$
(elle est aussi surconvergente). 
Un petit calcul montre que, si $\{f_i\}$ est la base de $\ddag{}(V)$
d{\'e}duite de $\{e_i\}$ par $(f_i)=W^{-1}(e_i)$, alors:  
\begin{align*} 
\on{Mat}_{\{f_i\}} (\phi) &= \phi(V)DV^{-1} = \phi(W^{-1})PW \\
\on{Mat}_{\{f_i\}} (\gamma_F) &= \gamma_F(V)V^{-1} = \gamma_F(W^{-1})GW 
\end{align*}
Les coefficients de la
matrice $\phi(V)DV^{-1}$ sont dans
$\on{Frac}\bhol{,F}$.
D'autre part les coefficients de $\phi(W^{-1})PW$ 
sont dans $\bdag{}_F$.
On en d{\'e}duit que
dans la base $\{f_i\}$, la matrice de $\phi$  
a ses coefficients dans $\on{Frac}\bhol{,F}$ 
d'une part et dans $\bdag{}_F$ d'autre part. 
C'est donc \cite[II.12]{Co99} que $\on{Mat}_{\{f_i\}} (\phi)
\in \on{M}_d(\on{Frac}\bplus_F)$.

De plus si $v=\det(V)$, alors
$\det(\on{Mat}_{\{f_i\}} (\phi))=\phi(v)\det(D)v^{-1}$ et comme 
$V \in \on{Id}+\pi\on{M}_d(\bhol{,F})$ on a $v \in  1+\pi \bhol{,F}$.
Les coefficients de $\on{Mat}_{\{f_i\}} (\phi)$ n'ont donc pas de
p{\^o}les en z{\'e}ro ce qui fait
qu'il existe $\lambda \in \bplus_F$ non-divisible par $\pi$
tel que $\lambda \on{Mat}_{\{f_i\}} (\phi) \in \on{M}_d(\bplus_F)$.
Pour les m{\^e}me raisons, $\on{Mat}_{\{f_i\}} (\gamma_F) \in
\on{M}_d(\on{Frac}\bplus_F)$.

Soit $D$ le $\bplus_F$-module engendr{\'e} par les $f_i$. On vient de
voir que $\lambda \phi(D) \subset D$ o{\`u}
$\lambda \in \bplus_F$ n'est pas divisible par $\pi$
et que le  $\on{Frac}(\bplus_F)$-module engendr{\'e} 
par $D$ est stable par $\Gamma_F$.
\begin{lemma}
Si $D$ est un $\bplus_F$-module libre de type fini 
tel que le $\on{Frac}(\bplus_F)$-module engendr{\'e} 
par $D$ est stable par $\Gamma_F$
et $\lambda \in \bplus_F-\pi \bplus_F$ est tel que 
$\lambda \phi(D) \subset D$, alors il existe $D' \subset D$ un
sous-module libre de type fini stable par $\phi$ et $\Gamma_F$
qui est de rang maximal.
\end{lemma}
Pour montrer que la repr{\'e}sentation $V$ est de hauteur finie il
suffit donc, gr{\^a}ce au r{\'e}sultat de Fontaine, 
de montrer le lemme ce que nous faisons maintenant.
Le lemme est d{\'e}montr{\'e} dans \cite{Co99} 
(c'est la r{\'e}union des {\'e}nonc{\'e}s III.8 {\`a} III.15)
mais l'une des {\'e}tapes 
utilise de mani{\`e}re cruciale que $k$ est fini et il nous faut la
contourner
\footnote{Il s'agit du lemme III.9. On remarquera d'ailleurs que ce
lemme aurait plut{\^o}t d{\^u} {\^e}tre {\'e}nonc{\'e} de la mani{\`e}re suivante:
``\textit{Soit $M$ un $\bplus_K$-module libre de rang
fini muni d'une action de $\Gamma$ tel que le
$\on{Frac}\bplus_K$-module $M\otimes_{\bplus_K}\on{Frac}\bplus_K$
soit muni d'une action de $\phi$ commutant {\`a} celle de $\Gamma$ et
telle qu'il existe $a\in\NN$ tel que l'on ait $\phi(M)\subset
\pi^{-a}M$. Alors $\phi(M)\subset M$.}''
C'est cet {\'e}nonc{\'e} qui est d{\'e}montr{\'e} et utilis{\'e}  
dans la suite de \cite{Co99}.}.

On va d'abord montrer que l'on peut supposer que $\Gamma_F(D) \subset D$.
Soit $G=\pscal{\gamma_F}$;
on renvoie {\`a} la cons{\'e}quence du corollaire III.15 de
\cite{Co99} pour la construction, sous
l'hypoth{\`e}se que le $\on{Frac}(\bplus_F)$-module engendr{\'e} 
par $D$ est stable par $\Gamma_F$, de $\alpha,\beta \in \bplus_F$ tels
que pour tout $\gamma\in G$ on ait $\gamma(\alpha D) \subset
\beta D$. Comme $\gamma(\pi)/\pi$ est inversible dans $\bplus_F$, on
peut supposer que $\pi$ ne divise pas $\alpha$; en effet, l'inclusion 
$\alpha D \subset \beta D$ implique que si $\pi$ divise $\alpha$,
alors $\pi$ divise $\beta$.
Soit alors $G \cdot \alpha D$ le
sous-$\bplus_F$-module de $\beta D$ engendr{\'e} par les $\gamma(\alpha d)$
o{\`u} $\gamma\in G$, $d\in D$. Comme $\bplus_F$ est noetherien (car principal), 
comme $\beta D$ est de
type fini, et comme $G \cdot \alpha D$ est r{\'e}union croissante de 
sous-modules de type
fini de $\beta D$, 
c'est que $G \cdot \alpha D$ est en fait r{\'e}union d'un nombre fini de
$\gamma(\alpha D)$. 
Il existe donc $n\in\NN$ et $\gamma_1,\cdots,\gamma_n \in G$
tels que $G \cdot \alpha D = \sum \gamma_i(\alpha D)$.
De plus $G \cdot \alpha D$ est stable par $G$ et
donc par continuit{\'e} il est stable par $\Gamma_F$.
Ensuite comme $\lambda \phi(D) \subset D$ on a
$\lambda\alpha \cdot \phi(\alpha D) \subset \alpha D$. Soit
$\mu=\prod_{i=1}^n \gamma_i(\lambda\alpha)$, on voit que
$\mu \phi(\gamma_i(\alpha D)) \subset \gamma_i(\alpha D)$ pour tout
$i$ et donc que $\mu \phi(G \cdot \alpha D) \subset G \cdot \alpha D$
et que $\pi$ ne divise pas $\mu$. Ceci montre que l'on peut se ramener
au cas o{\`u} $D$ est stable par $\Gamma_F$.

On suppose donc maintenant que 
$D$ est un $\bplus_F$-module libre de type fini 
stable par $\Gamma_F$
et que $\lambda \in \bplus_F-\pi \bplus_F$ est tel que 
$\lambda \phi(D) \subset D$.
L'id{\'e}al $I$ des $\delta\in\bplus_F$ tels que $\delta \phi(D) \subset
D$ est stable par $\Gamma_F$ et par \cite[III.8]{Co99} il est de la
forme $(\pi^a \prod_{i=0}^n \phi^i(q)^{\beta_i})$. De plus comme
$\lambda\in I$ et que $\pi$ ne divise pas $\lambda$ c'est que
$a=0$.
Soit alors $\alpha =\pi^{\beta_0+\cdots+\beta_n}
q^{\beta_1+\cdots+\beta_n} \cdots \phi^{n-1}(q)^{\beta_n}$. Un petit
calcul montre que $\phi(\alpha D) \subset \alpha D$
et que $\alpha D$ est stable par $\Gamma_F$. 
On peut donc
prendre $D'=\alpha D$ et ceci ach{\`e}ve la d{\'e}monstration du lemme et
donc du th{\'e}or{\`e}me \ref{hf}.
\qed\end{proof}

\subsection{Une autre d{\'e}finition de $\dnrig{}(V)$}\index{Bdagrig@$\bnrig{}{}$}
L'objet de ce paragraphe est de montrer qu'il existe un anneau  $\bnrig{}{}$,
tel que $\dnrig{}(V) = (\bnrig{}{}
\otimes_{\Qp} V)^{H_K}$.  On pose $\bnrig{}{} = \bnrig{}{,F}
\otimes_{\bdag{}_F} \bdag{}$. Le reste de ce paragraphe est consacr{\'e}
{\`a} la preuve de quelques propri{\'e}t{\'e}s de $\bnrig{}{}$, et  du
fait que $\dnrig{}(V) = (\bnrig{}{}
\otimes_{\Qp} V)^{H_K}$. Ces r{\'e}sultats ne sont 
pas utilis{\'e}s dans la suite de l'article.
\begin{lemma}
On a  $\bnrig{}{,K} = \bnrig{}{,F}
\otimes_{\bdag{}_F} \bdag{}_K$.
\end{lemma}

\begin{proof}
Voir la d{\'e}monstration du lemme 4.1 de \cite{An00}.
\qed\end{proof}

\begin{corollary}
On a $\bnrig{}{} = \bnrig{}{,K}
\otimes_{\bdag{}_K} \bdag{}$.
\end{corollary}

\begin{proposition}
On a $(\bnrig{}{})^{H_K} = \bnrig{}{,K}$.
\end{proposition}

\begin{proof}
Comme $\bdag{}_F$ est un corps, il existe une base $\{e_{\alpha}\}$ de 
$\bnrig{}{,F}$ sur $\bdag{}_F$. Soit $x \in \bnrig{}{}$, il s'{\'e}crit
de mani{\`e}re unique $x = \sum  \lambda_{\alpha} e_{\alpha}$, avec
$\lambda_{\alpha} \in \bdag{}$. Si $x$ est fixe par $H_K$, c'est donc
(comme les $e_{\alpha}$ sont fixes par $H_K$)
que $\lambda_{\alpha} \in (\bdag{})^{H_K}=\bdag{}_K$. Le lemme en r{\'e}sulte.
\qed\end{proof}

\begin{lemma}
L'application naturelle $\bnrig{}{,F}
\otimes_{\bdag{}_F} \bdag{}$ dans $\btrig{}{}$ est injective.
\end{lemma}

\begin{proof}
Soit $k \geq 1$, ayant la propri{\'e}t{\'e} suivante: c'est le plus petit
entier tel qu'il existe
$\lambda_1,\cdots,\lambda_k \in
\bnrig{}{,F}$, et $b_1,\cdots,b_k \in \bdag{}$, tels que $\sum_i
\lambda_i b_i = 0$ dans $\btrig{}{}$, avec $\sum \lambda_i \otimes b_i
\neq 0$.
 
Alors, si $h \in H_F$, on a $\sum_i \lambda_i h(b_i) = 0$, et donc,
$\sum_i \lambda_i (h(b_i/b_1)-b_i/b_1) = 0$. Par l'hypoth{\`e}se de
minimalit{\'e} sur $k$, $h(b_i/b_1)-b_i/b_1 = 0$, et donc
$b_i/b_1 \in \bdag{}_F$.
On a alors $\sum_i
\lambda_i b_i/b_1 =0$, et donc finalement, $\sum_i
\lambda_i \otimes b_i 
= b_1 \otimes (\sum_i
\lambda_i b_i/b_1)
=0$. L'application
$\bnrig{}{,F} \otimes_{\bdag{}_F} \bdag{}$ dans $\btrig{}{}$
est donc injective.
\qed\end{proof}

\begin{proposition}
On a un isomorphisme $\dnrig{}(V) = (\bnrig{}{}
\otimes_{\Qp} V)^{H_K}$.
\end{proposition}

\begin{proof}
Comme $\bdag{} \otimes_{\Qp} V = \bdag{} \otimes_{\bdag{}_K}
\ddag{}(V)$, 
il suffit de montrer que $\dnrig{}(V) = (\bnrig{}{} 
\otimes_{\bdag{}_K} \ddag{}(V))^{H_K}$. Comme $\ddag{}(V)$ est
un $\bdag{}_K$-module libre de rang $d$, 
on a: \[ (\bnrig{}{} 
\otimes_{\bdag{}_K} \ddag{}(V))^{H_K}= \bnrig{}{,K} \otimes_{\bdag{}_K}
\ddag{}(V)
= \dnrig{}(V). \]
\qed\end{proof}

\begin{remark}\index{Bdaglog@$\bnst{}{}$}
Si l'on pose $\bnst{}{}=\bnrig{}{}[\log(\pi)]$, alors
$\dnst{}(V)=(\bnst{}{} \otimes_{\Qp} V)^{H_K}$.
\end{remark}

\section{Propri{\'e}t{\'e}s de $\bnrig{,r}{,K}$}
L'anneau $\bnrig{}{,K}$ est isomorphe 
(canoniquement si $K=F$) {\`a} l'anneau de Robba utilis{\'e}
dans la th{\'e}orie des {\'e}quations diff{\'e}rentielles $p$-adiques. Ce
chapitre est consacr{\'e} {\`a} la d{\'e}monstration de quelques unes de ses
propri{\'e}t{\'e}s  relatives {\`a} l'action de $\Gamma_K$: on 
d{\'e}finit aussi des op{\'e}rateurs diff{\'e}rentiels qui seront utiles 
pour la suite. Enfin, on montre un th{\'e}or{\`e}me de structure pour les
modules sur l'anneau de Robba.

\subsection{L'op{\'e}rateur $\nabla$}
Dans ce paragraphe, $\gamma$ est un {\'e}l{\'e}ment de $\Gamma_K$ et $n(\gamma)=
v_p(1-\chi(\gamma))$\index{ngamma@$n(\gamma)$}. 
On suppose que $n(\gamma)\geq 1$ et que
$r \geq r_{n(K)}$.

\begin{lemma}\label{unmoinsgampetit} 
Si $I=[r;s]$ est un intervalle ferm{\'e} de $[r;+\infty[$, alors
il existe $n(I) \in \NN$, tel que pour 
$x\in\bnrig{,r}{,K}$, on ait
$V_I((1-\gamma)x)\geq V_I(x)+1$ d{\`e}s que $n(\gamma) \geq n(I)$.
\end{lemma}

\begin{proof}
Par densit{\'e} et lin{\'e}arit{\'e} 
on se ram{\`e}ne au cas o{\`u} $x=\pi_K^k$ avec
$k\in \ZZ$. Alors si $k \geq 0$ 
\begin{align*}
\gamma(\pi_K^k)-\pi_K^k &= \pi_K^k\left(\frac{\gamma(\pi_K^k)}
{\pi_K^k}-1\right) \\
&= \pi_K^k\left(\frac{\gamma(\pi_K)}{\pi_K}-1\right)
\left(\frac{\gamma(\pi_K^{k-1})}{\pi_K^{k-1}}+\cdots+1\right) 
\end{align*}
et sinon
\begin{align*} 
\gamma(\pi_K^{-k})-\pi_K^{-k} &= \pi_K^{-k}\left(\frac
{\pi_K^k}{\gamma(\pi_K^k)}-1\right) \\
&= \pi_K^{-k}\left(\frac{\pi_K}{\gamma(\pi_K)}-1\right)
\left(\frac{\pi_K^{k-1}}{\gamma(\pi_K^{k-1})}+\cdots+1\right) 
\end{align*}
Comme on a $V_I(xy)\geq V_I(x)+V_I(y)$ (voir la remarque
\ref{valanno}), le lemme r{\'e}sulte du fait que
$V_I(\gamma(\pi_K)/\pi_K-1) \geq 1$ 
si $n(\gamma)$ est assez grand; en effet, $\gamma(\pi_K)/\pi_K-1$  
tend vers $0$ quand $\gamma$ tend vers $1$. 
\qed\end{proof}

Si $I=[r;s]$,
soit $\mathbf{B}^I_K$ l'anneau des s{\'e}ries formelles en $\pi_K$ qui
convergent sur la couronne de rayons int{\'e}rieurs et ext{\'e}rieurs
$\alpha(K,r)$ et $\alpha(K,s)$. 
Le lemme pr{\'e}c{\'e}dent montre que si $\gamma$ est assez proche de $1$,
la s{\'e}rie d'op{\'e}rateurs \[ \frac{\log(\gamma)}{\log(\chi(\gamma))}
=-\frac{1}{\log(\chi(\gamma))}\sum_{n\geq 1}\frac{(1-\gamma)^n}{n}\]
converge vers un op{\'e}rateur continu $\nabla_I : \bnrig{,r}{,K} \ra 
\mathbf{B}^I_K$, et un petit argument de s{\'e}ries formelles montre que 
$\log(\gamma)/\log(\chi(\gamma))$ ne d{\'e}pend pas du choix de
$\gamma$. Notamment $\nabla_{I_1}(x)=\nabla_{I_2}(x)$ si
$I_k=[r;s_k]$. On en d{\'e}duit que la valeur commune des $\nabla_I(x)$
appartient {\`a} $\bnrig{,r}{,K}$, et que l'op{\'e}rateur $x \mapsto
\nabla(x)$, o{\`u} $\nabla(x)$ est la valeur commune des $\nabla_I(x)$,
est continu pour la topologie de Fr{\'e}chet\index{nabla@$\nabla$}.

De m{\^e}me, si $\gamma$ est assez proche de $1$,
alors la s{\'e}rie d'op{\'e}rateurs \[ \frac{\log(\gamma)}{\log(\chi(\gamma))(1-\gamma)}
=-\frac{1}{\log(\chi(\gamma))}\sum_{n\geq
 1}\frac{(1-\gamma)^{n-1}}{n}\] 
converge vers un op{\'e}rateur continu $\nabla/(1-\gamma) : \bnrig{,r}{,K} \ra 
\mathbf{B}^I_K$, et un petit argument de s{\'e}ries formelles montre que 
\[ \frac{\nabla}{1-\gamma_K}=
\frac{1-\gamma}{1-\gamma_K}\cdot\frac{\nabla}{1-\gamma} \]
ne d{\'e}pend pas du choix de
$\gamma$ et d{\'e}finit aussi un op{\'e}rateur  de $\bnrig{,r}{,K}$ dans
lui-m{\^e}me, continu pour la topologie de Fr{\'e}chet. Il est clair que \[ 
(1-\gamma_K) \frac{\nabla}{1-\gamma_K} =\nabla. \]

\begin{lemma}
La restriction de $\nabla$ {\`a} $\bnrig{,r}{,K}$
v{\'e}rifie $\nabla=t\cdot\partial$ o{\`u} 
$\partial(x)=(1+\pi)dx/d\pi$\index{dronde@$\partial$}.
\end{lemma}

\begin{remark}
Attention au fait que ces notations ne sont pas compatibles avec  
\cite{CC99}. Ce qui est not{\'e} $\partial$ chez nous
est not{\'e} $\nabla$ dans \cite{CC99}. 
\end{remark}

\begin{proof}
Dans \cite{CC99} il est d{\'e}montr{\'e} que l'image de $\bdag{,r_n}_K$ par 
$\iota_n=\phi^{-n}$ dans $\bdR^+$ 
est contenue dans $K_n[[t]]$ si $n\geq n(K)$.
L'image par $\phi^{-n}$ de $\bnrig{,r_n}{,K}$ a la m{\^e}me propri{\'e}t{\'e}
par continuit{\'e}.
Comme $\phi^{-n}$ est injectif et que $\nabla$, $t\partial$ commutent {\`a}
$\phi$ il suffit de montrer que $\nabla-t\partial$ est nul sur $K_n[[t]]$
ce qui est {\'e}vident.
\qed\end{proof}

Montrons que $\partial$ est presque surjective:
\begin{proposition}\label{invconn}
La connexion $\partial$ r{\'e}alise une surjection 
\begin{enumerate}
\item de $\bnrig{}{,K}+F\cdot\log(\pi_K)$ dans $\bnrig{}{,K}$;
\item de $\bnst{}{,K}$ dans lui-m{\^e}me.
\end{enumerate}
\end{proposition}

\begin{proof}
Tout d'abord si $e_K$ est l'indice
de ramification de $K_{\infty}/F_{\infty}$, alors
$\log(\pi/\pi_K^{e_K})$ est une s{\'e}rie surconvergente
ce qui fait que $\bnst{}{,K}=\bnrig{}{,K}[\log(\pi_K)]$.
Rappelons que $\bdag{}_K$ est un corps, et on peut donc se fixer 
$r \geq r_{n(K)}$ tel que $\partial(\pi_K)$ et 
$1/\partial(\pi_K)$ appartiennent {\`a} $\bdag{,r}_K$.

Soit $F(\pi_K) \in \bnrig{,r}{,K}$. 
Posons $F(\pi_K)/\partial(\pi_K)=\sum a_n \pi_K^n$, alors
\begin{align*} 
F(\pi_K) &= \sum_{n \neq -1} \frac{a_n}{n+1} (n+1)\pi_K^n \partial \pi_K
+ a_{-1} \frac{\partial\pi_K}{\pi_K} \\
&= \partial\left(\sum_{n \neq -1} \frac{a_n}{n+1}\pi_K^{n+1}
+ a_{-1}\log(\pi_K)\right) 
\end{align*}
ce qui montre que $F(\pi_K)=\partial G(\pi_K)$ avec $G(\pi_K) \in
\bnrig{,s}{,K}+F \cdot \log(\pi_K)$ pour $s>r$.

Ensuite la formule 
\begin{multline*} 
\partial( G_j(\pi_K)\log^j(\pi_K) ) = \\ 
(\partial G_j)(\pi_K)\partial \pi_K \log^j(\pi_K)+
G_j(\pi_K) j \log^{j-1}(\pi_K) \partial \pi_K / \pi_K 
\end{multline*}
montre que le (1) de la proposition 
implique le (2) par une r{\'e}currence imm{\'e}diate.
\qed\end{proof}

\begin{remark}
Si $\partial(y)=x$ et $x\in\bnst{,r}{,K}$, alors on ne peut pas dire
que $y\in\bnst{,r}{,K}$ mais en revanche
$y\in\bnst{,s}{,K}$ pour tout $s>r$.
\end{remark}

Les trois lemmes suivants seront utiles par la suite.

\begin{lemma}\label{divt}
Soit $x\in\bnrig{,r}{,K}$ tel que pour tout $n \gg 0$ on ait 
$x \in \phi^{n-1}(q)\bnrig{,r_n}{,K}$. Alors $x\in t \bnrig{}{,K}$.
\end{lemma}

\begin{proof}
Soit $s=r_{n_0}$ tel que $x\in \phi^{n-1}(q)\bnrig{,r_n}{,K}$ pour
tout $n \geq n_0$. Rappelons que $\bnrig{,s}{,K}$ s'identifie (non
canoniquement) {\`a} un anneau de s{\'e}ries formelles en $\pi_K$.
Soit $Q_n(\pi_K)=\phi^{n-1}(q)$.
On a $\bnrig{,s}{,K} \cap
\phi^{n-1}(q)\bnrig{,r_n}{,K}=\phi^{n-1}(q)\bnrig{,s}{,K}$
car une s{\'e}rie est divisible $Q_n(\pi_K)$ si et seulement si les
z{\'e}ros de la premi{\`e}re ont un ordre $\geq$ {\`a} ceux de la seconde,
ce que l'on peut v{\'e}rifier localement. De m{\^e}me,
comme les $\phi^{n-1}(q)$ sont premiers entre eux,
l'hypoth{\`e}se du lemme est {\'e}quivalente au fait que,
pour tout $n \gg 0$, on ait $x\in \phi^n(\pi) / \phi^{n_0-1}(\pi)
\bnrig{,s}{,K}=\phi^n(\pi)\bnrig{,s}{,K}$ 
puisque $\phi^{n_0-1}(\pi)$ est inversible dans
$\bnrig{,s}{,K}$.
On peut donc 
{\'e}crire $x=\phi^n(\pi)p^{-n}x_n$ et on va montrer que la suite $\{x_n\}$
converge dans $\bnrig{,s}{,K}$. 
Les $x_n \in \bnrig{,s}{,K}$ v{\'e}rifient 
la relation:
\[ \frac{\phi^n(\pi)}{p^n}x_n-\frac{\phi^{n+1}(\pi)}{p^{n+1}}x_{n+1}=0 \]
et donc
\[ (x_n-x_{n+1})+x_{n+1}\left(1-\frac{\phi^n(q)}{p}\right)=0 \]
On en d{\'e}duit que si $I$ est un intervalle contenu dans $[s;+\infty[$, alors
\begin{align*}
V_I(x_n-x_{n+1}) &= V_I\left(x_{n+1}\left(1-\frac{\phi^n(q)}{p}\right)\right) \\
&\geq V_I(x_{n+1})+V_I\left(1-\frac{\phi^n(q)}{p}\right) 
\end{align*}
On se fixe $I$ un intervalle compact 
et comme $\lim_{n \ra +\infty}V_I(1-\phi^n(q)/p)=+\infty$
sur toute couronne du
type $\{ |z| \in I \}$, on peut
supposer que $V_I(1-\phi^n(q)/p)>0$
pour $n$ assez grand. Alors
l'in{\'e}galit{\'e} ci-dessus montre que $V_I(x_n-x_{n+1}) > V_I(x_{n+1})$ et donc
que $V_I(x_n)=V_I(x_{n+1})$. Les $x_n$ ont donc
tous la m{\^e}me valuation pour $n$ assez grand. Enfin l'in{\'e}galit{\'e}
ci-dessus et le fait que $\lim_{n \ra +\infty}V_I(1-\phi^n(q)/p) =+\infty$
impliquent que $\lim_{n \ra +\infty}V_I(x_n-x_{n+1})=+\infty$ et donc
que la suite $x_n$ converge.

Ceci {\'e}tant vrai pour tout $I$, la suite $x_n$ converge pour la
topologie de Fr{\'e}chet vers une limite $y \in \bnrig{,s}{,K}$ et
un calcul imm{\'e}diat montre que $x=ty$.
\qed\end{proof}

\begin{remark}
On peut aussi dire que $xt^{-1}$ est une fonction m{\'e}romorphe sans
p{\^o}les, et cela la force {\`a} {\^e}tre holomorphe, c'est-{\`a}-dire que
$xt^{-1}\in \bnrig{}{,K}$. Cela suit de r{\'e}sultats g{\'e}n{\'e}raux de
Lazard \cite{La62}
sur les fonctions analytiques, que l'on trouvera dans le paragraphe suivant.
\end{remark}

\begin{proposition}\label{kertheta}
Soit $r>0$ et $n\geq n(r)$. Alors
\[ \ker(\theta\circ\iota_n: \btrig{,r}{,K} \ra \Cp)
=\phi^{n-1}(q)\btrig{,r}{,K} \]
\[ \ker(\theta\circ\iota_n: \bnrig{,r}{,K} \ra \Cp)
=\phi^{n-1}(q)\bnrig{,r}{,K} \]
\end{proposition}

\begin{proof}
Le premier point
est une cons{\'e}quence imm{\'e}diate de \ref{kertheta1}. Pour le
deuxi{\`e}me si $x \in \bnrig{,r}{,K}$ est tel que
$\theta\circ\iota_n(x)=0$, alors $x=\phi^{n-1}(q)y$ avec $y \in
\btrig{,r}{,K}$ et alors $x=\phi^{n-1}(q)R_0(y)$ ce qui fait que
$y=R_0(y)$ et donc que $y \in \bnrig{,r}{,K}$.
\qed\end{proof}

\begin{lemma}\label{iotasurj}
Si $r>0$ et $n \gg n(r)$, alors
l'application injective \[ \theta\circ\iota_n :
\bnrig{,r}{,K}/\phi^{n-1}(q) \ra K_n \] est une bijection.
\end{lemma}

\begin{proof}
{\'E}tant donn{\'e} que $\phi^{n-1}(q)\bnrig{,r}{,K} \cap 
\bnrig{,r}{,F} = \phi^{n-1}(q)\bnrig{,r}{,F}$, l'extension de corps
$\bnrig{,r}{,K}/\phi^{n-1}(q)$ sur $\bnrig{,r}{,F}/\phi^{n-1}(q)$ est de 
degr{\'e} $\bnrig{,r}{,K}/\bnrig{,r}{,F}=e_K=[K_{\infty}:F_{\infty}]$. 

La proposition est triviale dans le cas $K=F$
(on a d{\'e}j{\`a} $\bnrig{,r}{,F}/\phi^{n-1}(q)=F_n$), 
et dans le cas
g{\'e}n{\'e}ral, elle r{\'e}sulte alors du fait que
$\bnrig{,r}{,K}/\phi^{n-1}(q)$ est une extension de $F_n$ de degr{\'e}
$e_K=[K_n:F_n]$ pour $n \gg 0$. 
\qed\end{proof}

\subsection{Modules sur $\bnrig{,s}{,K}$}
On \'etablit ici des r\'esultats techniques qui 
serviront \`a d\'emontrer un r\'esultat sur la
structure de $\dnrig{}(V)$ quand $V$ est de de Rham.

Soit $s \in \RR$, tel que $n(s) \geq n(K)$.
L'anneau $\bnrig{,s}{,K}$ est muni de sa topologie de Fr{\'e}chet,
d{\'e}finie par les $\{V_I\}$, et si $M$ est un $\bnrig{,s}{,K}$-module
libre de rang $d$, alors le choix d'un isomorphisme $M =
(\bnrig{,s}{,K})^d$ permet de munir $M$ d'une topologie. Le
th{\'e}or{\`e}me de l'image ouverte pour les espaces de Fr{\'e}chet montre
que cette topologie ne d{\'e}pend pas du choix d'une base de $M$.
L'objet de ce paragraphe est de montrer le th{\'e}or{\`e}me suivant et son corollaire:
\begin{theorem}\label{libre}
Soit $M$ un $\bnrig{,s}{,K}$-module libre de rang fini $d$, et $N \subset M$ 
un sous-module ferm{\'e} pour la topologie de Fr{\'e}chet de $M$, tel 
que \[ N \otimes_{\bnrig{,s}{,K}}
 \on{Frac}(\bnrig{,s}{,K}) = M  \otimes_{\bnrig{,s}{,K}}
 \on{Frac}(\bnrig{,s}{,K}). \]
Alors $N$ est libre de rang $d$.
\end{theorem}

\begin{corollary}\label{libre2}
Soit $M$ un $\bnrig{,s}{,K}$-module libre de rang fini $d$, et $N \subset M$ 
un sous-module ferm{\'e} pour la topologie de Fr{\'e}chet de $M$.
Alors $N$ est libre de rang $e \leq d$, o{\`u} \[ e=\dim_{\on{Frac}(\bnrig{,s}{,K})}
N \otimes_{\bnrig{,s}{,K}} \on{Frac}(\bnrig{,s}{,K}). \]
\end{corollary}

Ces deux r{\'e}sultats sont, apparamment, bien connus des experts, mais
nous en reproduisons la d{\'e}monstration pour la commodit{\'e} du
lecteur.
La d{\'e}monstration va n{\'e}cessiter quelques
r{\'e}sultats pr{\'e}paratoires sur la structure de $\bnrig{,s}{,K}$.
Ces r{\'e}sultats se trouvent dans \cite{An00,La62,CM95}, moyennant 
une identification de l'anneau $\bnrig{,s}{,K}$ 
{\`a} l'anneau des s{\'e}ries formelles convergeant sur la
couronne $C(I_s)=\{ x \in \Cp,\ \alpha(s,K) \leq |x|_p < 1 \}$, 
avec $\alpha(s,K)=p^{-1/e_K s}$. 
Si $I$ est un intervalle de $[0;1[$, soit $\mathbf{B}^I_K$ 
l'anneau des s{\'e}ries de Laurent {\`a} coefficients dans $F$, qui 
convergent sur la couronne $C(I)=\{ x \in \Cp,\ 
|x|_p \in I \}$. 
On fixe une identification 
$\bnrig{,s}{,K}=\mathbf{B}^{I_s}_K$.

\begin{proposition}
Soit $I$ un intervalle de $[0;1[$. L'anneau $\mathbf{B}^I_K$ a alors
les propri{\'e}t{\'e}s suivantes:
\begin{enumerate}
\item si $I$ est un intervalle compact, 
alors $\mathbf{B}^I_K$ est un anneau principal; 
\item tout id{\'e}al de type fini de $\mathbf{B}^I_K$ est 
principal (c'est donc un anneau de Bezout);
\item tout id{\'e}al ferm{\'e} de $\mathbf{B}^I_K$ est principal;
\item toute fonction m{\'e}romorphe $f \in \on{Frac}(\mathbf{B}^I_K)$, 
qui n'a pas de p{\^o}les, est en fait holomorphe, c'est-{\`a}-dire que $f \in  
\mathbf{B}^I_K$;
\item toute matrice {\`a} $d$ lignes et $e$ colonnes, 
$M \in \on{M}_{d \times e}(\mathbf{B}^I_K)$, 
peut s'{\'e}crire $M=P\on{diag}(f_i)Q$
o{\`u} $P\in \on{GL}_d(\mathbf{B}^I_K)$, 
$Q\in \on{GL}_e(\mathbf{B}^I_K)$
et $\on{diag}(f_i) \in \on{M}_{d \times e}(\mathbf{B}^I_K)$ est nulle sauf sur 
la diagonale, et 
$f_1 | \cdots | f_d$. En d'autres termes,
$\mathbf{B}^I_K$ admet une th{\'e}orie des diviseurs {\'e}l{\'e}mentaires.
\end{enumerate}
\end{proposition}

\begin{proof}
Pour les 4 premiers points, on se reportera {\`a} l'article de Lazard
\cite{La62}.
Pour le cinqui{\`e}me point: on commence par montrer que
$\mathbf{B}^I_K$ est un anneau \textit{ad{\'e}quat}, 
c'est-{\`a}-dire que $\mathbf{B}^I_K$ est un anneau int{\`e}gre o{\`u} 
tout id{\'e}al de type fini est principal, et qui v{\'e}rifie la condition 
technique suivante: pour tous
$a,b \in \mathbf{B}^I_K$, on peut {\'e}crire $a=a_1\cdot a_2$ 
avec $(a_1,b)=1$ et $(a_3,b) \neq 1$ pour tout $a_3$ qui divise $a_2$ (et 
qui n'est pas une unit{\'e}). C'est une cons{\'e}quence des r{\'e}sultats de
\cite{La62}, qui permettent de construire des fonctions holomorphes
dont l'ensemble des z{\'e}ros est fix{\'e} {\`a} l'avance.
Le (5) est alors une cons{\'e}quence directe de \cite[Theorem 3]{He49}. 
\qed\end{proof}

Enfin, on aura besoin du lemme suivant:

\begin{lemma}\label{sousbez}
Soit $M$ un $\mathbf{B}^I_K$-module, libre de rang $d$, et $M' \subset M$ un 
sous-module de type fini. Alors $M'$ est libre de rang $\leq d$. Si
$M' \otimes_{\bnrig{,s}{,K}} 
\on{Frac}(\mathbf{B}^I_K) = M  \otimes_{\bnrig{,s}{,K}} 
\on{Frac}(\mathbf{B}^I_K)$, 
alors $M'$ est 
libre de rang $d$.
\end{lemma}

\begin{proof}
Si l'on exprime les coordonn{\'e}es d'une famille g{\'e}n{\'e}ratrice de $M'$ 
dans une base de $M$, c'est une cons{\'e}quence imm{\'e}diate de la th{\'e}orie 
des diviseurs {\'e}l{\'e}mentaires.
\qed\end{proof}

\begin{proof}[du th{\'e}or{\`e}me \ref{libre}]
La d{\'e}monstration est largement inspir{\'e}e de \cite{Kr69}.
Soit donc $M$ un $\mathbf{B}^{I_s}_K$-module libre de type fini, et $N
\subset M$ un sous-module ferm{\'e}.
Soit $e \geq d$.
Si $I$ est un intervalle compact, soit $\mathcal{G}(e,I)$ l'ensemble des 
$(n_1,\cdots,n_e) \in N^e$ qui engendrent le $\mathbf{B}^I_K$-module $N 
\otimes_{\bnrig{,s}{,K}} \mathbf{B}^I_K$. 
On remarquera que si $I$ est ferm{\'e}, alors $\mathbf{B}^I_K$ est principal, 
et donc $N \otimes_{\bnrig{,s}{,K}} \mathbf{B}^I_K$ 
est un $\mathbf{B}^I_K$-module libre de 
rang fini $\leq d$. Comme
\[ N \otimes_{\bnrig{,s}{,K}} \on{Frac}(\bnrig{,s}{,K}) = 
M \otimes_{\bnrig{,s}{,K}} \on{Frac}(\bnrig{,s}{,K}), \] 
alors $N \otimes_{\bnrig{,s}{,K}} \mathbf{B}^I_K$ est en fait 
un $\mathbf{B}^I_K$-module libre de 
rang {\'e}gal {\`a} $d$.

Montrons premi{\`e}rement que $\mathcal{G}(e,I)$ est non-vide: comme l'anneau
$\mathbf{B}^I_K$ est 
principal, il existe $n'_1,\cdots,n'_d \in N \otimes_{\bnrig{,s}{,K}} 
\mathbf{B}^I_K$ 
qui forment une 
base de $N \otimes_{\bnrig{,s}{,K}} \mathbf{B}^I_K$. 
On a $n'_i = \sum a_{ij} p_j$ o{\`u} $a_{ij} \in 
\mathbf{B}^I_K$, et $p_1,\cdots p_f$ est une famille de $N$. 
Le $\bnrig{,s}{,K}$-module engendr{\'e} 
par les $p_i$ est un sous-module de type fini de $M$, il est donc libre de 
rang $d$, de base $q_1,\cdots,q_d$. Il n'est pas difficile de voir que 
les $q_i$ forment une base de $N \otimes_{\bnrig{,s}{,K}} \mathbf{B}^I_K$ 
sur $\mathbf{B}^I_K$ et donc que 
$\mathcal{G}(d,I)$ est non-vide. Ceci implique que $\mathcal{G}(e,I)$ 
est non-vide pour $e 
\geq d$.

Montrons ensuite que $\mathcal{G}(e,I)$ est 
un ouvert de $N^e$. Comme $\mathcal{G}(d,I)$ est 
non-vide, on peut se fixer $(n_1,\cdots,n_d) \in \mathcal{G}(d,I)$. Soit 
$(p_1,\cdots,p_e) \in \mathcal{G}(e,I)$, et $P$ 
la matrice des $p_i$ dans la base 
$n_j$. Comme $\bnrig{,s}{,K}$ admet la th{\'e}orie des diviseurs 
{\'e}l{\'e}mentaires, on 
peut {\'e}crire $P=ADC$, avec $A \in \on{GL}_e(\bnrig{,s}{,K})$, 
$C \in \on{GL}_d(\bnrig{,s}{,K})$, 
et $D$ nulle hors de la diagonale.
La matrice $A$ induit un hom{\'e}omorphisme de $N^e$.
Quitte {\`a} modifier la base $n_i$ par $C$, et {\`a} faire le changement de 
coordonn{\'e}es sur $N^e$ induit par $A$, on se ram{\`e}ne {\`a} supposer que
$p_i = \alpha_i n_i$ et $p_i=0$ si $i\geq d+1$. 

Soit $Q$ 
la matrice $d \times d$ des $\alpha_i$ avec $i \leq d$. Par hypoth{\`e}se, 
$\det(Q) \in (\mathbf{B}^I_K)^*$. 
Or l'ensemble des $R \in \on{M}_d(\bnrig{,s}{,K})$ dont le d{\'e}terminant
v{\'e}rifie $\det(R) \in (\mathbf{B}^I_K)^*$ est ouvert
(le d{\'e}terminant est une application continue et 
$(\mathbf{B}^I_K)^*$ est ouvert dans $\mathbf{B}^I_K$, puisque, $I$
{\'e}tant compact, 
$\mathbf{B}^I_K$ est une alg{\`e}bre de Banach) et contient $Q$. 
L'ensemble des $(r_1,\cdots,r_e)$ obtenus via $(r_i)=P(n_i)$ est ouvert et 
contient $(p_1,\cdots,p_e)$, 
Si l'on prend pour $P$ toutes les matrices 
$e \times d$ telles que: 
\begin{enumerate}
\item le d\'eterminant de
la matrice $d \times d$ des $d$ premi{\`e}res lignes 
de $P$ est inversible dans $\mathbf{B}^I_K$;
\item la matrice $e-d \times d$
des $e-d$ derni{\`e}res lignes de $P$ est quelconque.
\end{enumerate}
alors
l'ensemble des $(r_1,\cdots,r_e)$ obtenus via $(r_i)=P(n_i)$ est ouvert et 
contient $(p_1,\cdots,p_e)$.
Ceci montre que $\mathcal{G}(e,I)$ est ouvert.

Montrons enfin que si $e\geq 2d$, alors 
$\mathcal{G}(e,I)$ est dense dans $N^e$. 
Soient $(n_1,\cdots,n_d) \in \mathcal{G}(d,I)$, 
et $(p_1,\cdots,p_e) 
\in N^e$.
Comme auparavant, la th{\'e}orie des diviseurs {\'e}l{\'e}mentaires permet de 
supposer 
que $p_i = \alpha_i n_i$ et $p_i=0$ si $i \geq d+1$. 
Soit $n>0$, et $q_i=p_i$ pour $1 \leq i \leq d$, $q_i = p^n n_{i-d}$ 
pour $d+1 \leq i \leq 2d$, et $q_i=0$ si $i \geq 2d+1$.
La famille $(q_i)$ est une famille g{\'e}n{\'e}ratrice de $N 
\otimes_{\bnrig{,s}{,K}} 
\mathbf{B}^I_K$, aussi proche de $(p_i)$ que l'on veut (le choix de $n$ est 
libre).

L'ensemble $\mathcal{G}(2d,I)$ est donc un 
ouvert dense de $N^{2d}$, et comme $N$ 
est un ferm{\'e} d'un Fr{\'e}chet, c'est un m{\'e}trique complet: il a la 
propri{\'e}t{\'e} de Baire. L'intersection \[ \cap_{n \gg 0}^{\infty} 
\mathcal{G}(2d,[r;1-1/n]) \] 
est donc dense dans $N^{2d}$, et notamment non-vide. 
Soit  
$(n'_1,\cdots,n'_{2d})$ dans l'intersection. Le sous 
$\bnrig{,s}{,K}$-module de $N$ 
engendr{\'e} par les $n'_i$ est de type fini et $\subset M$, il est donc 
libre de rang $d$, engendr{\'e} par $n_1,\cdots,n_d$. On voit que ces $n_i$ 
sont une base de $N \otimes_{\bnrig{,s}{,K}} \mathbf{B}^I_K$ 
pour tout $I$. Soit $n \in N$. On peut 
donc {\'e}crire $n=\sum_{i=1}^d a_i n_i$ avec $a_i \in 
\mathbf{B}^I_K$. Par unicit{\'e}, on a 
$a_i \in \cap \mathbf{B}^I_K = \bnrig{,s}{,K}$. 

Le module $N$ est donc libre de rang $d$, {\'e}tant engendr{\'e} par les $n_i$.
\qed\end{proof}

\begin{proof}[du corollaire \ref{libre2}]
Soit $Q= N \otimes_{\bnrig{,s}{,K}}
\on{Frac}(\bnrig{,s}{,K}) \cap M$. Montrons que c'est un 
$\bnrig{,s}{,K}$-module libre de rang $e \leq d$. On a une suite exacte
$0 \ra Q \ra M \ra M/Q \ra 0$, et par construction, $M/Q$ est sans-torsion: 
c'est un module de type fini et sans torsion, donc coh{\'e}rent, ce qui
fait \cite[I.2 exercice 11]{BOAC}
que $Q$ lui-m{\^e}me est coh{\'e}rent. Comme il 
est de type fini, le lemme  \ref{sousbez} montre 
qu'il est libre de rang fini $e$.
 
Ensuite, il est imm{\'e}diat 
que $N \otimes_{\bnrig{,s}{,K}} \on{Frac}(\bnrig{,s}{,K}) = 
Q \otimes_{\bnrig{,s}{,K}} \on{Frac}(\bnrig{,s}{,K})$, et 
on peut alors appliquer le th{\'e}or{\`e}me \ref{libre}.
\qed\end{proof}

\section{Structures diff{\'e}rentielles sur les 
$(\phi,\Gamma_K)$-modules et monodromie $p$-adique}
Dans le chapitre pr{\'e}c{\'e}dent on a fait l'{\'e}tude de
$(\bnrig{}{,K},\nabla)$. Ce chapitre est consacr{\'e} {\`a} la
construction de l'{\'e}quation diff{\'e}rentielle associ{\'e}e {\`a} une
repr{\'e}sentation $V$, qui est un module {\`a} connexion 
au-dessus de $(\bnrig{}{,K},\nabla)$. On donne ensuite des
applications aux repr{\'e}sentations semi-stables,
{\`a} la th{\'e}orie de Sen, et {\`a} la caract{\'e}risation des repr{\'e}sentations
de de Rham. Comme cons\'equence de cela, on d\'emontre la conjecture de 
monodromie $p$-adique.

\subsection{L'op{\'e}rateur $\nabla_V$}
Dans tout ce paragraphe, $V$ est une repr{\'e}sentation $p$-adique de
$G_K$. Rappelons que l'on a pos{\'e} $\dnrig{}(V)
=\bnrig{}{,K} \otimes_{\bdag{}_K} \ddag{}(V)$ et notamment que 
$\dnrig{}(V) \subset (\btrig{}{} \otimes_{\Qp} V)^{H_K}
=\btrig{}{,K} \otimes_{\bdag{}_K} \ddag{}(V)$. 
Le lemme de r{\'e}gularisation
par le Frobenius montre qu'on a $(\btrig{}{})^{\phi=1} \subset 
(\btrigplus{})^{\phi=1}=\Qp$
et donc que l'on peut r{\'e}cup{\'e}rer $V$ par la formule 
\[ V=(\btrig{}{} \otimes_{\bnrig{}{,K}}
\dnrig{}(V))^{\phi=1}. \] 
On choisit une base $\{e_i\}$ de $\ddag{}(V)$ et on prolonge
les op{\'e}rateurs $R_m$ {\`a} $\btrig{}{,K} \otimes_{\bdag{}_K}
\ddag{}(V)$ par $R_m(\sum \lambda_i \otimes e_i) =  \sum R_m(\lambda_i)
\otimes e_i$. Comme $R_m$ est $\bdag{}_K$-lin{\'e}aire, le prolongement
ne d{\'e}pend pas du choix de la base $\{e_i\}$. 
Nous aurons besoin du r{\'e}sultat suivant qui est
un corollaire imm{\'e}diat de la proposition \ref{kertheta1}:

\begin{lemma}\label{kerthetadhol}
Soit $r>0$ et $n\geq n(r)$. Alors  
\begin{align*}
\ker(\theta\circ\iota_n: \dnrig{,r}(V) \ra \Cp\otimes_{\Qp} V)
&=\phi^{n-1}(q)\dnrig{,r}(V) \\
\ker(\theta\circ\iota_n: \ddag{,r}(V) \ra \Cp\otimes_{\Qp} V)
&=\phi^{n-1}(q)\ddag{,r}(V) 
\end{align*}
\end{lemma}
 
\begin{proof}
La proposition \ref{kertheta1} montre que si $\theta\circ\iota_n(x)=0$,
alors $x=\phi^{n-1}(q)y$ avec $y \in \btrig{}{,K}\otimes_{\bnrig{}{,K}}
\dnrig{}(V)$. On a
alors $x=\phi^{n-1}(q)R_0(y)$ et donc $y=R_0(y) \in \dnrig{}(V)$. Enfin
pour le deuxi{\`e}me point il suffit d'utiliser le fait qu'un {\'e}l{\'e}ment de
$\bdag{,r_n}_K$ divisible par $\phi^{n-1}(q)$ dans $\bnrig{,r_n}{,K}$ l'est
dans $\bdag{,r_n}_K$.
\qed\end{proof}

On va munir $\dnrig{}(V)$  d'une connexion
$\nabla_V$ au-dessus de $(\bnrig{}{,K},\nabla)$, c'est-{\`a}-dire d'un 
op{\'e}rateur $\nabla_V$ tel que
$\nabla_V(\lambda\cdot x)= \nabla (\lambda) \cdot x +
\lambda\cdot\nabla_V(x)$.

On se fixe une base $\{e_i\}$ de $\ddag{}(V)$ et si $x\in\dnrig{,r}(V)$
s'{\'e}crit $x=\sum x_i \otimes e_i$ 
avec $x_i \in \bnrig{,r}{,K}$, alors on pose 
(si $I$ est tel que $r \in I$):
$V_I(x)= \inf_i V_I(x_i)$. La d{\'e}finition d{\'e}pend de la base mais
pas la topologie d{\'e}finie par les $V_I$ qui est {\'e}quivalente {\`a} la 
topologie induite sur $\dnrig{}(V)$ par $\btrig{}{}\otimes_{\Qp} V$, et
$\dnrig{,r}(V)$ est complet.
Dans cette partie $\gamma$ est un {\'e}l{\'e}ment de $\Gamma_K$ et $n(\gamma)=
v_p(1-\chi(\gamma))$. On suppose que $n(\gamma)\geq 1$ et que
$r \geq r_{n(K)}$.

\begin{lemma}
Si $I=[r;s]$ est un intervalle ferm{\'e} de $[r;+\infty[$, alors
il existe $n(I,V) \in \NN$, tel que pour 
$x\in\dnrig{,r}(V)$, on ait
$V_I((1-\gamma)x)\geq V_I(x)+1$ d{\`e}s que $n(\gamma) \geq n(I,V)$.
\end{lemma}

\begin{proof}
Si $y\in \ddag{,r}(V)$, alors
$V_I((1-\gamma)y) \geq V_I(y)+1$, quand $n(\gamma)$ est assez grand,
puisque l'action de $G_K$ est continue. On se fixe $n(\gamma)$
tel que $V_I((1-\gamma)e_i) \geq V_I(e_i)+1$ pour tout $i$ 
et le lemme r{\'e}sulte alors
du lemme \ref{unmoinsgampetit}
(le cas de la repr{\'e}sentation triviale) puisque
\[ (1-\gamma)(x_i \otimes e_i)=(1-\gamma)x_i \otimes \gamma(e_i)
+x_i \otimes (1-\gamma)e_i \]
\qed\end{proof}

Si $I=[r;s]$,
soit $\mathbf{B}^I_K$ l'anneau des s{\'e}ries formelles en $\pi_K$ qui
convergent sur la couronne de rayons int{\'e}rieurs et ext{\'e}rieurs
$\alpha(K,r)$ et $\alpha(K,s)$. 
Le lemme pr{\'e}c{\'e}dent montre que si $\gamma$ est assez proche de $1$,
la s{\'e}rie d'op{\'e}rateurs \[ \frac{\log(\gamma)}{\log(\chi(\gamma))}
=-\frac{1}{\log(\chi(\gamma))}\sum_{n\geq 1}\frac{(1-\gamma)^n}{n}\]
converge vers un op{\'e}rateur continu $\nabla_{I,V} : \dnrig{,r}(V) \ra 
\dnrig{,r}(V) \otimes_{\bnrig{,r}{,K}} 
\mathbf{B}^I_K$, et un petit argument de s{\'e}ries formelles montre que 
$\log(\gamma)/\log(\chi(\gamma))$ ne d{\'e}pend pas du choix de
$\gamma$. Notamment $\nabla_{I_1,V}(x)=\nabla_{I_2,V}(x)$ si
$I_k=[r;s_k]$. On en d{\'e}duit que la valeur commune des $\nabla_{I,V}(x)$
appartient {\`a} $\dnrig{,r}(V)$, et que l'op{\'e}rateur $x \mapsto
\nabla_V(x)$, o{\`u} $\nabla_V(x)$ est la valeur commune des $\nabla_{I,V}(x)$,
est continu pour la topologie de Fr{\'e}chet\index{nablaV@$\nabla_V$}.

Ensuite, un argument standard montre que si $x \in \dnrig{,r}(V)$,
alors \[ \nabla_V(x) = \lim_{\gamma \ra 1} \frac{\gamma(x) -
  1}{\chi(\gamma)-1} \]
ce qui fait que, comme $(1-\gamma)(\lambda x)=(1-\gamma)\lambda \cdot
\gamma(x) + \lambda \cdot (1-\gamma)x$, $\nabla_V$ v{\'e}rifie $\nabla_V(\lambda
x)= \nabla(\lambda)x+\lambda \nabla_V(x)$. C'est donc une connexion
au-dessus de l'op{\'e}rateur $\nabla: \bnrig{,r}{,K} \ra \bnrig{,r}{,K}$.

\begin{example}
Un calcul facile montre par exemple que $\nabla_{\Qp(r)}=t\partial+r$
sur $\dnrig{}(\Qp(r))$.
\end{example}
 
\begin{lemma}\label{divtv}
Soit $x$ un \'el\'ement de $\dnrig{,r}(V)$, tel que pour tout $n \gg 0$, on ait 
$x \in \phi^{n-1}(q)\dnrig{,r_n}(V)$. Alors $x\in t \dnrig{}(V)$.
\end{lemma}

\begin{proof}
Apr{\`e}s le choix d'une base de $\dnrig{}(V)$ cela suit imm{\'e}diatement
de \ref{divt}.
\qed\end{proof}

\subsection{Application aux repr{\'e}sentations semi-stables}
Un cristal sur $\bnrig{}{,K}$ est un $\bnrig{}{,K}$-module libre
muni d'un Frobenius et d'une connexion qui commutent, la connexion
{\'e}tant au-dessus de $\nabla$.
Com\-men\c{c}ons par d{\'e}finir ce qu'est un cristal
unipotent.
On remarquera que cela ne d{\'e}pend pas de l'{\'e}ventuelle
structure de Frobenius. 
Soit donc $\mbf$ un $\bnrig{}{,K}$-module libre de rang fini $d$ muni d'une connexion
$\nabla_M$ au-dessus de $\nabla$.

\begin{proposition}
Les propri{\'e}t{\'e}s suivantes sont {\'e}quivalentes:
\begin{enumerate}
\item $\nabla_M$ est triviale sur $\mbf\otimes_{\bnrig{}{,K}}
\bnst{}{,K}$, 
c'est-{\`a}-dire qu'il existe 
\[ e_0,\cdots,e_{d-1} \in\mbf\otimes_{\bnrig{}{,K}} \bnst{}{,K} \]
tels que $\nabla_M e_i=0$ et
$\oplus e_i\bnst{}{,K} [1/t] =
\mbf\otimes_{\bnrig{}{,K}}\bnst{}{,K} [1/t]$;
\item il existe $d$ {\'e}l{\'e}ments $f_0, \cdots,f_{d-1}$ de $\mbf$
qui forment
une base de $\mbf\otimes_{\bnrig{}{,K}}
\bnrig{}{,K} [1/t]$ sur $\bnrig{}{,K}[1/t]$ et tels 
que $\nabla_M(f_i) \in t\cdot\pscal{f_{i-1},\cdots,f_0}$
o{\`u} $\pscal{\cdot}$ d{\'e}note le $\bnrig{}{,K}$-module engendr{\'e}.
\end{enumerate}
\end{proposition}

\begin{proof}
Commen{\c c}ons par montrer que (1) implique (2). 
L'op{\'e}\-rateur de monodro\-mie $N=1\otimes N$ 
laisse stable le $F$-espace vectoriel 
engendr{\'e} par les $e_i$, car le (1) implique que les $e_i$ engendrent le
noyau de $\nabla_M$ agissant sur $\mbf\otimes_{\bnrig{}{,K}}
\bnst{}{,K}$, et un calcul direct montre que $N$ commute {\`a} $\nabla_M$
(car $N$ commute {\`a} l'action de $G_F$)
et donc stabilise son noyau:
on peut alors supposer que $N(e_i) \in \pscal{e_{i-1},\cdots,e_0}$
car $N$ est nilpotent.
On peut {\'e}crire de mani{\`e}re unique 
$e_i = \sum_{j=0}^{d-1} \log^j(\pi) d_{ji}$.
Montrons que $f_i=d_{0,i}$ est une famille qui satisfait la condition de 
(2). Le fait que $\nabla_M(e_i)=0$ et 
$N(e_i) \in \pscal{e_{i-1},\cdots,e_0}$
implique respectivement que 
\[ \nabla_M(d_{0,i})=d_{1,i} \frac{t(1+\pi)}{\pi} \text{ et }
d_{1,i} \in \pscal{d_{0,i-1},\cdots,d_{0,0}} \]
ce qui fait que $\nabla_M(d_{0,i}) \in t
 \pscal{d_{0,i-1},\cdots,d_{0,0}}$.
Montrons enfin que les $f_i$ engendrent
$\mbf\otimes_{\bnrig{}{,K}}  
\bnrig{}{,K} [1/t]$ sur $\bnrig{}{,K}[1/t]$. Soit
$m \in \mbf\otimes_{\bnrig{}{,K}}
 \bnrig{}{,K}[1/t]$, par hypoth{\`e}se on peut {\'e}crire
$m =\sum \mu_i e_i$ et donc $m = \sum \lambda_i f_i +n \cdot \log(\pi)$ 
o{\`u} $\lambda_i$ est le terme constant de $\mu_i$ et
comme $m \in \mbf\otimes_{\bnrig{}{,K}}
 \bnrig{}{,K}[1/t]$ on a n{\'e}c{\'e}ssairement $n=0$.

Montrons maintenant l'implication r{\'e}ciproque.
On va montrer par r{\'e}cur\-rence que l'on peut prendre 
$e_i = \sum_{j=0}^i f_j a_{ji}$ avec $a_{ji} \in \bnst{}{,K}$
et $a_{ii}=1$. 
En rang $1$ il n'y a rien
{\`a} montrer. En rang $d$, $\nabla_M$ induit une connexion
sur $(\oplus f_i\bnrig{}{,K})/f_0 \bnrig{}{,K}$ 
qui satisfait les m{\^e}mes conditions 
et il existe donc $e'_1,\cdots,e'_{d-1}$ tels que $\nabla_M(e'_i)=\alpha_i 
f_0$. Un calcul imm{\'e}diat montre que
$\alpha_i \in t \bnst{}{,K}$ et donc qu'il 
existe $\beta_i \in \bnst{}{,K}$ tel que $\nabla(\beta_i)=\alpha_i$
(par la proposition \ref{invconn}, il existe 
$\beta_i \in \bnst{}{,K}$ tel que $\partial(\beta_i)=t^{-1}\alpha_i$).
On pose alors $e_0=f_0$ et $e_i=e'_i-\beta_i f_0$, ce qui ach{\`e}ve la 
r{\'e}currence. La matrice de passage des $f_i$ aux $e_i$ est triangulaire
avec des $1$ sur la diagonale et on en d{\'e}duit que les $e_i$ engendrent 
bien $\mbf\otimes_{\bnrig{}{,K}}
 \bnst{}{,K}[1/t]$ sur $\bnst{}{,K}[1/t]$.
\qed\end{proof}

Un cristal qui satisfait les conditions de la proposition pr{\'e}c{\'e}dente
est dit unipotent. 

\begin{proposition}\label{conilp}
Soit $V$ une repr{\'e}sentation $p$-adique et $\dnrig{}(V)$ le cristal qu'on
lui a associ{\'e}. Alors il existe $n$ tel que la restriction de $V$ {\`a}
$G_{K_n}$ est semi-stable (respectivement cristalline) {\`a} poids n{\'e}gatifs
si et seulement si $\dnrig{}(V)$ est unipotent (respectivement trivial).
\end{proposition}

\begin{proof}
On a vu que $V$ est une repr{\'e}sentation
semi-stable de $G_{K_n}$
si et seulement si \[ \dnst{}(V)=\dnrig{}(V) \otimes_{\bnrig{}{,K}}
\bnst{}{,K}[1/t] \] a une base d'{\'e}l{\'e}ments stables par 
$\gamma_{K_n}$.
Si $V$ est semi-stable, alors $\nabla_V$ est donc triviale, c'est-{\`a}-dire 
que si $e_0,\cdots,e_{d-1}$ est une base de $\dst(V)$, alors ils satisfont
le (1) de la proposition pr{\'e}c{\'e}dente en raison du th{\'e}or{\`e}me de 
comparaison:
\[ \dst(V)\otimes_F \bnst{}{,K}[1/t]=
\dnrig{}(V)\otimes_{\bnrig{}{,K}}
\bnst{}{,K}[1/t] \]

R{\'e}ciproquement si $\dnrig{}(V)$ est unipotent, alors les $e_0,\cdots,e_{d-1}$
engendrent un $F$-espace vectoriel sur lequel $\log(\gamma)$ agit
trivialement ce qui fait que $\Gamma_K$ agit {\`a} travers un quotient
fini et donc que les $e_i$
sont stables par $\gamma_K^{p^n}$ pour $n$ assez grand et
forment alors une base de $\dst(V_n)$ si $V_n$ est la restriction de 
$V$ {\`a} $G_{K_n}$.  

On dit aussi dans ce cas que la connexion est unipotente.
De plus $V$ est cristalline si et seulement si on peut choisir les $f_i$
tels que $\nabla_M(f_i)=0$ et la connexion est alors triviale.
\qed\end{proof}

\subsection{Les modules $\dsen(V)$ et $\ddif(V)$}
Rappelons que Sen a montr{\'e} \cite{Sn73}
que si $V$ est une repr{\'e}\-sentation
$p$-adique, alors l'ensemble des sous-$K_{\infty}$-espaces vectoriels
de dimension finie
de $(\Cp\otimes_{\Qp} V)^{H_K}$ stables par $\Gamma_K$ admet un plus grand
{\'e}l{\'e}ment $\dsen(V)$\index{Dsen@$\dsen(V)$} 
et que
$\Cp\otimes_{K_{\infty}}\dsen(V)=\Cp\otimes_{\Qp} V$. De plus si $\gamma \in
\Gamma_K$ est suffisament proche de $1$, alors la s{\'e}rie qui d{\'e}finit
$\log(\gamma)$ converge en tant que s{\'e}rie d'op{\'e}rateurs
$\Qp$-lin{\'e}aires de $\dsen(V)$ et l'op{\'e}rateur $\Theta_V =
\log(\gamma)/\log(\chi(\gamma))$\index{ThetaV@$\Theta_V$} 
est un op{\'e}rateur
$K_{\infty}$-lin{\'e}aire qui ne d{\'e}pend pas de $\gamma$.

D'autre part Fontaine a montr{\'e} dans \cite{Fo00} que l'ensemble
des sous-$K_{\infty}[[t]]$-modules
libres de type fini
de $(\bdR^+\otimes_{\Qp} 
V)^{H_K}$ stables par $\Gamma_K$ contient un
plus grand {\'e}l{\'e}ment; nous le noterons $\ddif^+(V)$\index{Dddifplus@$\ddif^+(V)$}
dans cet article
et on a $\bdR^+\otimes_{K_{\infty}[[t]]}\ddif^+(V)=
\bdR^+ \otimes_{\Qp} V$. De plus si $\gamma\in
\Gamma_K$ est suffisament proche de $1$, alors la s{\'e}rie qui d{\'e}finit
$\log(\gamma)$ converge en tant que s{\'e}rie d'op{\'e}rateurs
$\Qp$-lin{\'e}aires de $\ddif^+(V)$ et l'op{\'e}rateur $\nabla_V =
\log(\gamma)/\log(\chi(\gamma))$\index{nablaV@$\nabla_V$} 
est un op{\'e}rateur
qui ne d{\'e}pend pas de $\gamma$ et qui v{\'e}rifie $\nabla_V(ax)=a
\nabla_V(x)+ \nabla(a) x$ ce qui montre que $\nabla_V$ est une
connexion sur $\ddif^+(V)$; on l'{\'e}tend {\`a} 
$\ddif(V)=K_{\infty}((t)) \otimes_{K_{\infty}[[t]]} 
\ddif^+(V)$\index{nablaV@$\nabla_V$}.

On peut retrouver $(\dsen(V),\Theta_V)$ {\`a} partir de
$(\ddif^+(V),\nabla_V)$ via l'appli\-cation $\theta:\bdR^+ \ra \Cp$
comme nous le verrons ci-dessous.
Rappelons \cite {CC99} que l'application $\iota_n=\phi^{-n}$ 
envoie $\bdag{,r_n}_K$ dans $K_n[[t]] \subset \bdR^+$ et
envoie donc $\ddag{,r_n}(V)$ dans un sous $K_n[[t]]$-module de
$\ddif^+(V)$.

\begin{proposition}\label{phiconn}
L'application de $K_{\infty}((t))\otimes_{\bdag{,r_n}_K}\ddag{,r_n}(V)$ dans
$\ddif(V)$ d{\'e}duite de $\iota_n$ est un isomorphisme de 
$K_{\infty}((t))$-modules avec connexion, si $n$ est assez grand.
\end{proposition}

\begin{proof}
On prend $n \geq n_0$ tel que $\iota_n(\bdag{,r_n}_K)\subset K_n[[t]]$.
Il est clair qu'alors 
\[ K_{\infty}[[t]] \otimes_{\iota_n(\bdag{,r_n}_K)} 
\iota_n(\ddag{,r_n}(V))\]
est un sous $K_{\infty}[[t]]$-module libre de type fini
de $(\bdR^+ \otimes_{\Qp} V)^{H_K}$ stable
par $\Gamma_K$.

Montrons que c'est $\ddif^+(V)$ pour $n \gg 0$. 
On a une application
$\theta : \ddif^+(V) \ra \dsen(V)$
et $\dsen(V)$
est un $K_{\infty}$-espace vectoriel
de dimension $d=\dim_{\Qp}(V)$. 
On en d{\'e}duit une application \[ \theta \circ \iota_n: \ddag{,r_n}(V)
\ra \dsen(V) \] dont le noyau est
$\phi^{n-1}(q) \ddag{,r_n}(V)$ par le lemme \ref{kerthetadhol}
ce qui fait que $\theta \circ \iota_n$
r{\'e}alise une injection d'un $\bdag{,r_n}_K/\phi^{n-1}(q)$ module de rang
$d$ dans $\dsen(V)$ et son image est un
$K_n$-espace vectoriel $V_n$ de dimension $d$ stable par 
$\Gamma_K$. Un petit calcul montre que des {\'e}l{\'e}ments de 
$V_n$ qui sont li\'es dans $\dsen(V)$ par une relation
\`a coefficients dans $K_{\infty}$,
sont li\'es par une relation \`a coefficients dans $K_n$,
et donc que l'application naturelle de
$K_{\infty} \otimes_{K_n} V_n$ dans $\dsen(V)$ est injective; son image
est de dimension maximale $d$, et est donc {\'e}gale
{\`a} $\dsen(V)$.
Cela implique que le d{\'e}terminant de l'injection
\[ K_{\infty}[[t]] \otimes_{\iota_n(\bdag{,r_n}_K)} \iota_n(\ddag{,r_n}(V))
\hookrightarrow \ddif^+(V) \] n'est pas divisible par $t$ et donc (comme
$K_{\infty}[[t]]$ est un anneau local d'id{\'e}al maximal $(t)$) par le lemme
de Nakayama que l'injection ci-dessus est en fait un isomorphisme.
\qed\end{proof}

\begin{corollary}
Par extension des scalaires on en d{\'e}duit que
l'application 
\[ K_{\infty}((t))\otimes_{\bnrig{,r_n}{,K}}\dnrig{,r_n}(V) \ra \ddif(V) \]
d{\'e}duite de $\iota_n$ est aussi un 
isomorphisme de $K_{\infty}((t))$-modules avec connexion, 
si $n$ est assez grand. 
\end{corollary}

Gr{\^a}ce {\`a} ces calculs on peut retrouver $\ddR(V)$ {\`a} partir de
$\ddag{}(V)$; la proposition suivante se trouve dans \cite[prop. 3.25]{Fo00}.

\begin{proposition}\label{kernabla}
Si $V$ est une repr{\'e}sentation $p$-adique de $G_K$, alors 
$K_{\infty} \otimes_K \ddR(V)$ est le noyau de la connexion $\nabla_V$ 
op{\'e}rant sur $\ddif(V)$. 
En particulier, $V$ est de de Rham si et seulement 
si $\nabla_V$ est la connexion triviale.
\end{proposition}

\begin{proof}
L'action de $\nabla_V$ sur $K_{\infty} \otimes_K \ddR(V)$
est triviale et de plus si $r \gg 0$, alors
$t^r K_{\infty}[[t]] \otimes_K \ddR(V)$ est un
sous-$K_{\infty}[[t]]$-module 
de $(\bdR^+ \otimes_{\Qp} V)^{H_K}$ libre de type fini et stable par
$\Gamma_K$ ce qui montre que $K_{\infty} \otimes_K \ddR(V)$ est bien
inclus dans le noyau de $\nabla_V$ agissant sur $\ddif(V)$. 

La th{\'e}orie g{\'e}n{\'e}rale des modules {\`a} connexion montre que le
noyau de $\nabla_V$ sur $\ddif(V)$ est un $K_{\infty}$-espace
vectoriel de dimension finie au plus $\dim_{K_{\infty}((t))}(\ddif(V))$;
il est aussi 
stable par $\Gamma_K$ et provient donc
par extension des scalaires d'un $K_n$-espace vectoriel stable par
$\Gamma_K$ pour $n \gg 0$. L'action de $\Gamma_K$ sur ce 
$K_n$-espace vectoriel est discr{\`e}te car 
l'alg{\`e}bre de Lie de $\Gamma_K$ agit trivialement, puisque
$\nabla_V=0$ 
(cf \cite[chap. V]{Se64}), et donc quitte
{\`a} augmenter $n$ l'action de $\Gamma_{K_n}$ est triviale et cet
espace est donc inclus dans 
$\ddR(V_n)=K_n \otimes \ddR(V)$ o{\`u} $V_n$ est 
la restriction {\`a} $G_{K_n}$ de $V$, ce qui fait que le noyau de la
connexion agissant sur $\ddif(V)$ est inclus dans $K_{\infty} \otimes_K \ddR(V)$.
On a donc
bien que $K_{\infty} \otimes_K \ddR(V)$ est le noyau de 
$\nabla_V$ sur $\ddif(V)$.
\qed\end{proof}

\subsection{Repr{\'e}sentations de de Rham}\label{ndrsec}
Soit $\partial_V=t^{-1}\nabla_V$. 
L'objet de ce paragraphe est de d{\'e}montrer le th{\'e}or{\`e}me suivant, et sa
r{\'e}ciproque, la proposition \ref{recipndr}:

\begin{theorem}\label{ndr}
Soit $V$ une repr{\'e}sentation de de Rham de $G_K$.
On suppose que les poids de 
Hodge-Tate de $V$ sont n{\'e}gatifs.
Alors il existe un unique sous
$\bnrig{}{,K}$-module libre de rang $d$
de $\dnrig{}(V)$ stable par $\partial_V$, $\ndr(V)$. 
Il v{\'e}rifie de plus les propri{\'e}t{\'e}s suivantes:
\begin{enumerate}
\item $\ndr(V)$ est stable par $\phi$ et par $\Gamma_K$;
\item si $s$ est assez grand, alors il existe $N_s \subset \dnrig{,s}(V)$,
libre de rang $d$, stable par $\partial_V$ et $\Gamma_K$, tel que
$\ndr(V) = N_s \otimes \bnrig{}{,K}$.
\end{enumerate}
\end{theorem}\index{NdRV@$\ndr(V)$}

On va montrer que $\ndr(V)= N_s \otimes \bnrig{}{,K}$, o{\`u}
$N_s$ l'ensemble des $x \in \dnrig{,s}(V)$, tels que pour tout $n
\geq n_0$, on ait $\iota_n(x) \in K_n[[t]] \otimes_K \ddR(V)$. 
Le module $\ndr(V)$ que l'on construit ainsi
est l'analogue du module $N$
construit par N. Wach \cite{Wa96}
dans le cas d'une repr{\'e}sentation absolument cristalline, ce qui
explique la notation.

La d{\'e}monstration de ce r{\'e}sultat va n{\'e}cessiter plusieurs {\'e}tapes
interm{\'e}diaires, consacr{\'e}es {\`a} la construction d'{\'e}l{\'e}ments de
$\dnrig{}(V)$. Pour l'instant, $V$ est une repr{\'e}sentation $p$-adique
quelconque.
On commence par le cas de la repr{\'e}sentation triviale:
\begin{lemma}\label{dense}
Soient $s \in \RR$, et $n \in \NN$, tels que $n \geq n(s) \geq n(K)$, et $w \in \NN$.
Si $f(t) \in K_n[[t]]$, alors il existe $\mu \in \bnrig{,s}{,K}$ tel que 
$\iota_n(\mu)-f(t) \in t^wK_n[[t]]$.
\end{lemma}

\begin{proof}
On cherche {\`a} montrer que $\iota_n(\bnrig{,s}{,K})$ est dense dans
$K_n[[t]]$ pour la topologie $t$-adique.
Comme $t \in \bnrig{,s}{,K}$, il suffit de montrer que l'application 
naturelle $\iota_n : \bnrig{,s}{,K} \ra K_n = K_n[[t]]/(t)$ est surjective. 
Or, cette application co{\"\i}ncide avec $\theta \circ \iota_n : \bnrig{,s}{,K}
\ra K_n$, qui est surjective par le lemme \ref{iotasurj}.
\qed\end{proof}

\begin{lemma}\label{partunit}
Soit $t_n= p^n t / \phi^{n-1}(q)$. On a \[ \theta \circ \iota_m (t_n) =
\begin{cases} 
0 \text{ si $n \neq m$;} \\ 
\eps^{(1)}-1 \text{ si $n=m$.} \end{cases} \]
\end{lemma}

\begin{proof}
Si $m=n$, alors 
\[ \theta \circ \iota_m (t_n) = 
\theta(t \frac{\pi_1}{\pi})=\theta(\pi_1 \cdot \frac
{t}{\pi})=\eps^{(1)}-1. \] 
Sinon, si $m < n$ 
\[ \theta \circ \iota_m (t_n) = 
\theta(\phi^{n-m-1}(\pi 
\frac{t}{\phi(\pi)}))=0 \]
et si $m > n$ 
\[  \theta \circ \iota_m (t_n) = \theta(p^{n-m} \frac{t}{\pi}  \frac{\pi}{\pi_1} \frac 
{\pi_1}{\pi_{m-n}} \pi_{m-n+1})=0, \] puisque $\theta(\frac{\pi}{\pi_1})=0$.
\qed\end{proof}

Soit $\{e_i\}$ une base de $\ddag{}(V)$. On suppose dans la suite que
$s$ est un r{\'e}el tel que $e_i \in \ddag{,s}(V)$, et que la matrice de
$\phi$ dans cette base a ses coefficients dans
$\bnrig{,s}{,K}$. Enfin, on suppose aussi que $\partial(\pi_K)$ est
inversible dans $\bnrig{,s}{,K}$. Bien entendu, toutes ces conditions
sont v{\'e}rifi{\'e}es si $s$ est assez grand.

\begin{lemma}
Si $x\in\ddif^+(V)$ et $w \in \NN$, alors il existe  $n_0=n_0(x,w)$, 
et $x_{n_0} \in \dnrig{,s}(V)$ tels que:
$\iota_{n_0}(x_{n_0})- x \in t^w \ddif^+(V)$.
\end{lemma}

\begin{proof}
On a vu en \ref{phiconn} que si $n$ est assez grand, alors
\[ K_{\infty}[[t]] \otimes_{\iota_n(\bnrig{,r_n}{,K})} 
\iota_n(\dnrig{,r_n}(V)) = 
\ddif^+(V).  \] 
Si $n$ v{\'e}rifie de plus $n \geq n(s)$, 
alors comme on a suppos{\'e} que 
$\dnrig{,r_n}(V) = \bnrig{,r_n}{,K} \otimes_{\bnrig{,s}{,K}}
\dnrig{,s}(V)$, on a 
$K_{\infty}[[t]] \otimes_{\iota_n(\bnrig{,s}{,K})}
\iota_n(\dnrig{,s}(V)) = 
\ddif^+(V)$. 

On peut donc {\'e}crire $x = \sum_{i=1}^d f_i(t) \otimes \iota_n(e_i)$ avec $f_i(t) 
\in K_{\infty}[[t]]$. 
Soit $n_0=n_0(x,w) \geq n$, tel que pour tout $i$, 
$f_i(t) \in K_{n_0}[t]+t^w K_{\infty}[[t]]$.
Soit $P_{n_0-n}=(p_{ij}) \in \on{M}_d(\bnrig{,p^{n_0-n-1}s}{,K})$ 
la matrice de $\phi^{n_0-n}$ dans la base $e_i$
(si $P=P_1$ est la matrice de $\phi$, on a \[ P_{n_0-n}=\phi^{n_0-n-1}(P)\cdots \phi(P)P, \] 
et si $P \in \on{M}_d(\bnrig{,s}{,K})$, alors $P_{n_0-n} \in
\on{M}_d(\bnrig{,p^{n_0-n-1}s}{,K})$).
On a 
\begin{align*} 
x &= \sum_{i=1}^d f_i(t) \otimes \iota_{n_0}( \phi^{n_0-n}(e_i) ) \\
 &= \sum_{i=1}^d 
\sum_{j=1}^d f_j(t) \iota_{n_0}(p_{ij}) \otimes \iota_{n_0}(e_i)  \\
 &= \sum_{i=1}^d g_i(t) \otimes \iota_{n_0}(e_i) 
\end{align*}
o{\`u} $g_i(t) \in K_{n_0}[t]+t^w K_{\infty}[[t]]$, comme on le v{\'e}rifie rapidement.

Par le lemme \ref{dense}, 
il existe $\mu_i \in \bnrig{,s}{,K}$, 
tels que $\iota_{n_0}(\mu_i)-g_i(t) \in t^w K_{\infty}[[t]]$. 
On peut alors poser $x_{n_0} = \sum \mu_i  e_i$.
\qed\end{proof}

On conserve les notations du lemme pr{\'e}c{\'e}dent.

\begin{proposition}\label{constrdrig}
Si $x \in \ddif^+(V)$, alors pour tout $n \geq n_0$,
il existe $x_n \in \dnrig{,s}(V)$ tel que 
$\iota_n(x_n)-x \in t^w \ddif^+(V)$.
\end{proposition}

\begin{proof}
Soit $x_{n_0}$ l'{\'e}l{\'e}ment construit dans le lemme pr{\'e}c{\'e}dent. Il 
satisfait donc: 
$\iota_{n_0}(x_{n_0})- x \in t^w \ddif^+(V)$. Notamment, $\iota_n( 
\phi^{n-n_0}(x_{n_0}) ) - x \in t^w \ddif^+(V)$. Soit 
$x'_n=\phi^{n-n_0}(x_{n_0})$, alors $x'_n \in \dnrig{,p^{n-n_0}s}(V)$, et 
donc $x'_n= \sum \mu'_n  e_i$ avec $\mu'_n \in 
\bnrig{,p^{n-n_0}s}{,K}$ et $e_i \in \dnrig{,s}(V)$. Le lemme \ref{dense}
fournit des $\mu_n \in \bnrig{,s}{,K}$ tels que 
$\iota_n(\mu_n)-\iota_n(\mu'_n) \in t^w K_{\infty}[[t]]$. On pose alors 
$x_n = \sum \mu_n e_i$.
\qed\end{proof}

On suppose maintenant que la repr{\'e}sentation $V$ est de de Rham, et
que ses poids de Hodge-Tate sont n{\'e}gatifs, ce qui fait qu'il existe
$w \in \NN$ tel que
\[ t^w \ddif^+(V) \subset K_{\infty}[[t]] \otimes_K \ddR(V)
\subset \ddif^+(V) \] 
Ensuite, il existe $n_1$ tel que, si $n\geq n_1$, alors
$\iota_n(\dnrig{,s}(V)) \subset K_n((t)) \otimes_K \ddR(V)$.
Soient $r_1, \cdots, r_d$ une base du $K$-espace
vectoriel $\ddR(V)$, et $n_0=\sup \{ \{ n_0(r_i,w) \}_i ,n_1 \}$.

\begin{proposition}
Soit $N_s$ l'ensemble des $x \in \dnrig{,s}(V)$, tels que pour tout $n
\geq n_0$, on ait $\iota_n(x) \in K_n[[t]] \otimes_K \ddR(V)$. Alors
$N_s$ est un $\bnrig{,s}{,K}$-module libre de rang $d$, et 
$\nabla_V(N_s) \subset t N_s$.\index{Ns@$N_s$}
\end{proposition}

\begin{proof}
Commen{\c c}ons par remarquer que $N_s$ contient $t^w \dnrig{,s}(V)$,
en particulier il est non-vide, et il contient un sous $\bnrig{,s}{,K}$-module
libre de rang $d$. Ensuite, il est ferm{\'e} pour la topologie de
Fr{\'e}chet de $\dnrig{,s}(V)$, car $\iota_n$ est continu pour tout
$n\geq n_0$. Le th{\'e}or{\`e}me \ref{libre} 
montre alors que $N_s$ est libre de rang $d$. Il en existe donc une
base $f_1,\cdots,f_d$. Pour $n \geq n_0$, la famille $\iota_n(f_i)$
forme une base de $K_n[[t]] \otimes_K \ddR(V)$: en effet, il suffit de
v{\'e}rifier que $\iota_n(N_s)= K_n[[t]] \otimes_K \ddR(V)$, c'est {\`a}
dire de montrer que pour tout $i$, $r_i \in \iota_n(N_s)$.
La proposition \ref{constrdrig},
appliqu{\'e}e {\`a} $r_i$, fournit $s_{i,n} \in \dnrig{,s}(V)$, tel que 
$\iota_n(s_{i,n})-r_i \in t^w \ddif^+(V)$. Soit $t_{i,n}=s_{i,n}
\left(\frac{p^nt}{\phi^{n-1}(q)}\right)^w$. 
Le lemme \ref{partunit} montre qu'il existe $u_{i,n} \in K_n[[t]]$, tel que
$\iota_n(t_{i,n})-((\eps^{(1)}-1)^w +t u_{i,n}) r_i \in t^w \ddif^+(V)$, et que si
$m \neq n$, alors $\iota_m(t_{i,n}) \in t^w \ddif^+(V)$.
On en d{\'e}duit que $t_{i,n} \in N_s$, et donc que $r_i \in \iota_n(N_s)$.

Ensuite, il est clair que $\nabla_V(N_s) \subset N_s$, car $\nabla_V$
commute {\`a} $\iota_n$, et $K_n[[t]] \otimes_K \ddR(V)$ est stable par
$\nabla_V$. Soit $D=(d_{ij})$ la matrice de la connexion $\nabla_V$ dans la
base $\{f_i\}$. Le fait que $\nabla_V(K_n[[t]] \otimes_K \ddR(V))
\subset t K_n[[t]] \otimes_K \ddR(V)$, montre que pour tout $n \geq
n_0$, $\iota_n(D) \in \on{M}_d(tK_n[[t]])$, et donc  par le lemme
\ref{kertheta1}, $D \in  \on{M}_d(\phi^{n-1}(q) \bnrig{,s}{,K})$.
Ceci {\'e}tant valable pour tout $n \geq n_0$, le lemme \ref{divt}
montre que $D \in  \on{M}_d(t \bnrig{,s}{,K})$, et donc que $\nabla_V(N_s)
\subset t N_s$. 
\qed\end{proof}

Montrons maintenant qu'un sous-module de $\dnrig{,s}(V)$, dont le rang
est maximal, et qui est stable par $\partial_V$, est unique.

\begin{lemma}\label{drdet}
Soient $V$ une repr{\'e}sentation de de Rham, 
dont les poids de Hodge-Tate sont n{\'e}gatifs,
et $M_s$ un $\bnrig{,s}{,K}$-module libre de rang $d$, inclus dans
$\dnrig{,s}(V)$. Si $\partial_V(M_s) \subset M_s$, alors
$\det_{\bnrig{,s}{,K}}(M_s)=
t^r \det_{\bnrig{,s}{,K}}(\dnrig{,s}(V))$, o\`u $r \geq 0$.
\end{lemma}

\begin{proof}
La repr{\'e}sentation $\det(V)$ est de de Rham, et donc 
\cite{Bu88sst} de la forme
$\chi^{-r}\omega$
($r \geq 0$, c'est l'oppos{\'e} de la somme des poids de Hodge-Tate de
$V$), o{\`u} la restriction de $\omega$ {\`a} $I_K$\index{IK@$I_K$}
(l'inertie de $G_K$) est finie. 
Notamment, il existe une base $e$ de $\dnrig{,s}(\det(V))$, telle
que l'action de $\Gamma_K$ sur $e$  est
donn{\'e}e par $\chi^{-r}\omega'$, o{\`u} $\omega'$ est d'ordre fini.
On a donc
$\gamma(e)=\chi^{-r}(\gamma)e$, si $\gamma$ est assez proche de $1$.
Il existe $\lambda \in \bnrig{,s}{,K}$, tel que
$\det(M_s)=\lambda \bnrig{,s}{,K} e$, et c'est un 
sous-module de 
$\det(\dnrig{,s}(V)) = \dnrig{,s}(\det(V))$ 
de rang $1$, stable par $\partial_V$. 
On doit avoir:
$\partial_V(\lambda e)=\alpha \lambda e$, 
avec $\alpha \in \bnrig{,s}{,K}$
et donc $\alpha \lambda e = 
\partial(\lambda) e + \lambda \partial_V(e)$. Comme 
$\partial_V(e)=-rt^{-1}e$, cela montre que $t$ divise $\lambda$, si $r \geq 1$.  
Une r{\'e}currence facile, sur $r$, montre 
qu'il existe $\mu \in \bnrig{,s}{,K}$, tel que
$\lambda=t^r \mu$. 
On doit alors avoir $\partial(\mu)=\beta \mu$, avec $\beta \in 
\bnrig{,s}{,K}$. Identifions cet anneau avec l'anneau des s{\'e}ries en $T$,
holomorphes sur la couronne $\{ z \in \Cp,\ \alpha(K,r) \leq |z|_p < 1\}$. 
Comme on a suppos{\'e} que 
$\partial(\pi_K)$ est inversible dans $\bnrig{,s}{,K}$, si $\mu(x)=0$, 
alors $(d\mu/dT)(x)=0$, et par r{\'e}currence $(d^k\mu/dT^k)(x)=0$
pour tout $k$. Si $\mu$ 
s'annule quelque part, c'est donc que $\mu=0$. On en conclut que $\mu$ est 
inversible et donc que $\det(M_s)=t^r \det(\dnrig{,s}(V))$.
\qed\end{proof}

\begin{corollary}
Si $N_s^1$ et $N_s^2$ sont deux sous-modules libres de rang $d$ de
$\dnrig{,s}(V)$, stables par $\partial_V$, alors $N_s^1=N_s^2$.
\end{corollary}

\begin{proof}
Le module $N_s^1+N_s^2$ v{\'e}rifie les m{\^e}mes hypoth{\`e}ses, et $N_s^1
\subset N_s^1+N_s^2$. Le lemme \ref{drdet} montre que
$\det(N_s^1)=\det(N_s^1+N_s^2)$, et donc que $N_s^1 = N_s^1+N_s^2  =
N_s^2$. 
\qed\end{proof}

\begin{proof}[du th{\'e}or{\`e}me \ref{ndr}]
Le lemme pr{\'e}c{\'e}dent montre que si $s_1 \leq s_2$, alors l'inclusion 
$\bnrig{,s_2}{,K}\otimes_{\bnrig{,s_1}{,K}}
N_{s_1} \subset N_{s_2}$ est en fait un isomorphisme. Posons
$\ndr(V)=\bnrig{}{,K} \otimes_{\bnrig{,s}{,K}} N_s$: on vient de voir que
cette d{\'e}finition ne d{\'e}pend pas du choix de $s$. De plus, $\phi(N_s)$
est un sous-module libre de rang $d$ de $\dnrig{,ps}(V)$, stable
par $\partial_V$, et par unicit{\'e} on a donc $\phi(N_s)=N_{ps}$. 
Par suite:
\[ \phi(\ndr(V))=\phi(\bnrig{}{,K} 
\otimes_{\bnrig{,s}{,K}} N_s) \subset \bnrig{}{,K} 
\otimes_{\bnrig{,ps}{,K}} N_{ps} = \ndr(V) \]
et $\ndr(V)$ est donc stable par $\phi$. Il est clair que $N_s$, et donc
$\ndr(V)$, est stable par $\Gamma_K$; par construction il est stable
par $\partial_V$. Enfin il n'est pas difficile de voir que $\ndr(V)$
d{\'e}termine $N_s$, qui est unique, et donc $\ndr(V)$ lui-m{\^e}me est unique.
\qed\end{proof}

Montrons maintenant la r{\'e}ciproque du th{\'e}or{\`e}me \ref{ndr}:
\begin{proposition}\label{recipndr}
Soit $V$ une repr{\'e}sentation $p$-adique de $G_K$.
On suppose qu'il existe $N_s \subset \dnrig{,s}(V)$,
libre de rang $d$ et stable par $\partial_V$.
Alors $V$ est de de Rham, et ses poids de Hodge-Tate sont n{\'e}gatifs.
\end{proposition}

\begin{proof}
Si $n$ est assez grand, alors $M_n=K_n[[t]] \otimes_{\iota_n(\bnrig{,s}{,K})}
\iota_n(N_s)$ est un sous 
$K_n[[t]]$-module de $\ddif^+(V)$, et
$M=K_{\infty}[[t]] \otimes_{K_n[[t]]} M_n$  
est un sous $K_{\infty}[[t]]$-module de $\ddif^+(V)$, 
de rang $\leq d$.
De plus, comme \[ \on{Frac}(\bnrig{,s}{,K}) \otimes_{\bnrig{,s}{,K}} N_s
= \on{Frac}(\bnrig{,s}{,K}) \otimes_{\bnrig{,s}{,K}} \dnrig{,s}(V), \]
$M$ est un  $K_{\infty}[[t]]$-module  libre de
rang $d$. Il est stable par $\partial_V$, soit $C$ la matrice de
$\partial_V$ dans une base de $M$. Alors l'{\'e}quation diff{\'e}rentielle 
$\frac{d}{dt}A+CA=0$ a une solution $A \in \on{GL}(d,K_{\infty}[[t]])$,
qui d{\'e}termine une famille libre $\{r_i\}$ de $M$, telle que
$\partial_V(r_i)=0$. La proposition
\ref{kernabla} montre que $V$ est de de Rham, et comme $\ddR(V)=\ker(\nabla_V)
\subset \ddif^+(V)$, c'est que ses poids de Hodge-Tate sont n{\'e}gatifs.
\qed\end{proof}

\subsection{Monodromie $p$-adique}
L'objet de ce paragraphe est de d\'emontrer le r\'esultat suivant, connu sous le nom de 
``conjecture de monodromie $p$-adique''. Cette conjecture, formul\'ee par Fontaine
\cite[6.2]{Bu88sst}, est l'analogue pour les repr\'esentations $p$-adiques du
th\'eor\`eme de monodromie $\ell$-adique de Grothendieck.
\begin{theorem}\label{monofont}
Si $V$ est une repr\'esentation $p$-adique de de Rham, alors $V$ est potentiellement semi-stable.
\end{theorem}

Les deux ingr\'edients principaux pour la d\'emonstration de ce th\'eor\`eme sont d'une part la
construction de $\ndr(V)$, et d'autre part un profond th\'eor\`eme d'Y. Andr\'e sur la structure
des \'equations diff\'erentielles $p$-adiques. Avant de d\'emontrer le th\'eor\`eme
\ref{monofont}, faisons quelques brefs rappels sur les \'equations diff\'erentielles $p$-adiques,
afin de fixer les notations et d'\'enoncer les r\'esultats de \cite{An01}.

Supposons dans la suite que $k$ est alg\'ebriquement clos.
Soit $\mathcal{R}_F=\bnrig{}{,F}$\index{RobbaF@$\mathcal{R}_F$}
l'anneau de Robba \`a coefficients dans $F$. Cet anneau est muni
d'un op\'erateur de Frobenius $\phi$ et d'une d\'erivation $\partial$, qui proviennent tous les
deux d'op\'erateurs d\'efinis sur $\adag{}_F$. De plus,  $\adag{}_F$ est un anneau de s\'eries
formelles en $\pi$, 
et $\phi(\pi)=\pi^p \mod{p}$. Parmi les extensions finies de l'an\-neau de Robba,
on s'int\'eresse
\`a celles qui proviennent d'une extension finie s\'eparable de $\bdag{}_F$, c'est-\`a-dire \`a
celles de la forme $\mathcal{R}_K=\bnrig{}{,K}$. 
Elles sont canoniquement munies du Frobenius de $\bnrig{}{,K}$,
et de la d\'erivation qui \'etend $\partial$.
Via la th\'eorie du corps de normes, 
on voit que ces extensions sont
celles qui proviennent d'une extension finie s\'eparable de $k((\pi))$.

Un $\phi$-module $M$ sur $\mathcal{R}_K$ est un
$\mathcal{R}_K$-module libre de rang fini, muni d'une application $\phi: M \ra M$ semi-lin\'eaire,
telle que $\phi^*(M)=M$ (c'est-\`a-dire que le $\mathcal{R}_K$-module engendr\'e par $\phi(M)$ est
\'egal \`a $M$). 
Une \'equation diff\'erentielle $p$-adique est un $\mathcal{R}_K$-module de
pr\'esentation finie muni d'une connexion $\nabla_M : M \ra M \otimes \Omega^1_{\mathcal{R}_K}$,
ou, ce qui revient au m\^eme car $\Omega^1_{\mathcal{R}_K} \simeq \mathcal{R}_K
\frac{d\pi_K}{\pi_K}$, d'une application $\partial_M: M \ra M$, qui \'etend la d\'erivation
$\partial:\mathcal{R}_K \ra
\mathcal{R}_K$ donn\'ee sur $\mathcal{R}_F$ par $\partial f(\pi)=(1+\pi)df/d\pi$, 
et qui satisfait la r\`egle de Leibniz. Une \'equation diff\'erentielle $p$-adique
$M$ est dite \^etre munie d'une structure de Frobenius, si $M$ est un $\phi$-module et si
le diagramme suivant est commutatif:
\[ \begin{CD}
M @>{\nabla_M}>> M \otimes \Omega^1_{\mathcal{R}_K} \\
@V{\phi}VV  @VV{\phi \otimes d\phi}V \\
M @>{\nabla_M}>> M \otimes \Omega^1_{\mathcal{R}_K}
\end{CD} \]
\'Etant donn\'e que sur $\mathcal{R}_F$, on a
\[ dx=\frac{\pi\partial x}{1+\pi}\frac{d\pi}{\pi}, \] 
cela revient \`a demander que $\partial_M \circ \phi = p \phi \circ \partial_M$.

Les \'equations diff\'erentielles $p$-adiques, munies ou non d'une structure de Frobenius, ont
\'et\'e \'etudi\'ees en d\'etail par Christol et Mebkhout; nous renvoyons \`a \cite{CM00} pour un
survol des r\'esultats et 
de quelques applications. En utilisant ces r\'esultats, Andr\'e a r\'eussi \`a
d\'emontrer \cite{An01} le th\'eor\`eme \ref{monoandre} ci-dessous, 
lui aussi connu sous le nom de ``conjecture de monodromie
$p$-adique''. Signalons que Kedlaya a plus r\'ecemment obtenu une d\'emonstration de ce m\^eme
r\'esultat dans \cite{KK00}, par des m\'ethodes tr\`es diff\'erentes (plus proches de l'esprit
du pr\'esent article). Enfin Mebkhout a aussi anonc\'e une d\'emonstration \cite{Mk01}.
\begin{theorem}\label{monoandre}
Tout module diff\'erentiel $M$ de pr\'esentation finie sur $\mathcal{R}_K$ admettant une structure
de Frobenius poss\`ede une base de solutions dans $\mathcal{R}'[\log \pi]$ o\`u $\mathcal{R}'$ est
l'extension finie \'etale de $\mathcal{R}_K$ issue d'une extension finie s\'eparable convenable de
$k_K((\pi_K))$, c'est-\`a-dire que $M$ poss\`ede une base de 
solutions dans $\mathcal{R}_L[\log\pi]$, pour une extension finie $L$ de $K$.
\end{theorem}

On est maintenant en mesure de d\'emontrer le th\'eor\`eme \ref{monofont}.
\begin{proof}[du th\'eor\`eme \ref{monofont}]
Soit $V$ une repr\'esentation $p$-adique de $G_K$ qui est de de Rham. On cherche \`a montrer que
$V$ est potentiellement semi-stable, et on peut supposer que
\begin{enumerate}
\item les poids de Hodge-Tate de $V$ sont n\'egatifs, en effet $V$ est potentiellement semi-stable
si et seulement si $V(-r)$ l'est pour $r \in \ZZ$;
\item $k$ est alg\'ebriquement clos, car $V$ est de Rham ou potentiellement
semi-stable si et seulement si sa restriction \`a l'inertie l'est (voir 
\cite[5.1.5]{Bu88sst} pour une
d\'emonstration de cela dans le cas semi-stable. Dans les autres cas, la d\'emonstration est
compl\`etement identique).
\end{enumerate}
On peut alors associer \`a $V$ le module $\ndr(V)$, qui est une \'equation diff\'e\-rentielle 
$p$-adique munie d'une structure de Frobenius. Le seul point \`a v\'erifier est que
$\phi^*(\ndr(V)) =\ndr(V)$. Mais $\phi^*(\ndr(V))$ est un sous-module libre de rang $d$ de
$\dnrig{}(V)$ qui est stable par $\partial_V$: par le th\'eor\`eme \ref{ndr}, un tel objet est
n\'ecessairement \'egal \`a $\ndr(V)$.

Par le th\'eor\`eme d'Andr\'e, et la remarque qui
le suit, il existe une extension finie $L/K$ telle que $(\ndr(V) \otimes_{\bnrig{}{,K}}
\bnst{}{,L})^{\partial_V = 0}$ est un $L \cap F^{nr}$-espace vectoriel de dimension $d$. Comme
l'alg\`ebre de Lie de $\Gamma_L$ agit trivialement sur celui-ci, il existe un sous groupe ouvert
de $\Gamma_L$ qui agit trvialement dessus \cite[chap V]{Se64}, 
et donc quitte \`a remplacer $L$ par
$L(\eps^{(n)})$, avec $n\gg 0$, on peut supposer que $\Gamma_L$ agit trivialement sur 
$(\ndr(V) \otimes_{\bnrig{}{,K}}
\bnst{}{,L})^{\partial_V = 0}$. Comme $\ndr(V) \otimes_{\bnrig{}{,K}} \bnst{}{,L} \subset
\dnst{}(V_L)$, o\`u $V_L$ est $V$ vue comme repr\'esentation de $G_L$, on en d\'eduit que 
$\dnst{}(V_L)^{\Gamma_L}$ est un  $L \cap F^{nr}$-espace vectoriel de dimension $d$. Par le
th\'eor\`eme \ref{isomcomp}, $V_L$ est semi-stable. Ceci implique que $V$ est
potentiellement semi-stable.
\qed\end{proof}

\subsection{Repr{\'e}sentations $\Cp$-admissibles}
Pour illustrer les m\'ethodes que nous avons d\'evelopp\'ees, nous donnons 
dans ce paragraphe une nouvelle
d\'emonstration d'un th\'eor\`eme de Sen sur les repr{\'e}\-sentations $\Cp$-admissibles: ce sont
celles dont la restriction \`a l'inertie est potentiellement triviale.

Tout d'abord, la th{\'e}orie de Sen permet de caract{\'e}riser facilement les
repr{\'e}sentations $\Cp$-admissibles en termes de l'op\'erateur $\Theta_V$:  
\begin{proposition}
Si $V$ est une repr{\'e}sentation $p$-adique, alors $V$ est
$\Cp$-admissible si et seulement si $\Theta_V=0$.
\end{proposition}
\begin{proof}
Si $V$ est $\Cp$-admissible, alors on a $(\Cp \otimes_{\Qp}
V)^{H_K}=\hat{K}_{\infty} \otimes_K (\Cp \otimes_{\Qp}
V)^{G_K}$, il contient le $K_{\infty}$-espace vectoriel 
$K_{\infty} \otimes_K (\Cp \otimes_{\Qp}
V)^{G_K}$, qui est stable par $\Gamma_K$ et de dimension $d$: c'est
donc $\dsen(V)$. L'action de $\nabla_V$ est bien triviale sur
cet espace.

R{\'e}ciproquement, si $\dsen(V)$ est muni de la connexion triviale, alors 
l'action de $\Gamma_K$ est discr{\`e}te sur $(\Cp \otimes_{\Qp}
V)^{G_K}$, et par Hilbert 90,
$V$ est $\Cp$-admissible.
\qed\end{proof}

La proposition suivante est due {\`a} 
Sen \cite{Sn73} et est un cas particulier de
r{\'e}sultats assez g{\'e}n{\'e}raux sur la caract{\'e}risation de l'alg{\`e}bre
de Lie de l'image de $I_K$ pour une repr{\'e}sentation $p$-adique en
terme de l'op{\'e}rateur $\Theta_V$:
\begin{proposition}\label{sen}
Soit $V$ une repr{\'e}sentation $p$-adique $\Cp$-admissible. Alors
la restriction de $V$ {\`a} $I_K$ est potentiellement triviale.
\end{proposition}

Le but de ce paragraphe est de donner une d{\'e}monstration de cette
proposition qui repose sur la th{\'e}orie des {\'e}quations
diff{\'e}rentielles $p$-adiques. Comme on ne s'int{\'e}resse qu'{\`a} la
restriction de $V$ {\`a} $I_K$, on suppose dans tout ce paragraphe que $k$ est
alg{\'e}briquement clos. Soient 
\[ \delta_{\pi_K}=\frac{\pi_K}{\partial\pi_K} \partial
\ \text{ et }\ 
\mu=p
\phi\left(\frac{\partial\pi_K}{\pi_K}\right)\frac{\pi_K}{\partial\pi_K},
\]\index{delatpiK@$\delta_{\pi_K}$}
ce qui fait que $\delta_{\pi_K}(f(\pi_K))=g(\pi_K)$ o{\`u} $g(X)=X\cdot df/dX$. 
Nous aurons besoin du r{\'e}sultat suivant de Tsuzuki \cite[5.1.1]{TS99}:
\begin{theorem}\label{tsutheo}
Si $k$ est alg\'ebriquement clos, et si on se donne 
deux matrices $A \in \on{GL}(d,\bdag{}_K)$ et
$C \in \on{M}(d,\bdag{}_K)$, telles que $\delta_{\pi_K}A+AC=\mu \phi(C)A$, telles que le
$(\phi,\partial)$-module correspondant $M$ est \'etale, 
alors il existe une extension finie $L/K$,
et une matrice $Y\in\on{GL}(d,\bdag{}_L)$, telle que $\delta_{\pi_K} Y+YC=0$ et
$Y=\phi(Y)A$. 

En d'autres termes, les $F$-isocristaux de pente nulle sont potentiellement triviaux.
\end{theorem}

On utilise ce th\'eor\`eme pour montrer le r\'esultat suivant:

\begin{lemma}\label{tsutsu}
Soit $\partial: \bdag{}_K \ra \bdag{}_K$ d{\'e}fini par
$\partial(x)=(1+\pi)
\frac{dx}{d\pi}$.
Si $V$ est une repr{\'e}sentation $p$-adique de $G_K$ et
$\partial_V:\ddag{}(V)\ra\ddag{}(V)$ est un op{\'e}rateur diff{\'e}rentiel au-dessus
de $(\bdag{}_K,\partial)$ tel que $\partial_V\circ\phi=p \phi\circ\partial_V$.
Alors le sous-groupe d'inertie de $H_K$ agit {\`a} travers un quotient
fini sur $V$.
\end{lemma}

\begin{proof}
Soient $A$ la matrice de $\phi$ et $C=\frac{\pi_K}{\partial(\pi_K)} \on{Mat}(\partial_V)$. La
relation de commutation entre $\phi$ et $\partial_V$ se traduit par $\delta_{\pi_K}A+AC=\mu
\phi(C)A$. On peut donc appliquer le th\'eor\`eme de Tsuzuki, et en d\'eduire qu'il existe $L/K$
finie telle que $\bdag{}_L\otimes_{\bdag{}_K} \ddag{}(V)$ est trivial: la matrice $Y$ nous
donne une base $\{y_i\}_i$ de $\bdag{}_L\otimes_{\bdag{}_K} \ddag{}(V)$,
telle que $\phi(y_i)=y_i$ ce qui fait que la restriction
de $V$ {\`a} $H_L$ est triviale et donc que la restriction de $V$
(comme repr{\'e}sentation de $G_K$) {\`a} l'inertie de $H_K$ est 
potentiellement triviale.
\qed\end{proof}

\begin{lemma}\label{controle}
Si la connexion $\nabla_V$ de $\dnrig{}(V)$ v{\'e}rifie $\nabla_V(\dnrig{}(V))
\subset t \dnrig{}(V)$, alors l'op{\'e}\-ra\-teur $\partial_V=t^{-1} \nabla_V$
v{\'e}rifie $\partial_V(\ddag{}(V)) \subset \ddag{}(V)$.
\end{lemma}

\begin{proof}
Soit $e_1,\hdots,e_d$ une base de $\ddag{}(V)$ sur $\bdag{}_K$, telle
que la matrice de $\phi$ v{\'e}rifie $Q=\on{Mat}(\phi) \in \on{GL}_d(\adag{}_K)$.
C'est aussi une base de $\dnrig{}(V)$ sur $\bnrig{}{,K}$. Soit 
$D=\on{Mat}(\partial_V) \in \on{M}_d(\bnrig{}{,K})$.
On remarquera que $\partial Q 
\in \on{M}_d(\bdag{}_K)$
puisque $\partial$ pr{\'e}serve $\bdag{}_K$. 
Le fait
que $\nabla_V$ et $\phi$ commutent montre que  $p\phi(D)Q=\partial Q+QD$.
Soit $H=D+Q^{-1}\partial Q$. L'{\'e}quation pr{\'e}cedente s'{\'e}crit alors
\[ H-pQ^{-1}\phi(H)Q
= - p Q^{-1}\phi(Q^{-1}\partial(Q))Q \]
et $- p Q^{-1}\phi(Q^{-1}\partial(Q))Q \in \on{M}_d(\bnrig{}{,K})$.
Le r{\'e}sultat que l'on 
cherche {\`a} {\'e}tablir est que $H$ et donc $D$ a ses coefficients dans
$\bdag{}_K$. Il suffit donc de montrer le fait suivant:
si $H \in \on{M}_d(\bnrig{}{,K})$ v{\'e}rifie $H-pQ^{-1}\phi(H)Q = R \in
\on{M}_d(\bdag{}_K)$, 
alors en fait $H \in \on{M}_d(\bdag{}_K)$, ce que nous faisons
maintenant. 
Soient $|\cdot|_I$ les
normes correspondant aux valuations $V_I$; on les {\'e}tend {\`a}
$\on{M}_d(\bdag{}_K)$ en d{\'e}cidant que $|M|_I=\sup |m_{ij}|_I$.
Un {\'e}l{\'e}ment $x$ de $\bnrig{,r}{,K}$
est dans $\bdag{,r}_K$ si et seulement si la suite des $|x|_I$ 
(avec $I \subset [r;+\infty[$) est born{\'e}e, et il en est donc de m{\^e}me
pour les matrices {\`a} coefficients dans ces anneaux.
Comme $Q \in \on{GL}_d(\adag{}_K)$, on a pour toute matrice $M$,
$|Q^{-1}M|_I=|MQ|_I=|M|_I$.
Fixons $r$ tel que $Q$ et $Q^{-1} \in \on{GL}_d(\adag{,r}{,K})$, et 
$R \in \on{M}_d(\bdag{,r}{,K})$.
Soit $I_h=[p^h r; p^{h+1}r]$.
Dans notre cas $H-pQ^{-1}\phi(H)Q=R$ et il existe une constante $C$
telle que pour tout $I$, on ait: $|R|_I \leq C$.
Comme on a $H=R+p Q^{-1}\phi(H)Q$, alors pour $h \geq 1$:
\[ |H|_{I_h} \leq \sup( C, p^{-1}|H|_{I_{h-1}} ) \leq \cdots \leq 
\sup( C, p^{-h}|H|_{I_0} ) \]
ce qui fait que $|H|_{I_h} \leq \sup ( C, |H|_{I_0} )$, pour tout $h$,
et donc que $H \in  \on{M}_d(\bdag{}_K)$.
\qed\end{proof}

\begin{lemma}\label{conntriv} 
La connexion $\Theta_V$ sur
$\dsen(V)$ est triviale si et seulement
si $\nabla_V(\ddag{}(V))
\subset t\cdot \dnrig{}(V)$.
\end{lemma}

\begin{proof}
Si $\nabla_V(\ddag{}(V))
\subset t\cdot \dnrig{}(V)$, alors comme $(\dsen(V),\Theta_V)$ est
l'image par $\theta\circ\iota_n$ de $(\dnrig{,r_n}(V),\nabla_V)$,
c'est que $\Theta_V=0$.

R{\'e}ciproquement si $\Theta_V=0$, alors 
pour tout $n \gg 0$ on a 
\begin{multline*}
\nabla_V(\ddag{,r_n}(V))\subset 
\ker(\theta\circ\iota_n:\dnrig{,r_n}(V) \ra \Cp \otimes_{\Qp} V) \\
= \phi^{n-1}(q)
\dnrig{,r_n}(V) \end{multline*} 
par le lemme \ref{kerthetadhol},
et le lemme \ref{divtv} permet de conclure.
\qed\end{proof}

\begin{proof}[de la proposition \ref{sen}]
Comme on s'int{\'e}resse {\`a} la restriction de $V$ {\`a} l'iner\-tie,
on peut supposer que le corps r{\'e}siduel $k$ de $K$
est alg{\'e}briquement clos, ce que l'on fait maintenant.
On a vu
qu'une repr{\'e}sen\-ta\-tion $V$ est $\Cp$-admissible si et
seulement si le module $\dsen(V)$ qui lui est associ{\'e} par la th{\'e}orie
de Sen est muni de la connexion triviale et 
on a montr{\'e} que cette connexion est l'image
par $\theta$ de celle qui existe
sur $\ddif^+(V)$. Le lemme \ref{conntriv} montre que dans ce cas 
$\nabla_V(\ddag{}(V))
\subset t\cdot \dnrig{}(V)$.

Le lemme \ref{controle} montre
alors que $t^{-1} \nabla_V$ est une connexion surconvergente sur $\ddag{}(V)$.
Le lemme \ref{tsutsu} montre ensuite qu'il existe une extension finie $L/K$ telle que
l'action de $H_L$ sur $V$ est
triviale sur $V$. Comme 
$\nabla_V=t \partial_V= 0$, l'action de $\Gamma_L$
est finie, ce qui fait qu'il existe $n \geq 0$ tel que $G_{L_n}$ agit
trivialement sur $V$, et donc 
que la restriction 
de $V$ {\`a} $I_K$ est potentiellement triviale.
\qed\end{proof}

\begin{remark}
Une repr\'esentation $\Cp$-admissible $V$ est de de Rham, et on a $\ndr(V)=\dnrig{}(V)$ dans ce
cas particulier.
\end{remark}

\begin{remark}\label{remasentsu}
En fait, le th\'eor\`eme de Sen implique le th\'eor\`eme 
de Tsuzuki (et lui est donc \'equivalent). Indiquons bri\`evement de quoi il s'agit:
si $N$ est un isocristal de pente nulle, on peut ``exponentier'' la connexion pour
construire une matrice que l'on utilise pour faire 
agir un sous-groupe ouvert de $\Gamma_K$. Le fait que $\phi$ est de pente nulle assure que cette
matrice est surconvergente. On a donc construit un $(\phi,\Gamma_K)$-module \'etale,
c'est-\`a-dire une repr\'esentation $p$-adique $V$. Il n'est pas dur de montrer que $V$ est
$\Cp$-admissible, et d'en conclure, via le th\'eor\`eme de Sen, que $V$, et donc $N$, devient
trivial apr\`es une extension finie.
\end{remark}

\section{Extensions de repr{\'e}sentations semi-stables}
Dans ce chapitre, on utilise
le th\'eor\`eme de monodromie $p$-adique, pour montrer que si $E$ est une
repr\'esentation de de Rham qui est une extension de deux repr\'esentations
semi-stables, alors $E$ est semi-stable. On va montrer que:
\begin{enumerate}
\item Toute repr{\'e}sentation ordinaire de $G_K$ est semi-stable;
\item une extension de deux repr{\'e}sentations semi-stables qui 
est de de Rham est elle-m{\^e}me semi-stable;
\item si $V$ est une repr\'esentation semi-stable dont les 
poids de Hodge-Tate sont tous $\geq 2$, alors 
$\exp_V : \ddR(V) \ra H^1(K,V)$ est un isomorphisme.
\end{enumerate}
Ces trois r{\'e}sultats {\'e}taient connus dans le cas d'un corps
r{\'e}siduel fini, o{\`u} on peut les d{\'e}duire de calculs de dimension de
groupes de cohomologie galoisienne (ce qui n'est plus possible si le
corps r{\'e}siduel n'est pas fini).
Le (1) avait {\'e}t{\'e} d{\'e}montr{\'e} dans ce cas l{\`a} par
Perrin-Riou \cite{Bu88ord,BP94,BP99}  
comme corollaire 
des calculs de Bloch et Kato 
\cite{BK91}, le (2) par Hyodo
\cite{Hy95,Ne95}, et le (3) par
Bloch et Kato.

\begin{remark}
Dans une version ant\'erieure de cet article, 
on utilisait le calcul de la cohomologie galoisienne via les
$(\phi,\Gamma_K)$-modules pour d{\'e}montrer ces r\'esultats.
\end{remark}

\subsection{Extensions de repr\'esentations semi-stables}
L'objet de ce paragraphe est de montrer le th\'eor\`eme suivant:
\begin{theorem}\label{extsst}
Si $E$ est une repr\'esentation de de Rham, qui est une extension de 
$W$ par $V$, o\`u $V$ et $W$ sont semi-stables, alors $E$ est semi-stable.
\end{theorem}
et son corollaire:
\begin{corollary}\label{ordisst}
Toute repr{\'e}sentation ordinaire de $G_K$ est semi-stable.
\end{corollary}

Commen\c{c}ons par montrer la proposition \ref{extsst};
par le th\'eor\`eme \ref{monofont}, il suffit de montrer le 
lemme\footnote{Je remercie le referee de m'avoir signal\'e ce lemme et sa d\'emonstration.} 
suivant:
\begin{lemma}
Si $E$ est une repr\'esentation potentiellement semi-stable de $G_K$, 
qui est une extension de $W$ par $V$, o\`u $V$ et $W$ sont 
deux repr\'esentations semi-stables de $G_K$, alors $E$ est semi-stable.
\end{lemma}

\begin{proof}
Soit $L$ une extension finie de $K$, telle que la restriction de $E$ \`a $G_L$ est semi-stable.
Soit $\dst^L(\cdot)=(\bst \otimes_{\Qp} \cdot)^{G_L}$. On a une suite exacte:
\[ 0 \ra \dst^L(V) \ra \dst^L(E) \ra \dst^L(W) \ra 0, \]
et en prenant les invariants par $\on{Gal}(L/K)$, on en d\'eduit:
\[ 0 \ra \dst(V) \ra \dst(E) \ra \dst(W) \ra H^1(\on{Gal}(L/K),\dst^L(V)). \]
Pour des raisons de dimension, il suffit de montrer que
\[ H^1(\on{Gal}(L/K),\dst^L(V))=0, \] ce qui r\'esulte du fait que 
$\dst^L(V)$ est un $\Qp$-espace vectoriel, et donc 
en particulier que $\#\on{Gal}(L/K)$ y est inversible.
\qed\end{proof}

Pour montrer le corollaire \ref{ordisst}, on va tout 
d'abord montrer le lemme suivant; si $V$ est
une repr\'esentation de de Rham, $H^1_g(K,V)$ 
est l'ensemble des classes d'extensions qui sont de de Rham.
\begin{lemma}\label{posdr}
Si $V$ est une repr\'esentation de de Rham dont les 
poids de Hodge-Tate sont tous $\geq 1$, alors
$H^1(K,V)=H^1_g(K,V)$.
\end{lemma}

\begin{proof}
Soit $E$ une extension de $\Qp$ par $V$: il suffit de montrer que $E$ est de de Rham.
L'hypoth\`ese sur les poids de Hodge-Tate de $V$ implique que $H^1(K,\bdR^+ \otimes_{\Qp} V)=0$
(ceci suit de la d\'ecomposition de Hodge-Tate de $V \otimes_{\Qp}(\oplus \Cp(i))$, du
th\'eor\`eme de Tate: $H^1(K,\Cp(j))=0$ si $j \neq 0$, et des suites exactes $0 \ra
t^{i+1}\bdR^+ \ra t^i \bdR^+ \ra \Cp(i) \ra 0$, que l'on utilise pour $i \geq 0$).  

On a donc un morceau de suite exacte
$(\bdR^+ \otimes_{\Qp} E)^{G_K} \ra (\bdR^+)^{G_K} \ra H^1(K,\bdR^+ \otimes_{\Qp} V) = 0$, ce qui
fait que $1 \in (\bdR^+)^{G_K}=K$ se rel\`eve en un \'el\'ement de
$(\bdR^+ \otimes_{\Qp} E)^{G_K}$, qui n'est pas dans $(\bdR^+ \otimes_{\Qp} V)^{G_K}$. On en
d\'eduit imm\'ediatement que la dimension de $(\bdR^+ \otimes_{\Qp} E)^{G_K}$ est la bonne, et
donc que $E$ est de de Rham. 
\qed\end{proof}

\begin{proof}[du corollaire \ref{ordisst}]
Rappelons \cite[1.1]{Bu88ord} qu'une repr\'esentation $V$ est ordinaire, s'il existe une
filtration $\on{Fil}^i V$ de $V$, d\'ecroissante 
exhaustive et s\'epar\'ee, par des sous-espaces 
$\on{Fil}^i V$ stables par $G_K$, et telle que le groupe d'inertie $I_K$ agit sur
$\on{Fil}^i V / \on{Fil}^{i+1} V$ par $\chi^i$.

Une r\'ecurrence imm\'ediate sur la longueur de la filtration montre que
qu'une repr\'esentation ordinaire est de de Rham, et comme elle est extension successive 
de repr\'esentations semi-stables (par hypoth\`ese de r\'ecurrence), elle est elle-m\^eme
semi-stable. 
\qed\end{proof}

\subsection{L'exponentielle de Bloch-Kato}
Dans ce paragraphe, $V$ est une repr{\'e}sentation semi-stable. 
Rappelons que les anneaux $\bmax$ et $\bdR$ sont reli{\'e}s par la suite
exacte fondamentale (cf \cite[III.5]{Co98} et \cite{Bu88per}):
\[ 0 \ra \Qp \ra \bmax^{\phi=1} \ra \bdR/\bdR^+ \ra 0 \]
En tensorisant avec $V$ et en prenant les invariants par $G_K$ on
obtient un d{\'e}but de suite exacte longue:
\begin{multline*} 
0 \ra V^{G_K} \ra \dcris(V)^{\phi=1} \ra 
((\bdR/\bdR^+)\otimes_{\Qp} V)^{G_K} \\ 
\ra H^1(K,V) \ra
H^1(K,\bmax^{\phi=1}\otimes_{\Qp} V) 
\end{multline*}
Soit $H^1_e(K,V)=\ker(H^1(K,V) \ra H^1(K,\bmax^{\phi=1}\otimes_{\Qp}
V))$. On d{\'e}duit de la suite exacte ci-dessus une application 
de $\ddR(V)$ dans $H^1_e(K,V)$ appel{\'e}e exponentielle de Bloch-Kato,
et not{\'e}e $\exp_V$\index{expV@$\exp_V$}. 
D'autre part comme $V$ est de de Rham on a
\[ ((\bdR/\bdR^+)\otimes_{\Qp} V)^{G_K}=\ddR(V)/ \on{Fil}^0 \ddR(V) \]
et \cite[lemma 3.8.1]{BK91} l'image de $\exp_V$ est $H^1_e(K,V)$ tout
entier. Dans le cas o{\`u} $k$ est fini des calculs de dimension
montrent que si $r\gg 0$, alors $H^1_e(K,V(r))=H^1(K,V(r))$ et donc
que $\exp_{V(r)}: \ddR(V(r)) \ra H^1(K,V(r))$ est un isomorphisme
(on remarquera que
si $r\gg 0$, alors $\on{Fil}^0 \ddR(V)=\{0\}$). 
Le but de ce paragraphe est de montrer, sans condition sur $k$, que
si tous les poids de Hodge-Tate de $V$ sont $\geq 2$, 
alors $\exp_V: \ddR(V) \ra H^1(K,V)$ 
est un isomorphisme.

\begin{lemma}
Si $V$ est une repr{\'e}sentation semi-stable telle que  
$1-\phi: \dst(V) \ra \dst(V)$ est surjectif
et $\dst(V)^{\phi=1/p}=0$, alors
$H^1_g(K,V)=H^1_e(K,V)$.
\end{lemma}

\begin{proof}
Soit $c \in H^1_g(K,V)$ et $E$ l'extension de $\Qp$ par $V$ qui 
correspond {\`a} $c$. Le th{\'e}or{\`e}me \ref{extsst} montre qu'il existe
$x\in \bst\otimes_{\Qp} V$ tel que $c(g)=(g-1)x$. Comme $c(g)\in V$ on a 
$(\phi-1)c(g)=0$ et donc $\phi(x)-x \in (\bst\otimes_{\Qp} V)^{G_K}$. Il existe
donc $y\in \dst(V)$ tel que $\phi(y)-y=\phi(x)-x$ ce qui revient {\`a}
$\phi(x-y)=x-y$ et 
quitte {\`a} remplacer $x$ par $x-y$
on peut donc supposer que
$x \in \bst^{\phi=1}\otimes_{\Qp} V$. Comme $N$ commute {\`a} $G_K$ on a
$N(c(g))=0=(g-1)(N(x))$ et donc $N(x)\in \dst(V)$. De plus comme
$x \in \bst^{\phi=1}\otimes_{\Qp} V$ on doit avoir $N(x) 
\in \dst(V)^{\phi=1/p}=0$ et donc $x\in\bmax^{\phi=1}\otimes_{\Qp} V$
ce qui montre que $c(g)\in H^1_e(K,V)$.
\qed\end{proof}

\begin{lemma}
Si $V$ a tous ses poids de Hodge-Tate $\geq 2$, alors
$V$ satisfait 
les conditions du lemme pr{\'e}c{\'e}dent.
\end{lemma}

\begin{proof}
Rappelons que le Frobenius $\phi$ sur
$\dst(V)$ est bijectif. De plus, si $W$ est une repr\'esentation 
semi-stable dont les poids de Hodge-Tate
sont positifs, alors il existe un r\'eseau $M_W$ de 
$\dst(W)$ tel que $M_W \subset \phi(M_W)$ (c'est le dual 
de l'assertion que si $X$ est semi-stable positive, il existe un
r\'eseau de $\dst(X)$ stable par $\phi$).
En tordant, on voit donc qu'il existe 
un r{\'e}seau $M$ de $\dst(V)$,
tel que $\phi^{-1}(M) \subset p^2 M$.
Il est alors {\'e}vident
que $1-\phi=-\phi(1-\phi^{-1})$ est
surjectif, et que $\dst(V)^{\phi=1/p}=0$.
\qed\end{proof}

\begin{theorem}\label{expbk}
Si $V$ est une repr{\'e}sentation semi-stable de $G_K$, dont les poids de Hodge-Tate
sont $\geq 2$, alors l'exponentielle de Bloch-Kato $\exp_V : \ddR(V) \ra
H^1(K,V)$ est un isomorphisme.
\end{theorem}

\begin{proof}
Les deux lemmes pr{\'e}c{\'e}dents, et le lemme \ref{posdr}, montrent que l'on a 
\[ H^1(K,V)=H^1_g(K,V)=H^1_e(K,V) \]
et $\exp_V$ est donc surjective. D'autre
part, on a $\on{Fil}^0 \ddR(V)=\{0\}$, 
ce qui montre que $\exp_V$ est injective.
C'est donc un isomorphisme.
\qed\end{proof}

\newpage
\section*{Diagramme des anneaux de p{\'e}riodes}\label{diag}
Le diagramme ci-dessous r{\'e}capitule les relations entre les
diff{\'e}rents anneaux de p{\'e}riodes. 
Les fl{\`e}ches qui se terminent par
$\xymatrix@1{&\ar@{->>}[r]&}$ sont surjectives, la fl{\`e}che 
en pointill{\'e}s $\xymatrix@1{&\ar@{.>}[r]&}$ est une limite inductive de
morphismes d{\'e}finis sur des sous-anneaux
($\iota_n: \btst{,r_n}{} \ra \bdR^+$), 
et toutes les autres sont injectives.
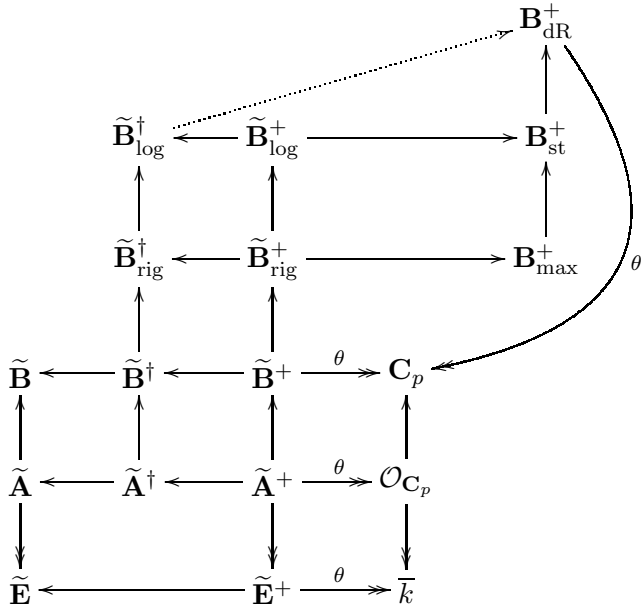
\[ \xymatrix{
& &  & & \bdR^+ \ar@{->>}@/^5pc/[lddd]^{\theta} & \\
& \btst{}{} \ar@{.>}[urrr] & \btstplus{} \ar[rr] 
\ar[l] & & \bst^+ \ar[u]  & \\
& \btrig{}{} \ar[u] & \btrigplus{} \ar[u] 
\ar[l] \ar[rr] & & \bmax^+ \ar[u]  & \\  
\bt & \btdag{} \ar[l] \ar[u] & \btplus \ar[l]
\ar[u] \ar@{->>}[r]^{\theta} & \Cp & & \\ 
\at \ar[u] \ar@{->>}[d] & \atdag{} \ar[l] \ar[u] 
& \atplus \ar[l] \ar[u] 
\ar@{->>}[d] \ar@{->>}[r]^{\theta} & 
\OO_{\Cp} \ar[u] \ar@{->>}[d] & & \\ 
\et & & \etplus \ar[ll] \ar@{->>}[r]^{\theta} &  \overline{k}  & &  \\
} \]
Tous les anneaux qui ont des tildes ($\tilde{{\quad}}$) ont aussi des
versions sans tilde: on passe de cette derni{\`e}re {\`a} la version avec
tilde en rendant Frobenius inversible, et en compl{\'e}tant. Par
exemple, $\et$ est le compl{\'e}t{\'e} de la perfection de $\e$.

Les trois anneaux de la colonne de gauche 
(du moins leurs versions sans tildes)
sont reli{\'e}s 
{\`a} la th{\'e}orie des $(\phi,\Gamma_K)$-modules. Les trois anneaux de
la colonne de droite sont reli{\'e}s {\`a} la th{\'e}orie de Hodge
$p$-adique. Pour faire le lien entre ces deux th{\'e}ories, on passe
d'un c{\^o}t{\'e} {\`a} un autre via tous les anneaux interm{\'e}diaires. La situation
id{\'e}ale est quand on peut se placer dans la colonne du milieu
(par exemple, de haut en bas: les repr{\'e}sentations 
semi-stables, cristallines, ou
de hauteur finie).

R{\'e}capitulons les diff{\'e}rents anneaux de s{\'e}ries formelles qui
interviennent; soit $C[r;1[$ la couronne $\{ z \in \Cp,\ p^{-1/r} \leq
|z|_p < 1 \}$.
On a alors:
\begin{align*}
{\mathbf{E}_F^+} &= {k[[T]]} \\
{\mathbf{A}_F^+} &= {\OO_F[[T]]} \\
{\mathbf{B}_F^+} &= {F \otimes_{\OO_F} \OO_F[[T]]}  \\
{} & {} \\
{\mathbf{E}_F} &= {k((T))} \\
{\mathbf{A}_F} &= {\hat{\OO_F[[T]][T^{-1}]}} \\
{\mathbf{B}_F} &= {F \otimes_{\OO_F} \hat{\OO_F[[T]][T^{-1}]}}  \\
{} & {} \\
{\adag{,r}_F} &= { \{\text{s{\'e}ries de Laurent $f(T)$, qui convergent 
sur $C[r;1[$,} } \\ 
& {\text{et y sont born{\'e}es par $1$}\} } \\
{\bdag{,r}_F} &= {\{ \text{s{\'e}ries de Laurent $f(T)$, qui convergent 
sur $C[r;1[$,} } \\ 
& { \text{et y sont born{\'e}es}\} } \\
{} & {} \\
{\bnrig{,r}{,F}} &= { \{\text{s{\'e}ries de Laurent $f(T)$, qui convergent 
sur $C[r;1[$}\} } \\
{\bnst{,r}{,F}} &= {\bnrig{,r}{,F}[\log(T)]} \\
{} & {} \\
{\bhol{,F}} &= { \{ f(T) \in F[[T]],\ 
\text{$f(T)$ converge sur le disque ouvert $D[0;1[$} \} } \\
{\blog{,F}} &= {\bhol{,F}[\log(T)]}
\end{align*}

\renewcommand{\indexname}{Index des notations}
\printindex

\renewcommand{\refname}{Bibliographie}

\end{document}